\documentclass[a4paper,10pt]{amsart}
\usepackage{fullpage}
\usepackage[T1]{fontenc}
\usepackage{amsmath,amsfonts,amssymb,amsthm}
\usepackage{pb-diagram,graphicx}
\usepackage{epsfig}
\usepackage{graphics}
\usepackage[pst-3d]{pstricks}
\usepackage{psfrag}
\usepackage{latexsym}
\usepackage{multicol}
\usepackage{times}
\usepackage{bbold}
\usepackage{dsfont}
\usepackage{stmaryrd}
\usepackage{subcaption}
\usepackage{hyperref}
\usepackage{booktabs}
\usepackage{listings}
\usepackage{color}
\usepackage[hyphenbreaks]{breakurl}
\newcommand{\transp}{\mathsf{T}}
\newcommand{\hold}{\mathcal{D}}

\newcommand{\e}{\epsilon}
\newcommand{\Om}{\Omega}
\newcommand{\R}{\mathds{R}}

\newcommand{\divv}{\operatorname{div}}

\newcommand{\tr}{\operatorname{tr}}
\newcommand{\vp}{\varphi}
\newcommand{\VV}{\theta}
\newcommand{\ben}{\begin{equation}}
\newcommand{\een}{\end{equation}}
\newcommand{\bproof}{\begin{proof}}
\newcommand{\eproof}{\end{proof}}
\newcommand{\Sb}{S}

\newcommand{\dt}{{\frac{d}{dt}}}

\newcommand{\Tt}{T_t}

\newcommand{\f}{\lstinline}

\newcommand{\DV}{dJ^{\tt{vol}}(\Om;\VV)}
\newcommand{\DS}{dJ^{\tt{surf}}(\Om;\VV)}
\newcommand{\singset}{L}
\newcommand{\Gnum}{K}
\newcommand{\ddo}{:}
\newcommand{\fd}{F_\hold}
\newcommand{\fg}{F_\Gamma}

\newcommand{\ca}{c_1}
\newcommand{\cb}{c_2}
\newcommand{\cc}{c_3}

\newtheorem{definition}{Definition}

\newtheorem{proposition}{Proposition}

\definecolor{dkgreen}{rgb}{0,0.6,0}
\definecolor{gray}{rgb}{0.5,0.5,0.5}
\definecolor{mauve}{rgb}{0.58,0,0.82}

\lstnewenvironment{ttlisting}{
\lstset{
basicstyle=\ttfamily,
frame=none,
numbers=none,
}
}{}

\title{A level set-based structural optimization code using FEniCS}
\author{Antoine Laurain}
\address{Instituto de Matem\'atica e Estat\'istica,Universidade de S\~{a}o Paulo,
Rua do Mat\~{a}o, 1010, 05508-090 - S\~{a}o Paulo, Brazil}
\email{laurain@ime.usp.br}
\begin{document}
\begin{abstract}
This paper presents an educational code written using FEniCS, based on the level set method, to perform compliance minimization in structural optimization.
We use the concept of distributed shape derivative to compute a descent direction for the compliance, which is defined as a shape functional. 
The use of the  distributed shape derivative is facilitated by  FEniCS, which allows to handle complicated partial differential equations with a simple implementation.
The code is written for compliance minimization in the framework of linearized elasticity, and can be easily adapted to tackle other functionals and partial differential equations.
We also provide an extension of the code for compliant mechanisms.
We start by explaining how to compute shape derivatives, and discuss the differences between the  distributed and boundary expressions of the shape derivative. 
Then we describe the implementation in details, and  show the application of this code to some classical benchmarks of topology optimization.
The code is available at \burl{http://antoinelaurain.com/compliance.htm}, 
and the  main file is also given in the appendix.
\end{abstract}

\maketitle

\section{Introduction}
The popular ``99 line'' Matlab code by Sigmund published in 2001 \cite{Sigmund2014} has started a trend of sharing and publishing educational codes for structural optimization. 
Since then, an upgrade of the ``99 line'' code has been published, improving speed and reducing the code size to 88 lines; see \cite{Andreassen2010}. 
The codes of \cite{Andreassen2010,Sigmund2014} are written for Matlab and are based on the solid isotropic microstructure with penalty  (SIMP) approach \cite{bendsoe1,ZHOU1991309}.
Various other codes have been published using different approaches and/or other platforms than Matlab. We review here several categories of approaches to tackle this problem.

In the SIMP approach the material is allowed to have intermediate values, and the optimization variables are the material densities of the mesh elements. The intermediate values are also penalized using a power law to enforce $0-1$ values. Using filtering techniques, it provides feasible designs.
Considering SIMP approaches as in \cite{Sigmund2014}, Talischi et al. have introduced \texttt{PolyMesher} \cite{Talischi2012} and \texttt{PolyTop} \cite{Talischi2012b} to provide a MATLAB implementation of topology optimization using a general framework for finite element discretization and analysis. 

Another category of approaches for topology optimization which has emerged after the SIMP approach are level set methods. 
They consist in representing the boundary of the moving domain $\Om$ as the zero level set of a function $\phi$. 
Level set methods were introduced by Osher and Sethian \cite{MR965860} in the context of the mean curvature flow to facilitate the modelization of topological changes during curve evolution.
Since then, they have been applied to many shape optimization and boundary perturbations problems. 
There is already a substantial literature for level set methods applied to structural optimization, see \cite{MR1911658,MR2033390,MR1843648,MR1783559,MR1951408}
for the pioneering works using this approach,  and \cite{MR3107583} for a review.
Early references for level set approaches include a code in FEMLAB \cite{Liu2005} by Liu et al. in 2005, a Matlab code \cite{Challis2009} in the spirit of the ``99 line'' code, by Challis in 2010, and a 
a 88 lines Matlab code \cite{Otomori2014} using a reaction-diffusion equation by Otomori et al in 2014. 

Other approaches to structural topology optimization include phase-field methods \cite{MR2108636},
level-set methods without Hamilton-Jacobi equations \cite{NME:NME824},
and an algorithm based on the notion of topological derivative \cite{MR2235384}.   
We also mention an early FreeFem++ code \cite{MR2252743} by Allaire and Pantz in 2006, implementing the boundary variation method and the homogenization. For a critical comparison of four different level-set approaches and one phase-field approach, see \cite{Gain2013}. 

The code presented in the present paper enters the category of level set methods.
In the usual approach, the concept of shape derivative \cite{MR2731611,SokZol92}  is used to compute the sensitivity of the objective functional.
Is is known that the shape derivative is a distribution on the boundary  of the domain, and algorithms are usually based on this property. 
This means that the shape derivative is expressed as a boundary integral, and then extended to the entire domain or to a narrow band for use in the level set method; see
\cite{MR1911658,MR2033390,Challis2009,MR2356899,MR2459656,HLsegmentation,MR1951408} for applications of this approach.
The shape derivative can also be written as a domain integral, 
which is called {\it distributed}, {\it volumetric} or {\it domain expression} of the shape derivative; see \cite{MR2642680,delfour1997computation,MR3348199,LaurainSturmM2AN}, and \cite{MR3350625,MR3013681,MR3493474} for applications.

From a numerical point of view, the distributed expression is often easier to implement than the boundary expression as it is a volume integral.
Other advantages of  the distributed expression are presented in \cite{MR2642680,MR3348199}. 
In \cite{MR2642680}, it is shown that the discretization 
and the shape differentiation processes commute for the volume expression but not for the
boundary expression; i.e., a discretization of the boundary expression does not generally lead to the same expression as the shape derivative computed after the problem is discretized.
In \cite{MR3348199}, the authors conclude that ``volume based
expressions for the shape gradient often offer better accuracy than the use of formulas
involving traces on boundaries''. 
See also \cite{MR3264217} for a discussion about the difficulty to use the boundary expression in the multi-material setting.
In the present paper, the main focus is the compact yet efficient implementation of the level set method for structural optimization allowed by the distributed shape derivative. 
We also show that it is useful to handle the ersatz material approach.
Combining these techniques, we obtain a straightforward and general  way of solving the shape optimization problem, from the rigorous calculation of the shape derivative to the numerical implementation.

The choice of FEniCS for the implementation is motivated by its ability to facilitate the implementation of complicated variational formulations, thanks to a near-mathematical notation. 
This is appropriate in our case since the expression of the distributed shape derivative is usually lengthy. 
The FEniCS Project (\burl{https://fenicsproject.org/}) is a collaborative project with a particular focus on automated solution of differential equations by finite element methods; see \cite{ans20553,fenics:book}.

The paper is structured as follows. 
In Section \ref{section1}, we recall the definition of the shape derivative, and show the relation between its distributed and boundary expression.
In Section \ref{sec:shape_der_compliance}, we compute the shape derivative in distributed and boundary form for a general functional in linear elasticity using a Lagrangian approach, and we discuss the particular cases of compliance and compliant mechanisms. 
In Section \ref{sec:descent}, we show how to obtain descent directions.
In section \ref{section5}, we explain the level set method used in the present paper, which is a variation of the usual level set method suited for the distributed shape derivative. 
In this section we also describe the discretization and reinitialization procedures. 
In section \ref{sec:numerical_implementation}, we explain in details the numerical implementation.
In section \ref{sec:num_res}, we show numerical results for several classical benchmarks.
Finally, in section \ref{sec:init_time}, we discuss the computation time and the influence of the initialization on the optimal design.
In the appendix, we give the code for the main file \f{compliance.py}.
\section{Volume and boundary expressions of the shape derivative}\label{section1}
In this section we recall basic notions about the shape derivative, the main tool used in this paper.
Let $\mathcal{P}(\hold)$ be the set of subsets of $\hold$, where the so-called {\it universe} $\hold\subset \R^m$ is assumed to be a piecewise smooth open and bounded set, and $\mathds{P}$ be a subset of $\mathcal P(\hold)$.
In our numerical application, $\hold$ is a rectangle.
Let $k\geq 1$ be an integer, $C^k_c(\R^m,\R^m)$ be the set of $k$-times continuously differentiable vector-valued functions with compact support.
Let $\singset\subset\partial\hold$ be the set of points where the normal $n$ is not defined, i.e. the set of singular points of $\partial\hold$, such as the corners of a rectangle.
Define 
$$ \Theta^k(\hold) = \{\VV\in C^k_c(\R^m,\R^m)| \VV\cdot n_{|\partial\hold\setminus \singset} = 0 \mbox{ and } \VV_{|\singset} = 0\}$$
equipped with the topology induced by $C^k_c(\R^m,\R^m)$.
Consider a vector field $\VV\in \Theta^k(\hold)$ and the associated flow
$\Tt^{\VV}:\R^m\rightarrow \R^m$, $t\in [0,\tau]$ defined for each $x_0\in \R^m$ as $\Tt^{\VV}(x_0):=x(t)$, where $x:[0,\tau]\rightarrow \R$ solves 
\begin{align}\label{Vxt}
\begin{split}
\dot{x}(t)&= \VV(x(t))    \quad \text{ for } t\in [0,\tau],\quad  x(0) =x_0.
\end{split}
\end{align}
We use the simpler notation $\Tt=\Tt^{\VV}$ when no confusion is possible.
Let $\Om\in \mathcal P(\hold)$ and denote $n$ the outward unit normal vector to $\Omega$.
We consider the family of perturbed domains  
\begin{equation}\label{domain}
\Omega_t := \Tt^{\VV}(\Omega). 
\end{equation}
The choice of $\Theta^k(\hold)$ guarantees that $\Tt^{\VV}$ maps $\overline{\hold}$ onto $\overline{\hold}$, so that $\overline{\Omega_t}\subset \overline{\hold}$; see \cite[Theorem 2.16]{SokZol92}.
\begin{definition}\label{def1}
Let $J : \mathds{P} \rightarrow \R$ be a shape function.
\begin{itemize}
\item[(i)] The Eulerian semiderivative of $J$ at $\Omega$ in direction $\theta \in \Theta^k(\hold)$, when the limit exists,
is defined by
\ben
dJ(\Omega ;\VV):= \lim_{t \searrow 0}\frac{J(\Omega_t)-J(\Omega)}{t}.
\een
\item[(ii)] $J$ is  \textit{shape differentiable} at $\Omega$ if it has a Eulerian semiderivative at $\Omega$ for all $\theta \in \Theta^k(\hold)$ and the mapping
\begin{align*}
dJ(\Omega): \Theta^k(\hold) &  \to \R,\; \VV     \mapsto dJ(\Omega ;\VV)
\end{align*}
is linear and continuous, in which case $dJ(\Omega)$ is called the \textit{shape derivative} at $\Omega$.
\end{itemize}
\end{definition}
When the shape derivative is computed as a volume integral, it is convenient to write it in the following particular form.
\begin{definition}\label{def:tensor2}
Let $\Om\in \mathds{P}$ be open.
A shape differentiable function $J$ admits  a tensor representation of order $1$ if
there exist tensors $\Sb_l\in  L^1(\hold, \mathcal L^l(\R^m,\R^m))$, $l =0,1$, such that
\begin{align}
\label{ea:volume_form2}
dJ(\Omega ;\VV) = \int_\hold \Sb_1 \ddo D\VV +  \Sb_0\cdot \VV\,dx,
\end{align}
for all  $\VV\in \Theta^k(\hold)$. Here $\mathcal L^l(\R^m, \R^m) $ denotes the space of multilinear maps from 
$(\R^m)^l$ to $\R^m$.
\end{definition}
Expression \eqref{ea:volume_form2} is called distributed, volumetric, or domain expression of the shape derivative.
Under natural regularity assumptions, the shape derivative only depends on the restriction of the normal component $\VV\cdot n$ to the interface $\partial\Om$.  
This fundamental result is known as the Hadamard-Zol\'esio structure theorem in shape optimization; see \cite[pp. 480-481]{MR2731611}.  
From the tensor representation \eqref{ea:volume_form2}, one immediately obtains such structure of the shape derivative as follows.
\begin{proposition}\label{tensor_relations}
Let $\Om\in \mathds{P}$ and assume $\partial \Omega$ is $C^2$. 
Suppose that $dJ(\Omega)$ has the tensor representation \eqref{ea:volume_form2}.
If $\Sb_l$, $l=0,1$ are of class $W^{1,1}$ in $\Omega$ and $\hold\setminus \overline \Omega$, 
then we obtain the so-called {\it boundary expression} of the shape derivative: 
\ben
\label{eq:general_boundary_exp}
dJ(\Omega)(\theta) = \int_{\partial \Omega} g\, \theta\cdot n\, ds,
\een
with
$
g :=  [(\Sb_1^+-\Sb_1^-)n]\cdot n,
$
where $+$ and $-$ denote the restrictions of the tensor to $\Omega$ and $\hold\setminus\overline \Omega$, respectively.
\end{proposition}
See \cite{LaurainSturmM2AN} for a proof of Proposition \ref{tensor_relations} in a more general case.
Usually the boundary expression \eqref{eq:general_boundary_exp} is used to devise level set-based numerical methods, but in this paper we present an alternative approach based on the volume expression \eqref{ea:volume_form2}, which allows a simple implementation.
We use a Lagrangian approach to compute the tensor representation \eqref{ea:volume_form2}.

Further, we sometimes denote the distributed expression \eqref{ea:volume_form2}  by $\DV$,  and the boundary expression \eqref{eq:general_boundary_exp} by $\DS$ when we compare them. 
Note that if the domain is $C^2$, Proposition \ref{tensor_relations} shows that 
$$\DV=\DS.$$ 
When $\Om$ is less regular than $C^2$, it may happen that $dJ(\Omega)(\theta)$ cannot be written in the form \eqref{eq:general_boundary_exp}. 
Note that even in this case, 
$\DV$ is a distribution with support on the boundary, even if written as a domain integral. 
\section{Shape derivatives in the framework of linear elasticity}
\label{sec:shape_der_compliance}
\subsection{Shape derivative of the volume}\label{sec:volume}
We introduce a parameterized domain $\Om_t =\Tt^{\VV}(\Om)$ as in \eqref{domain}.
We start with the simple case of the volume 
$$\mathcal{V}(\Om_t) := \int_{\Om_t} 1\, dx,$$
which is useful to become familiar with the computation of shape derivatives.  
Using the change of variable $x\mapsto T_t(x)$, we get
$$\mathcal{V}(\Om_t) := \int_{\Om} \xi(t)\, dx,$$
where $\xi(t): = |\det DT_t| = \det DT_t$ for $t$ small enough.
We have  $\xi'(0) = \frac{d}{dt}\det DT_t|_{t=0} = \divv\VV$; see for instance \cite[Theorem 4.1, pp. 182]{MR2731611}.
Thus the distributed expression of the shape derivative of the volume is given by
$$d\mathcal{V}^{\tt{vol}}(\Om;\theta) = \int_{\Om} \divv\VV = \int_{\Om} I_d\ddo D\theta ,$$
where $I_d$ is the identity matrix. 
We have obtained the distributed expression \eqref{ea:volume_form2} of the shape derivative of the volume with $S_1 = I_d$ and $S_0 = 0$.

Applying Proposition \ref{tensor_relations}, assuming $\partial\Om$ is $C^2$, we get the usual boundary expression of the shape derivative
$$d\mathcal{V}^{\tt{surf}}(\Om;\theta) =  \int_{\partial\Om} \VV\cdot n,$$
which, in this case, is the same as applying Stokes' theorem.

\subsection{The ersatz material approach}\label{sec:ersatz}
We use the framework of the ersatz material, which is common in level set-based topology optimization of structures; see for instance \cite{MR2033390,MR1951408}. 
It is convenient as it allows to work on a fixed domain $\hold$ instead of the variable domain $\Om$, but also can create instability issues as pointed out in \cite{MR2674629}.
The idea of the ersatz material method is that the fixed domain $\hold$ is filled with two homogeneous materials with different Hooke elasticity tensors $A_{0}$ and $A_{1}$ defined by
$$A_{i}\xi=2\mu_{i}\xi+\lambda_{i}(\text{Tr}\xi)I_d,\quad i=0,1.$$
with Lam\'e moduli $\lambda_{i}$ and $\mu_{i}$ for $i=0,1$, $I_d$ is the identity matrix and $\xi$ is a matrix. 
The first material lays in the open subset $\Omega$ of $\hold$ and the background material fills the complement so that Hooke's law is written in $\hold$ as 
\begin{align}
\label{eq:ersatz1} A_\Om= A_0 \chi_{\Om} + \e A_0\chi_{\hold\setminus\Om},
\end{align}
where $\chi_{\Omega}$ denotes the indicator function of $\Omega$, and $\epsilon$ is a given small parameter. Hence the region $\hold\setminus\Omega$ represents a ``weak phase'' whereas $\Om$ is the ``strong phase''. 
The optimization is still performed with respect to the variable set $\Om$, but here $\Om$ is embedded in the fixed, larger set $\hold$. 

Note that we compute the shape derivative for the PDE including the ersatz material, unlike what is usually done in the literature; see Section \ref{sec:bd_distri_shapeder} for a more detailed discussion of this point. 

Let $\Omega\subset\hold \subset\R^m$, $m=2,3$, where $\hold$ is a fixed domain whose boundary $\partial\hold$ is partitioned into four subsets $\Gamma_{d}$, 
$\Gamma_{n}$, $\Gamma_{s}$  and $\Gamma$. 
A homogeneous Dirichlet (respectively  Neumann) boundary condition is imposed on $\Gamma_d$ (resp. $\Gamma$). 
On $\Gamma_n$, a non-homogeneous Neumann condition is imposed, which represents a given surface load $g\in H^{-1/2}(\Gamma_n)^m$.
The free interface between the weak and strong phase is $\partial\Om$.
A spring  with stiffness $k_s$ is attached on the boundary $\Gamma_s$, which corresponds to a Robin boundary condition; this condition is used for mechanisms.
Let $H^1_{d}(\hold)^m$ be the space of vector fields in $H^1(\hold)^m$ which satisfy the homogeneous Dirichlet boundary conditions on $\Gamma_d$. 

We define a parameterized domain $\Om_t =\Tt^{\VV}(\Om)$ as in \eqref{domain}, and we assume additionally that $\Tt^{\VV} = \operatorname{id}$ on $\Gamma_{d}\cup\Gamma_{n}\cup\Gamma_{s}\cup\Gamma_m$, where $\operatorname{id}$ is the identity.

In the ersatz material approach, the displacement field $u\in H_d^1(\hold)^m$ is the solution of the linearized elasticity system 
\begin{align}
\label{ersatz1.1}-\divv A_\Om e(u) &= 0 \mbox{ in }\hold,\\
\label{ersatz1.2}u & = 0 \mbox{ on } \Gamma_d,\\
\label{ersatz1.3}A_\Om e(u) n &= g  \mbox{ on } \Gamma_n,\\
\label{ersatz1.4}A_\Om e(u) n & = 0 \mbox{ on } \Gamma,\\
\label{ersatz1.5}A_\Om e(u) n & = -k_s u \mbox{ on } \Gamma_s,
\end{align}
where the symmetrized gradient is
$e(u) = (Du + Du^\transp )/2 ,$
and $Du^\transp$ denotes the transpose of $Du$. 
We consider the following functional
\begin{align} 
\label{shape_funk}
J(\Om) &: =\ca \int_\hold A_\Om e(u(x))\ddo e(u(x))\, dx + \cb \int_\hold \fd(x,u(x)) \, dx + \cc \int_{\Gamma_m} \fg(x,u(x))\, ds_x. 
\end{align}
We assume that $\fd$ and $\fg$ are smooth functions of $u$, that $\fd$ is $C^1$ with respect to the first argument, and $\Gamma_m \subset\partial\hold$.
The set $\Gamma_m \subset\partial\hold$ is a region where the shape displacements are monitored and is used for mechanisms only; it is set to $\Gamma_m = \emptyset$ in other cases.
This general functional covers several important cases such as the compliance and certain functionals used for compliant mechanisms.
The case $(\ca,\cb,\cc) = (1,0,0)$ corresponds to the compliance. 
The case of compliant mechanisms may be achieved by an appropriate choice of $\fd$ and $\fg$, see Section \ref{sec:inverter}.

We denote by $u_t$ the solution of \eqref{ersatz1.1}-\eqref{ersatz1.5} with $\Om$ substituted by $\Om_t$. 
Defining $u^t := u_t\circ \Tt$ and using the chain rule we have the relation
\begin{equation}
\label{Dut}
Du^t = D(u_t\circ \Tt) = (Du_t)\circ \Tt\, D\Tt, 
\end{equation}
and consequently
\begin{align}
\label{eq04} 
E(t,u^t) &:= e(u_t)\circ \Tt  = (Du^t D\Tt^{-1} + D\Tt^{-\transp}(Du^t)^\transp)/2.
\end{align}
The variational formulation of the PDE is to find $u_t\in H^1_d(\hold)^m$  such that  
\begin{equation}
\label{eq05}
\int_{\hold} A_{\Om_t}e(u_t)\ddo e(v_t) +\int_{\Gamma_s} k_s u_t\cdot v_t = \int_{\Gamma_n} g\cdot v_t, 
\end{equation}
for all $v_t\in H^1_d(\hold)^m$. 
We proceed with the change of variable $x\mapsto\Tt(x)$ in \eqref{eq05} which yields
\begin{align}
\label{55}
& \int_{\hold} [A_{\Om_t} e(u_t)]\circ \Tt\ddo e(v_t)\circ \Tt\, \xi(t) +\int_{\Gamma_s} k_s u_t\cdot v_t = \int_{\Gamma_n} g\cdot v_t, \mbox{ for all } v_t\in H^1_d(\hold)^m,
\end{align}
where $\xi(t): = |\det DT_t|$ is the Jacobian of the transformation   $x\mapsto\Tt(x)$.
Note that neither the Jacobian nor $T_t$ needs to appear in the integrals on $\Gamma_n$ and $\Gamma_s$, since we have assumed $T_t = \operatorname{id}$ on $\Gamma_{d}\cup\Gamma_{n}\cup\Gamma_{s}\cup\Gamma_m$.
In view of \eqref{eq04}, we can rewrite \eqref{55} as
\begin{equation}
\label{Eq_transported}\int_{\hold} A_\Om E(t,u^t)\ddo E(t,v)\, \xi(t) +\int_{\Gamma_s} k_s u^t\cdot v = \int_{\Gamma_n} g\cdot v, 
\end{equation}
for all $v\in H^1_{d}(\hold)^m$. 
In a similar way, we have
\begin{align*} 
J(\Om_t) & =\ca \int_\hold A_{\Om_t} e(u_t)\ddo e(u_t) + \cb \int_\hold \fd(x,u_t(x)) \, dx + \cc \int_{\Gamma_m} \fg(x,u_t(x))\, ds_x, 
\end{align*}
and using the change of variable $x\mapsto T_t(x)$ yields
\begin{align} 
\label{Comp_transported}
\begin{split}
J(\Om_t) & =\ca \int_\hold A_\Om E(t,u^t)\ddo E(t,u^t)\, \xi(t) + \cb \int_\hold \fd(T_t(x),u^t(x))\xi(t)(x) \, dx \\
& + \cc \int_{\Gamma_m} \fg(x,u^t(x))\, ds_x. 
\end{split}
\end{align}

\subsection{Shape derivative using a Lagrangian approach}\label{sec:lag_approach}

To compute the shape derivative of $J(\Om)$, we use the averaged adjoint method, a Lagrangian-type method introduced in \cite{MR3374631}.
Formally, the Lagrangian $G$ is obtained by summing the expression \eqref{Comp_transported} of the cost functional and the variational formulation \eqref{Eq_transported} of the PDE constraint, and  $(u^t,v)$ is replaced with the variables $(\vp,\psi)$; see \cite{LaurainSturmM2AN,MR3374631}
for a rigorous mathematical presentation and detailed explanations. 
Writing $A$ instead of $A_\Om$ for simplicity, this yields 
\begin{align*} 
G(t,\vp,\psi) := &  \ca \int_\hold A E(t,\vp)\ddo E(t,\vp)\, \xi(t) 
+ \cb \int_\hold \fd(T_t(x),\vp(x))\xi(t)(x) \, dx
+ \cc \int_{\Gamma_m} \fg(x,\vp(x))\, ds_x\\
& + \int_\hold AE(t,\vp)\ddo E(t,\psi)\, \xi(t) +\int_{\Gamma_s} k_s \vp \cdot\psi - \int_{\Gamma_n} g\cdot \psi. 
\end{align*}
In view of \eqref{Eq_transported} and \eqref{Comp_transported}, we have $J(\Om_t) = G(t,u^t,\psi)$ for all $\psi\in H^1_{d}(\hold)^m$. 
Thus the shape derivative can be computed as
\begin{equation}\label{eq:shape_der}
dJ(\Om;\VV) = \dt(G(t,u^t,\psi))|_{t=0}.   
\end{equation}
The advantage of the Lagrangian is that, under suitable assumptions, one can show that
\begin{align}\label{eq:dt_G_single}
\dt(G(t,u^t,\psi))|_{t=0} & =\partial_t G(0,u^0,p^0).
\end{align}
which essentially means that it is not necessary to compute the derivative of $u^t$ to compute $dJ(\Om;\VV)$.
In this paper we assume for simplicity that \eqref{eq:dt_G_single} is true for the problem under consideration, but note that this result can be made mathematically rigorous using the averaged adjoint method; see \cite{doi:10.1137/15100477X,LaurainSturmM2AN,MR3374631}.
The adjoint is given as the solution of the following first-order optimality condition
$$ \partial_{\vp}G(0,u,p)(\hat\vp) = 0 \mbox{ for all }\hat\vp\in  H^1_{d}(\hold)^m,$$
which yields, using $A=A^\transp$, 
\begin{align*}
& 2\ca\int_\hold AE(0,u)\ddo E(0,\hat\vp)
+ \cb \int_\hold \partial_u \fd(x,u(x))\cdot\hat\vp(x) \, dx
+ \cc \int_{\Gamma_m} \partial_u\fg(x,u(x))\cdot\hat\vp(x)\, ds_x\\
& +\int_{\Gamma_s} k_s \hat\vp \cdot p
+ \int_\hold AE(0,\hat\vp)\ddo E(0,p) = 0, \mbox{ for all } \hat\vp\in  H^1_{d}(\hold)^m. 
\end{align*}
Since $E(0,v)=e(v)$ for $v\in H^1_{d}(\hold)^m$, and $A=A^\transp$,  we get the  adjoint equation
\begin{align}
\label{eq:203}
\begin{split}
&\int_\hold A e(p)\ddo e(\hat\vp) +\int_{\Gamma_s} k_s p\cdot\hat\vp  
=  -2\ca\int_\hold A e(u)\ddo e(\hat\vp)  - \cb \int_\hold \partial_u \fd(x,u(x))\cdot\hat\vp(x) \, dx\\
& - \cc \int_{\Gamma_m} \partial_u\fg(x,u(x))\cdot\hat\vp(x)\, ds_x\ \mbox{ for all }\hat\vp\in  H^1_{d}(\hold)^m.
\end{split}
\end{align}
Using \eqref{eq:shape_der} and \eqref{eq:dt_G_single}, we obtain
\begin{align*} 
dJ(\Om;\VV) =& \ca\int_\hold  A \partial_t E(0,u)\ddo E(0,u) 
+\ca\int_\hold  A  E(0,u)\ddo \partial_t E(0,u) + A E(0,u)\ddo E(0,u) \divv\VV \\
& + \cb \int_\hold \partial_x\fd(x,u(x))\cdot \theta(x) +  \fd(x,u(x))\divv\theta(x) \, dx\\
& +\int_\hold  A \partial_t E(0,u)\ddo E(0,p) +\int_\hold  A  E(0,u)\ddo \partial_t E(0,p) + A E(0,u)\ddo E(0,p) \divv\VV
\end{align*}
We also compute, using $E(0,v) = e(v)$ for $v\in H^1_{d}(\hold)^m$,
$$ \partial_t E(0,v) = (-Dv D\VV - D\VV^\transp Dv^\transp)/2.$$
Using $A=A^\transp$  we obtain
\begin{align*} 
dJ(\Om;\VV)  = &  -\int_\hold  \frac{1}{2} Du^\transp (Ae(p) +  (Ae(p))^\transp )\ddo D\VV
-\int_\hold  \frac{1}{2} Dp^\transp (Ae(u) +  (Ae(u))^\transp )\ddo D\VV\\
& + \int_\hold  (A e(u)\ddo e(p) + \ca A e(u)\ddo e(u)) \divv\VV 
-\ca\int_\hold  (Du^\transp Ae(u) +Du^\transp  (Ae(u))^\transp )\ddo D\VV\\
& + \cb \int_\hold \partial_x\fd(x,u(x))\cdot \theta(x) +  \fd(x,u(x))\divv\theta(x) \, dx.
\end{align*}
Using $(Ae(v))^\transp  = Ae(v)$ we get
\begin{align} 
\label{shape_der_ersatz0} 
\DV =&  \int_{\hold} \Sb_1 \ddo D\VV +\Sb_0\cdot\VV,
\end{align}
with 
\begin{align} 
\notag\Sb_1 =& -Du^\transp A_\Om e(p) - Dp^\transp A_\Om e(u) - 2c_1 Du^\transp A_\Om e(u)\\
\label{Sb1_ersatz}    &  +(A_\Om e(u)\ddo e(p) + \ca A_\Om e(u)\ddo e(u) +\cb \fd(\cdot,u)) I_d,\\
\label{Sb1_ersatz2}  \Sb_0 =& c_2\partial_x\fd(\cdot ,u) .
\end{align}
Formula \eqref{shape_der_ersatz0} is convenient for the numerics as it can be implemented in a straightforward way in FEniCS.
\subsection{Compliance}
\begin{figure}
\begin{center}
\begin{subfigure}[b]{0.4\textwidth}
\includegraphics[width=\textwidth]{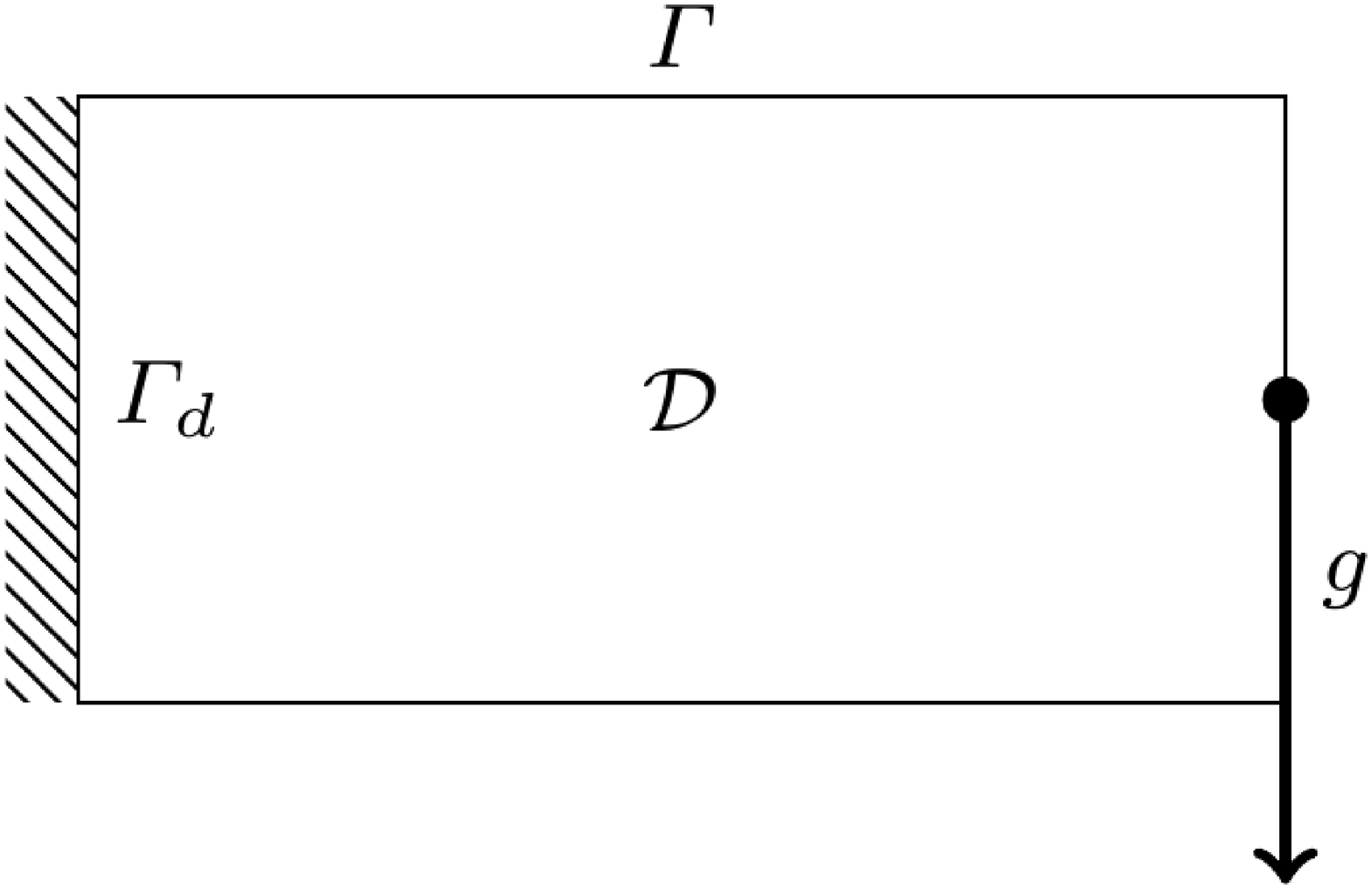}
\end{subfigure}
\begin{subfigure}[b]{0.4\textwidth}
\raisebox{5mm}{
\includegraphics[width=\textwidth]{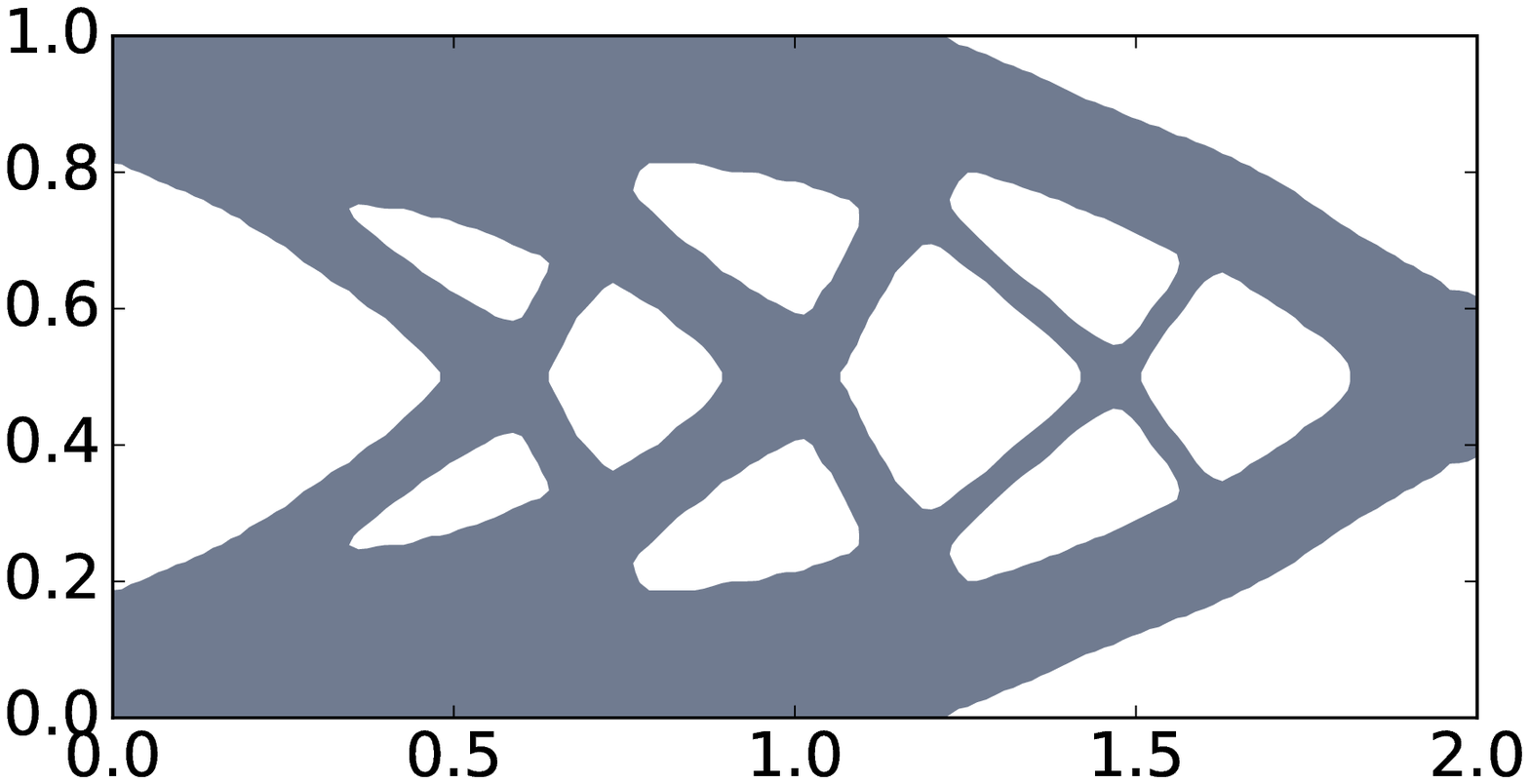}
}
\end{subfigure}
\end{center}
\caption{Design domain of cantilever (left) and example of optimal design (right). }\label{fig:cantilever_fig}
\end{figure}

The case of the compliance is obtained by setting $(\ca,\cb,\cc) = (1,0,0)$. 
This yields  $\Sb_0 \equiv 0$,  $p=-2u$ and \eqref{Sb1_ersatz} becomes  
\begin{equation} \label{Sb1_ersatz_compliance}
\Sb_1 = 2Du^\transp A_\Om e(u)-  A_\Om e(u)\ddo e(u) I_d. 
\end{equation}
See Figure \ref{fig:cantilever_fig} for an example of design domain, boundary conditions, and optimal design for minimization of the compliance. 

\subsection{Multiple load cases}\label{sec:multiple_load}
For multiple load cases and compliance, we consider the set of forces $\{g_i\}_{i\in I}$ and the compliance is the sum of the compliances associated to each force $g_i$:
\begin{align*} 
J(\Om) := \sum_{i\in I} \int_{\Gamma_n} g_i\cdot u_i, 
\end{align*}
where $u_i$ is the solution of the linearized elasticity system corresponding to $g_i$.
The shape derivative is in this case
\begin{align} \label{eq:35}
\DV =& \int_\hold  \Sb_1 \ddo D\VV,
\end{align}
with 
$
\Sb_1 := \sum_{i\in I} 2Du_i^\transp A_\Om e(u_i)-  A_\Om e(u_i)\ddo e(u_i) I_d . 
$
\subsection{Inverter mechanism}\label{sec:inverter}
\begin{figure}
\begin{center}
\begin{subfigure}[b]{0.4\textwidth}
\includegraphics[width=\textwidth]{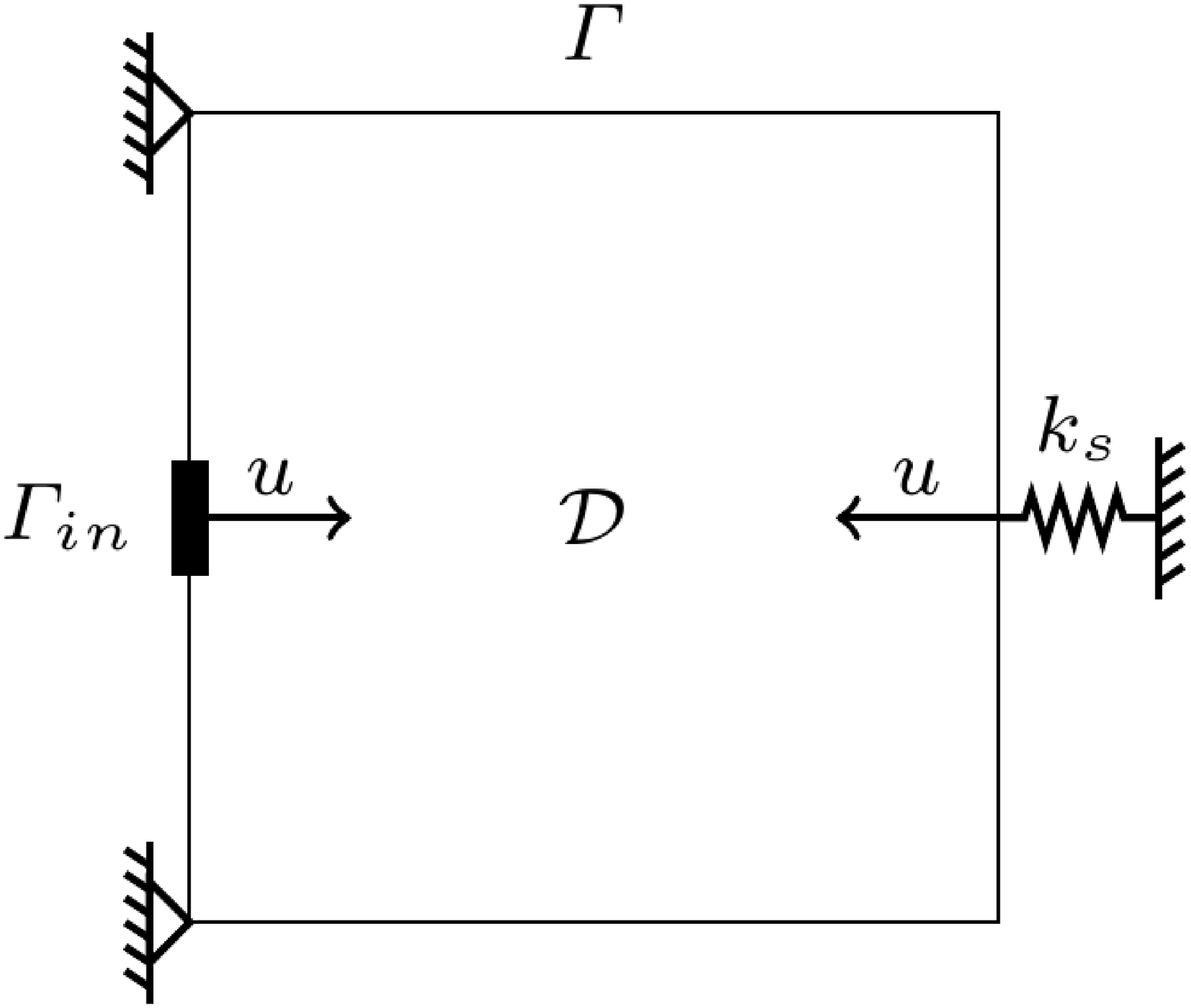}
\end{subfigure}
\begin{subfigure}[b]{0.3\textwidth}
\raisebox{2mm}{
\includegraphics[width=\textwidth]{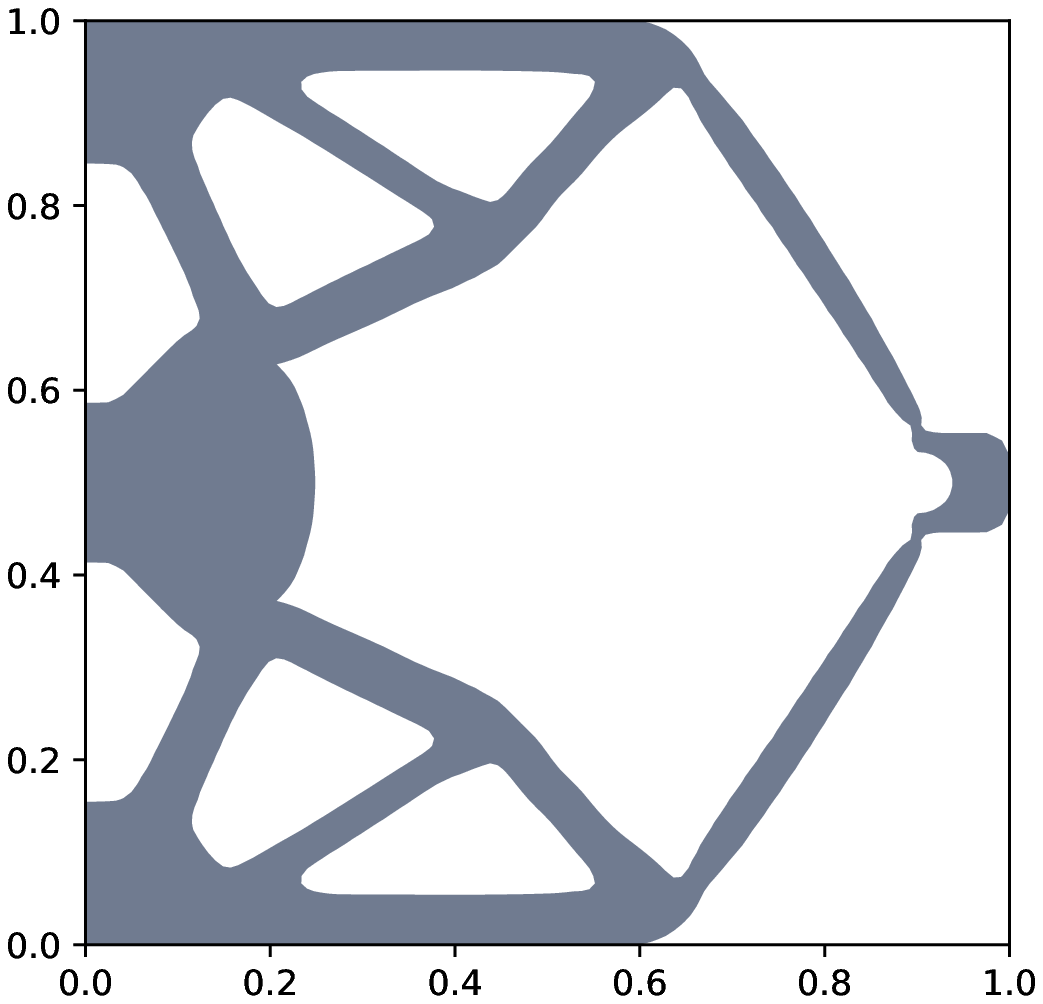}
}
\end{subfigure}
\end{center}
\caption{Design domain of inverter (left) and optimal design (right) for $(Nx,Ny)=(121,121)$. }\label{fig:mechanism_fig}
\end{figure}

The displacement inverter converts an input displacement on the
left edge to a displacement in the opposite direction on the right edge; see \cite{MR2008524} for a detailed description.
We take  $\hold = (0,1)^2$, and an actuation force $g = (g_x,0)$, $g_x>0$, is applied at the input point  $(0,0.5)$. 
We define the output boundary $\Gamma_{out}$ and input boundary $\Gamma_{in}$ such that $\Gamma_m = \Gamma_{out} \cup \Gamma_{in}$, $\Gamma_{in} = \{0\}\times (a_0,a_1)$, $\Gamma_{out} = \Gamma_s$ and  $\Gamma_s = \{1\}\times (b_0,b_1)$.
An artificial spring with stiffness $k_s>0$ is attached at the output $\Gamma_{out}$ to simulate the resistance of a workpiece.
In order to maximize output displacement, while limiting the input displacement, we minimize \eqref{shape_funk} with $(\ca,\cb,\cc) = (0,0,1)$ and
$$\fg(x,u(x)) = \eta_{in}u_1(x)\chi_{\Gamma_{in}}(x) + \eta_{out} u_1(x)\chi_{\Gamma_{out}}(x) ,$$
with $u = (u_1,u_2)$, $\eta_{in},\eta_{out}$ some positive constants. Note that $u_1>0$ on $\Gamma_{in}$ and $u_1<0$ on $\Gamma_{out}$.
We obtain $\Sb_0 \equiv 0$ and \eqref{Sb1_ersatz} becomes 
\begin{align*} 
\Sb_1 &= -Du^\transp A_\Om e(p) - Dp^\transp A_\Om e(u) + A_\Om e(u)\ddo e(p) I_d.
\end{align*}
For this functional the problem is not self-adjoint and in view of \eqref{eq:203}, the adjoint $p$ is the solution of
\begin{align}
\label{adjoint_inverter}
\int_\hold  A e(p)\ddo e(\hat\vp) +\int_{\Gamma_s} k_s p\cdot\hat\vp 
=    -\eta_{in}\int_{\Gamma_{in}} \hat\vp_1
- \eta_{out}\int_{\Gamma_{out}} \hat\vp_1,\quad \forall \hat\vp\in  H^1_{d}(\hold)^m.
\end{align}
where $\hat\vp = (\hat\vp_1,\hat\vp_2)$.
See Figure \ref{fig:mechanism_fig} for an example of design domain, boundary conditions and optimal design for the inverter. 

\subsection{Comparison of shape derivatives with and without ersatz material}\label{sec:bd_distri_shapeder}

Usually the boundary expression of the shape derivative is used in level set methods, and computed for the problem without ersatz material, although the elasticity system is solved using the ersatz material in the numerics. 
This small mismatch is justified by the fact that the tensor of the ersatz material has a small amplitude. 
The reason why this mismatch is tolerated in the numerics is probably because 
the boundary expression of the shape derivative in that case is unpractical to handle numerically, as it requires to use the jump of the gradient across the moving interface between the strong and the weak phases; see  \eqref{shape_der_ersatz}. %

In any case, it is more precise to use the proper shape derivative corresponding to the ersatz material framework for the numerics, in order to avoid this mismatch. 
Another advantage of using the exact formula for the ersatz approach is that this formula is actually valid for any value of $\e$, and not only for $\e$ small. 
This can be used for a mixture of two materials for instance.
We show in this section that computing and implementing the formula of the distributed shape derivative is not more difficult for the ersatz material approach.

First we compare the distributed shape derivative without ersatz material. 
The elasticity system is in this case
\begin{align}
\label{el1} -\divv A e(u) &= 0 \mbox{ in }\Om,\\
\label{el2} u & = 0 \mbox{ on } \Gamma_d,\\
\label{el3} Ae(u) n &= g  \mbox{ on } \Gamma_n,\\
\label{el4} Ae(u) n & = 0 \mbox{ on } \Gamma,\\
\label{el5} Ae(u) n & = -k_s u \mbox{ on } \Gamma_s.
\end{align}
As in Section \ref{sec:ersatz}, $\Gamma_d$, $\Gamma_n$ and $\Gamma_s$ are fixed, but in this case $\Gamma = \partial\Om\setminus (\Gamma_n\cup\Gamma_d \cup \Gamma_s)$ is the free boundary. 
We also assume that the interface between $\Gamma$ and the fixed boundaries is fixed.
The Hooke elasticity tensor $A$ satisfies 
$A\xi = 2\mu\xi +\lambda\tr(\xi) I_d$,
where $\xi\in\R^{m\times m}$, $\mu,\lambda$ are the Lam\'e parameters. 
The variational formulation of \eqref{el1}-\eqref{el5} consists in finding $u\in H^1_{d}(\Omega)^m$  such that  
$$\int_\Om Ae(u)\ddo e(v) + \int_{\Gamma_s} k_s u\cdot v = \int_{\Gamma_n} g\cdot v\mbox{ for all }v\in H^1_{d}(\Omega)^m.$$
The cost functional is in this case
\begin{align*} 
J_0(\Om_t) : =\ca \int_\Om Ae(u)\ddo e(u) + \cb \int_\Om \fd(x,u(x)) \, dx
+ \cc \int_{\Gamma_m} \fg(x,u(x))\, ds_x. 
\end{align*}
A similar calculation as in Section \ref{sec:lag_approach} yields
\begin{align*} 
dJ_0^{\tt{vol}}(\Om;\VV) =&  \int_{\Om} \Sb_1 \ddo D\VV +\Sb_0\cdot\VV,
\end{align*}
with 
\begin{align} 
\Sb_1 &= -Du^\transp A e(p) - Dp^\transp Ae(u) - 2c_1 Du^\transp A e(u)
\label{Sb1_ersatz_b}   +(A e(u)\ddo e(p) + \ca A e(u)\ddo e(u) +\cb \fd(\cdot,u)) I_d,\\
\label{Sb1_ersatz_b2}  \Sb_0 &= c_2 \partial_x\fd(\cdot ,u) .
\end{align}
In the case $(\ca,\cb,\cc) = (1,0,0)$, which  corresponds to the compliance, we have  $\Sb_0\equiv 0$ and $p=-2u$, yielding  
\begin{align} \label{Sb1}
\Sb_1 := 2Du^\transp Ae(u)-  A e(u)\ddo e(u) I_d \quad \mbox{ in }\Om. 
\end{align}
A similar formula can be found in \cite[Section 2.5]{MR3013681}, for a slightly different case, and where $S_1$ is identified as the energy-momentum tensor in continuum mechanics introduced by Eshelby  in \cite{MR0489190}.
Compare also \eqref{Sb1} with the shape derivative in \cite[Theorem 3.3]{MR3493474}, also in the framework of linearized elasticity but for a different functional.

Note that \eqref{Sb1} is similar to \eqref{Sb1_ersatz_compliance}, the main difference being that  \eqref{Sb1_ersatz_compliance} is defined in $\hold$ and \eqref{Sb1} is defined in $\Om$.
Thus from a numerical point of view, \eqref{Sb1_ersatz_compliance} is not more difficult to implement than \eqref{Sb1}.

Now we compare the boundary expression of the shape derivatives with and without ersatz material, in the case of the compliance.
Assuming $\Om$ is  $C^2$ and using \eqref{Sb1} and \eqref{eq:general_boundary_exp}  of Proposition \ref{tensor_relations},  we obtain the boundary expression of the shape derivative
\begin{align*} 
dJ_0^{\tt{surf}}(\Om;\VV)=& 
\int_{\partial\Om} (\Sb_1n\cdot n) \VV\cdot n = \int_{\Gamma} (\Sb_1n\cdot n) \VV\cdot n,
\end{align*}
since $\partial\Om\setminus\Gamma$ is fixed. 
Then we compute
\begin{align*} 
\Sb_1n\cdot n &= 2Du^\transp Ae(u)n\cdot n-  A e(u)\ddo e(u)
=2Ae(u)n\cdot Du n-  A e(u)\ddo e(u). 
\end{align*}
On $\Gamma$, we have $Ae(u)n = 0$ which yields
\begin{align} 
\label{shape_der_boundary_0}
\begin{split}
dJ_0^{\tt{surf}}(\Om;\VV) =&  \int_{\Gamma} -A e(u)\ddo e(u) \VV\cdot n,
\end{split}
\end{align}  
which is a particular case of the formula in \cite[Theorem 7]{MR2033390}.

In the case of the ersatz material, applying Proposition \ref{tensor_relations} and assuming  $\partial\Om$ is $C^2$, the distributed expression \eqref{shape_der_ersatz0} yields the boundary expression 
\begin{align*}
\DS=&  \int_{\partial\Om} [(S_1^+ - S_1^-)n]\cdot n\ \VV\cdot n,
\end{align*}
with 
\begin{align*}
&[(S_1^+ - S_1^-)n]\cdot n = -\llbracket A e(u)\ddo e(u)\rrbracket
 + 2A_0 e(u)^+ n\cdot Du^+ n - 2\e A_0 e(u)^- n\cdot Du^- n, 
\end{align*}
and where the exponents $(\cdot)^+$ and $(\cdot)^-$ denote the restrictions to $\Om$ and $\hold\setminus\Om$, respectively. Also
$$\llbracket v\rrbracket := \gamma_{\Om} (v) - \gamma_{\hold\setminus\Om} (v) $$ 
denotes the jump of a function $v$ across the interface $\partial\Om$; here $\gamma_{\Om} (v)$ is the trace of $v_{|{\Om}}$ on $\partial\Om$. 
Using the transmission condition 
$ A_0 e(u)^+ n =  \e A_0 e(u)^- n \mbox{ on }\partial\Om$,
we obtain
\begin{align} 
\label{shape_der_ersatz} 
\begin{split}
\DS=&  \int_{\partial\Om} (2\e A_0 e(u^-)n\cdot \llbracket Du\rrbracket n) \VV\cdot n 
 -  \int_{\partial\Om} \llbracket A e(u)\ddo e(u)\rrbracket  \VV\cdot n.
\end{split}
\end{align} 
The two main differences between \eqref{shape_der_ersatz} and \eqref{shape_der_boundary_0} are the small perturbation term
$(2\e A_0 e(u^-)n\cdot \llbracket Du \rrbracket n $
and the fact that $\llbracket A e(u)\ddo e(u)\rrbracket$ is a jump across the interface $\partial\Om$.
We observe that \eqref{shape_der_boundary_0} is easier to implement than \eqref{shape_der_ersatz} in a numerical method. 
\section{Descent direction}\label{sec:descent}
For the numerical method we need a descent direction $\VV$, i.e. a vector field satisfying  
$dJ(\Omega ;\VV) < 0$.
When $dJ(\Omega ;\VV)$ is written using the boundary expression 
$$\DS= \int_{\partial\Om} G(\Om)\, \VV\cdot n,$$
then a simple choice is to take $\VV = - G(\Om)\cdot n$. 
However, this choice assumes that $G(\Om)$ and $\partial\Om$ are quite regular, and in practice this may yield a $\VV$ with a poor regularity and lead to an unstable behaviour of the algorithm such as irregular or oscillating boundaries.
A better choice is to find a smoother descent direction by finding $\VV\in \mathds{H}(\partial\Om)$  such that    
\begin{equation}
\label{VP0}
\mathcal{B}(\VV,\xi)=- dJ^{\tt{surf}}(\Omega;\xi) \mbox{  for all  } \xi\in \mathds{H}(\partial\Om) ,
\end{equation}
where $\mathds{H}(\partial\Om)$ is an appropriate Sobolev space of vector fields on $\partial\Om$ and $\mathcal{B}: \mathds{H}(\partial\Om)\times \mathds{H}(\partial\Om)\to\R$, $k\geq 1$, is a positive definite bilinear form on $\partial\Om$. 

In the case of the present paper we use the distributed expression \eqref{shape_der_ersatz0} of the shape derivative, therefore we use a positive definite bilinear form $\mathcal{B}: \mathds{H}(\hold) \times \mathds{H}(\hold)\to\R$, where $\mathds{H}(\hold)$ is an appropriate Sobolev space of vector fields on $\hold$. Thus the problem is to find $\VV\in \mathds{H}(\hold)$  such that    
\begin{equation}
\label{VP_1}
\mathcal{B}(\VV,\xi)=-  dJ^{\tt{vol}}(\Omega;\xi) \mbox{  for all  } \xi\in \mathds{H}(\hold),
\end{equation}
With this choice, the solution $\VV$ of \eqref{VP_1} is defined on all of $\hold$ and is a descent direction since $dJ(\Omega ;\VV) = -\mathcal{B}(\VV,\VV)< 0$ if $\VV\neq 0$.

It is also possible to combine the two approaches by substituting $\DV$ with $\DS$ in \eqref{VP_1}. This was done in \cite{MR2225309} where a strong improvement of the rate of convergence of the level-set method was observed; see also \cite{MR1998617} for a thorough discussion of various possibilities for $\mathcal{B}$. Bilinear forms defined on $\hold$ are  useful for the level set method which requires $\VV$ on $\hold$; see Section \eqref{section5}.

In our algorithm we choose $\mathds{H}(\hold) = H^1(\hold)^m$ and
\begin{equation}
\label{VP_2}
\mathcal{B}(\VV,\xi) = \int_\hold \alpha_1 D\VV\ddo D\xi +  \alpha_2 \VV\cdot \xi,
\end{equation}
with $\alpha_1 = 1$ and $\alpha_2 = 0.1$. 
We also take the boundary conditions $\VV\cdot n = 0$ on $\partial\hold$; see  Section \ref{subsec:descent}.

\section{Level set method}\label{section5}
The level set method, originally introduced in \cite{MR965860}, gives a general framework for the computation of evolving interfaces using an implicit representation of these interfaces.  
We refer to the monographs \cite{MR1939127,MR1700751} for a complete description of the level set method.
The core idea of this method is to represent the boundary of the moving domain $\Omega_t\subset \hold \in \R^N$ as the zero level set of a continuous function $\phi(\cdot,t): \hold\to\R$.

Let us consider the family of domains $\Omega_t\subset \hold$ as defined in \eqref{domain}. Each domain $\Omega_t$ can be defined as
\begin{equation}
\Omega_t :=\{x\in \hold,\ \phi(x,t) < 0\},
\end{equation}
where $\phi: \hold\times \R^+ \to \R$ is Lipschitz continuous and called {\it level set function}. 
Indeed, if we assume $|\nabla\phi(\cdot,t)|\neq 0$ on the set  $\{x\in \hold,\ \phi(x,t) = 0\}$, then we have
\begin{equation}
\partial \Omega_t = \{x\in \hold,\ \phi(x,t) = 0\},
\end{equation}
i.e. the boundary $\partial \Omega_t$ is the zero level set of  $\phi(\cdot,t)$.

Let $x(t)$ be the position of a moving boundary point of
$\partial \Omega_t$, with velocity $\dot{x}(t)=\VV(x(t))$ according to \eqref{Vxt}.
Differentiating the relation $\phi(x(t),t)=0$ with respect to $t$
yields the Hamilton-Jacobi equation:
\begin{equation*}
\partial_t\phi (x(t),t)+ \VV(x(t))\cdot \nabla \phi(x(t),t) = 0  \quad \mbox{ in } \partial \Omega_t\times\R^+,
\end{equation*}
which is then extended to all of $\hold$ via the equation
\begin{equation}
\partial_t\phi(x,t) + \VV(x)\cdot \nabla \phi(x,t) = 0  \quad \mbox{ in } \hold\times\R^+,
\label{eq:transport}
\end{equation}
or alternatively to $U\times\R^+$ where $U$ is a neighbourhood of $\partial \Omega_t$.

When $\VV = \vartheta_n n$ is a normal vector field on $\partial\Omega_t$, noting that an extension to $\hold$ of the unit outward normal vector $n$ to $\Omega_t$ is given by $\nabla \phi /|\nabla \phi|$, and extending $\vartheta_n$ to all of $\hold$, one obtains from \eqref{eq:transport} the level set equation
\begin{equation}
\partial_t \phi + \vartheta_n |\nabla \phi| = 0 \quad \mbox{ in } \hold\times\R^+.
\label{eq:HJ1}
\end{equation}

The initial data $\phi(x,0)=\phi_0(x)$ accompanying the Hamilton-Jacobi equation \eqref{eq:transport} or \eqref{eq:HJ1} can be chosen as the signed distance function to the
initial boundary $\partial \Omega_0$ in order to satisfy the condition $|\nabla u|\neq 0$ on $\partial\Omega$, i.e.
\begin{equation}
\phi_0(x) = 
\left\{\begin{array}{rl}
d(x,\partial \Omega_0), & \mbox{ if } x\in (\Omega_0)^c, \\
-d(x,\partial \Omega_0), & \mbox{ if } x\in \Omega_0 .
\end{array}
\right.
\end{equation}
The fast marching method \cite{MR1700751} and the fast sweeping method \cite{MR2114640} are efficient methods to compute the signed distance function.

\subsection{Level set method and volume expression of the shape derivative}

In the case of the distributed shape derivative  \eqref{shape_der_ersatz0}, we do not extend $\vartheta_n$ to $\hold$, instead we obtain directly a descent direction $\VV$ defined in $\hold$ by solving \eqref{VP_1}, where $\DV$ is given by \eqref{shape_der_ersatz0}. 
Thus, unlike  the usual level set method, $\VV$ is not necessarily normal to $\partial\Om_t$ and $\phi$ is not governed by \eqref{eq:HJ1} but rather by the Hamilton-Jacobi equation \eqref{eq:transport}.

In shape optimization,  $\vartheta_n$ usually depends on the solution of one or several PDEs and their gradient. Since the boundary $\partial\Omega_t$  in general does not match the grid nodes where $\phi$ and the solutions of the  partial differential equations are defined in the numerical application, the computation and extension of $\vartheta_n$ may
require the  interpolation on $\partial\Omega_t$ of functions defined at the grid points only, complicating the numerical implementation  and introducing an 
additional interpolation error. This is an issue in particular for interface problems, such as the problem of elasticity with ersatz material studied in this paper, where $\vartheta_n$ is the jump of a function across the interface, as in \eqref{shape_der_ersatz}, which requires several interpolations and is error-prone. In the distributed shape derivative framework,  $\VV$ only needs to be defined at grid nodes.

\subsection{Discretization of the Hamilton-Jacobi equation}\label{sec:discretization_hj}
Let $\hold=(0,1)\times(0,1)$ to simplify the presentation.
For the discretization of the Hamilton-Jacobi equation \eqref{eq:transport},
we first define the mesh grid corresponding to $\hold$. We introduce the nodes $P_{ij}$
whose coordinates are given by $(i\Delta x, j \Delta y)$, $1\leq i,j\leq N$ where $\Delta x$ 
and $\Delta y$ are the steps of the discretization in the $x$ and $y$
directions, respectively. 
Let us write $t^k = k\Delta t$ for
the discrete time, with $k\in \mathds{N}$ and $\Delta t$ is the time
step. 
Denote the approximation
$\phi_{ij}^k \simeq \phi(P_{ij},t^k)$.

In the usual level set method, equation \eqref{eq:HJ1} is discretized using an explicit upwind scheme proposed by Osher and Sethian \cite{MR1939127,MR965860,MR1700751}. This scheme applies to the specific form \eqref{eq:HJ1} but is not suited to discretize \eqref{eq:transport} required for our application. Equation \eqref{eq:transport} is of the form
\begin{equation}
\partial_t\phi + H(\nabla \phi)= 0  \quad \mbox{ in } \hold\times\R^+,
\end{equation}
where $ H(\nabla \phi) := \VV\cdot \nabla \phi $ is the so-called Hamiltonian. 
We use a Lax-Friedrichs flux, see \cite{MR1111446}, which writes in our case:
\begin{align*}
& \hat H^{LF}(p^-,p^+,q^-,q^+)  
= H\left(\frac{p^- + p^+}{2},\frac{q^- +q^+}{2}\right) 
-\frac{1}{2} (p^+-p^-)\alpha^x -\frac{1}{2}(p^+-p^-)\alpha^y ,  
\end{align*}
where $\alpha^x = |\VV_x|$, $\alpha^y = |\VV_y|$, $\VV = (\VV_x,\VV_y)$ and
\begin{align}
\label{pq}
\begin{split}
p^- &= D_x^{-}\phi_{ij}=\dfrac{\phi_{ij}-\phi_{i-1,j}}{\Delta x},\qquad
p^+ = D_x^{+}\phi_{ij}=\dfrac{\phi_{i+1,j}-\phi_{ij}}{\Delta x},\\
q^- &= D_y^{-}\phi_{ij}=\dfrac{\phi_{ij}-\phi_{i,j-1}}{\Delta y},\qquad
q^+ = D_y^{+}\phi_{ij}=\dfrac{\phi_{i,j+1}-\phi_{ij}}{\Delta y},
\end{split}
\end{align}
are the backward and forward approximations of the $x$-derivative and $y$-derivative of 
$\phi$ at $P_{ij}$, respectively. Using a forward Euler time discretization, the numerical scheme corresponding to \eqref{eq:transport} is 
\begin{equation}
\label{num_scheme_HJ}
\phi_{ij}^{k+1}=\phi_{ij}^k - \Delta t \ \hat H^{LF}(p^-,p^+,q^-,q^+)
\end{equation}
where $p^-,p^+,q^-,q^+$ are computed for $\phi_{ij}^k$.
\subsection{Reinitialization}\label{reinit}
For numerical accuracy, the solution of the level set equation 
\eqref{eq:transport} should not be too flat or too steep. 
This is fulfilled for instance if $\phi$ is the distance function
i.e. $|\nabla \phi| = 1$. 
Even if one initializes $\phi$ using a
signed distance function, the solution
$\phi$ of the level set equation \eqref{eq:transport} does not generally
remain close to a distance function, 
thus we regularly perform a reinitialization
of $\phi$; see \cite{MR1214016}.

We present here briefly the procedure for the reinitialization introduced in \cite{MR1723321}.
The reinitialization at time $t$ is performed by solving to steady state the following Hamilton-Jacobi type equation
\begin{align*}
\partial_\tau\varphi + S(\phi)(|\nabla\varphi|-1) & = 0 \mbox{ in }\hold\times\R^+,\\
\varphi(x,0) & = \phi(x,t), \ x\in \hold,
\end{align*}
where $S(\phi)$ is an approximation of the sign function
\begin{equation}\label{signum}
S(\phi) = \frac{\phi}{\sqrt{\phi^2 + |\nabla\phi|^2\epsilon_s^2}}, 
\end{equation}
with $\epsilon_s = \min(\Delta x,\Delta y)$.

For the discretization we use the standard explicit upwind
scheme; see \cite{MR1939127,MR965860,MR1700751},
\begin{align}\label{reinit_phi_ij}
\varphi_{ij}^{k+1} = \varphi_{ij}^{k} - \Delta t\ \Gnum (p^-,p^+,q^-,q^+),   
\end{align}
where 
\begin{align}\label{computeG}
& \Gnum(p^-,p^+,q^-,q^+) 
 = \max(S(\phi_{ij}),0) \Gnum^+ + \min(S(\phi_{ij}),0) \Gnum^-,  
\end{align}
and 
\begin{align}
\label{Gplus} \Gnum^+ & = \left[ \max(p^-,0)^2 + \min(p^+,0)^2  + \max(q^-,0)^2 + \min(q^+,0)^2\right]^{1/2},\\
\label{Gminus} \Gnum^- & = \left[ \min(p^-,0)^2 + \max(p^+,0)^2 + \min(q^-,0)^2 + \max(q^+,0)^2\right]^{1/2},
\end{align}
and where $p^-,p^+,q^-,q^+$ are computed for $\phi_{ij}^k$ using \eqref{pq}.

\section{Implementation}\label{sec:numerical_implementation}
In this section we explain the implementation step by step. 
The code presented in this paper has been written for FEniCS 2017.1, and is compatible with FEniCS 2016.2.
With a small number of modifications, the code may also run with earlier versions of FEniCS.
The code can be downloaded at  \burl{http://antoinelaurain.com/compliance.htm}.
\subsection{Introduction}
We explain the code for the case of the compliance, i.e. $(\ca,\cb,\cc)=(1,0,0)$ in \eqref{shape_funk} and $\Gamma_s =\emptyset$, and we consider an additional volume constraint, so the functional that we minimize is
\begin{equation}\label{Jnum}
\mathcal{J}(\Om) := J(\Om) + \Lambda  \mathcal{V}(\Om),
\end{equation}
where $\Lambda$ is a constant and $\mathcal{V}(\Om)$ is the volume of $\Om$.

The main file \f{compliance.py} can be found in the appendix, and 
we use a file \f{init.py} to initialize the data which depend on the chosen case. 
The user can choose between the six following cases:  \f{half_wheel}, \f{bridge}, \f{cantilever}, \f{cantilever_asymmetric}, \f{MBB_beam},  and \f{cantilever_twoforces}. 
For instance, to run the cantilever case, the command line is 
\begin{lstlisting}[numbers=none,xleftmargin=0cm]
python compliance.py cantilever
\end{lstlisting}

An important feature of the code is that we use two separate grids.
On one hand, $\hold$ is discretized using a structured grid \f{mesh} made of isosceles triangles (each square is divided into four triangles), which is used to compute the solution \f{U} of the elasticity system, and also to compute the descent direction \f{th} corresponding to $\VV$.
The spaces \f{V} and \f{Vvec} are spaces of scalar and vector-valued functions on \f{mesh}, respectively. 
On the other hand, we use an additional Cartesian grid, whose vertices are included in the set of vertices of \f{mesh}, to implement the numerical scheme \eqref{num_scheme_HJ} to solve the Hamilton-Jacobi equation, and also to perform the reinitialization  \eqref{reinit_phi_ij}.
In \f{compliance.py}, the quantities defined on the Cartesian grid are matrices and therefore distinguished by the suffix \f{mat}.
For instance \f{phi} is a function defined on \f{mesh}, while \f{phi_mat} is the corresponding function defined on the Cartesian grid. 
We need a mechanism to alternate between functions defined on \f{mesh} and functions defined on the Cartesian grid. 
This is explained in detail in Section \ref{subsec:lsm}.

In the first few lines of the code, we import the modules \f{dolfin}, \f{init}, \f{cm} and \f{pyplot} from \f{matplotlib},  \f{numpy}, \f{sys} and \f{os}.
The module \f{matplotlib} (\url{http://matplotlib.org/}) is used for plotting the design.
The module \f{dolfin} is a problem-solving environment required by FEniCS. 
The purpose of the line 
\begin{lstlisting}[firstnumber=9]
pp.switch_backend('Agg') 
\end{lstlisting}
is to use the \f{Agg} back end instead of the default \f{WebAgg} back end.
With the \f{Agg} back end, the figures do not appear on the screen, but are saved to a file; see lines  106-113.

\subsection{Initialization of case-dependent parameters}
In this section we describe the content of the file \f{init.py}, which provides initial data.
The outputs of \f{init.py} are the case-dependent variables, i.e. 
\f{Lag}, \f{Nx}, \f{Ny}, \f{lx}, \f{ly}, \f{Load}, \f{Name}, \f{ds}, \f{bcd}, \f{mesh}, \f{phi_mat}. 
The space \f{Vvec} is not case-dependent but is required to define the boundary conditions \f{bcd}.

The Lagrange multiplier $\Lambda$ for the volume constraint is called here \f{Lag}.
The variable \f{Load} is the position of the pointwise load, for example  \f{Load = [Point(lx, 0.5)]} for the cantilever, which  means that the load is applied at the point $(l_x,0.5)$. 
For the asymmetric cantilever we have \f{Load = [Point(lx, 0.0)]}. 

The fixed domain $\hold$ is a rectangle $\hold=[0,l_x]\times [0,l_y]$. 
In \f{init.py} this corresponds to the variables \f{lx,ly}. 
The mesh is built using the line 
\begin{lstlisting}[numbers=none,xleftmargin=0cm]
mesh  = RectangleMesh(Point(0.0,0.0),Point(lx,ly),Nx,Ny,'crossed')  
\end{lstlisting}
The class \f{RectangleMesh} creates a mesh in a 2D rectangle spanned by two points (opposing corners) of the rectangle. The arguments \f{Nx,Ny} specify the number of divisions in the $x$- and $y$-directions, and the optional argument \f{crossed} means that each square of the grid is divided in four triangles, defined by the crossing diagonals of the square. 
We choose \texttt{lx,ly,Nx,Ny} with the constraint  $\mbox{\texttt{lx Nx}}^{-1} = \mbox{\texttt{ly Ny}}^{-1}$.
The choice of the argument  \f{crossed} is necessary to have a symmetric displacement $u$ and in turn to keep a symmetric design throughout the iterations if the problem is symmetric, for instance in the case of the cantilever. 
Note that to preserve the symmetry of solutions at all time, one must choose an odd number of divisions \f{Nx} or \f{Ny}, depending on the orientation of the symmetry.
For instance, in the case of the symmetric cantilever, one can choose \f{Ny}$ = 75$ since the symmetry axis is the line $y = 1/2$, and \f{Nx}$ = 150$.

Since we chose a mesh with crossed diagonals, each square has an additional vertex at its center, where the diagonals meet. Therefore the total number of vertices is  
\begin{lstlisting}[firstnumber = 37]
dofsV_max =(Nx+1)*(Ny+1) + Nx*Ny 
\end{lstlisting}
We also define
\f{dofsVvec_max = 2*dofsV_max} 
in line 37, this represents the degrees of freedom for the vector function space \f{Vvec}. 

The case-dependent boundary $\Gamma_d$ is defined using the class \f{DirBd}, and instantiated by \f{dirBd = DirBd()} in \f{init.py}. 
We tag \f{dirBd} with the number $1$, the other boundaries with $0$, and introduce the boundary measure \f{ds}.
The Dirichlet boundary condition on $\Gamma_d$ is defined using  
\begin{lstlisting}[numbers=none,xleftmargin=0cm]
DirichletBC(Vvec,(0.0,0.0),boundaries,1) 
\end{lstlisting}
When several types of Dirichlet boundary conditions are required, as in the case of the half-wheel for instance,
the variable \f{bcd} is defined as a list of boundary conditions.
For the cantilever, \f{bcd} has only one element.
For the cases of the half-wheel and MBB-beam, we also define an additional class \f{DirBd2} to define the boundary conditions \f{bcd} because there are two different types of Dirichlet conditions; see Sections \ref{subsec:halfwheel} and \ref{subsec:mbb_beam}.

\subsection{Other initialization parameters}
The ersatz material coefficient $\epsilon$ is called \f{eps_er}; see  \eqref{eq:ersatz1}.
The elasticity parameters $E,\nu,\mu,\lambda$ are given lines 14-15.
In lines 17-19, a directory is created to save the results. 
In line 21, \f{ls_max = 3} is the maximum number of line searches for one iteration of the main loop, \f{ls} is an iteration counter for the line search, and the step size used in the gradient method is \f{beta},  initialized as \f{beta0_init}.
We choose \f{beta0_init = 0.5}.
We also choose \f{gamma = 0.8} and \f{gamma2 = 0.8}, which are used to modify the step size in the line search; see Section \ref{sec:ls}.

The counter \f{It} in line 25 keeps track of the iterations of the main loop. 
In line 25, we also fix a maximum number of iterations \f{ItMax = int(1.5*Nx)}, which depends on the mesh size \f{Nx} due to the fact that the time step \f{dt} is a decreasing function of \f{Nx}; see Section \ref{subsec:update_lsm}.

\subsection{Finite elements}
In line 28 and in the file \f{init.py}, we define the following finite element spaces associated with \f{mesh}:
\begin{lstlisting}[numbers=none,xleftmargin=0cm]
V = FunctionSpace(mesh, 'CG', 1)
Vvec = VectorFunctionSpace(mesh, 'CG', 1)
\end{lstlisting}
Here, \f{CG} is short for ``continuous Galerkin'', and the last argument is the degree of the element, meaning we have chosen the standard piecewise linear Lagrange elements. 
Note that the type of elements and degree can be easily modified using this command, and FEniCS offers a variety of them.
However, the level set part of our code has been written for this particular type of elements, so changing it would require to modify other parts of the code, such as the function \f{_comp_lsf} line 173, so one should be aware that it would not be a straightforward modification. 

\subsection{The function \f{_comp_lsf}}\label{subsec:lsm}

Here we explain the mechanism to get \f{phi} from \f{phi_mat}. 
Indeed \f{phi_mat} is updated every iteration by the function \f{_hj} in line 100, and we need \f{phi} to define the new set \f{Omega} in lines 42-44. 
Observe that the set of vertices of the Cartesian grid is included in the set of vertices of \f{mesh}, indeed the vertices of \f{mesh} are precisely the vertices of the Cartesian grid, plus the vertices in the center of the squares where the diagonals meet, due to the choice of the argument \f{crossed} in \f{mesh}. 
Thus we compute the values of \f{phi} at the center of the squares using interpolation.

This is done in the function \f{_comp_lsf} (lines 173-182) in the following way.
First of all, in lines 33-34, \f{dofsV} and \f{dofsVvec} are the coordinates of the vertices associated with the degrees of freedom. 
They are used in lines 35-36 to define \f{px}, \f{py} and \f{pxvec}, \f{pyvec}, which have integer values and are used by \f{_comp_lsf} to find the correspondence between the entries of the matrix \f{phi_mat} and the entries of \f{phi}. 
In  \f{_comp_lsf}, precisely line 175, we check if the vertex associated with \f{px}, \f{py} corresponds to a vertex on the Cartesian grid.
If this is the case, we set the values of \f{phi} in lines 176-177 to be equal to the values of \f{phi_mat} at the vertices which are common between \f{mesh} and the Cartesian grid.
Otherwise, the vertex is at the center of a square, and we set the value at this vertex to be the mean value of the four vertices of the surrounding square; see lines 179-181.

Thus the output of \f{_comp_lsf} is the function \f{phi} defined on \f{mesh}.
Note that if we had chosen squares with just one diagonal  (choosing \f{left} or \f{right} instead of \f{crossed} in \f{RectangleMesh} in the file \f{init.py}) instead of two, there would be an exact correspondence between \f{phi} and \f{phi_mat}, so that switching between the two would be easier. 

\subsection{Initialization of the level set function}\label{sec:init_lsf}
In lines 39-40, we initialize \f{phi} as a function in the space \f{V}, and using \f{_comp_lsf} we determine its entries using \f{phi_mat}.
The matrix \f{phi_mat} is initialized in the file \f{init.py}, since it is case-dependent.
For instance, for the cantilever we can choose
\begin{align} 
\label{init_phi}
\begin{split}
\phi(x,y) = & -\cos(8\pi x/l_x)\cos(4\pi y) - 0.4
+ \max(200(0.01 - x^2 - (y-l_y/2)^2),0)\\
& + \max(100(x +y -l_x-l_y+0.1),0)
+ \max(100(x -y -l_x+0.1),0),
\end{split}
\end{align}
which is the initialization yielding the result in Figure \ref{fig:init_cantilevers1}. 
The coefficients inside the cosine determine the initial number of ``holes'' inside the domain (i.e. the number of connected components of $\hold\setminus\Om$). 
Here \eqref{init_phi} corresponds to ten initial holes inside the domains (plus some half-holes on the boundary of $\hold$).

The reason for the additional three $\max$ terms in \eqref{init_phi} is specific of our approach.
Since $\VV\in\Theta^k(\hold)$, we have $\VV = 0$ in the corners of the rectangle $\hold$.
Therefore the shape in a small neighbourhood of the corners will not change, and if we start with an inappropriate initialization, we will end with a small set of unwanted material in certain corners.
Therefore the r\^ole of the $\max$-terms in \eqref{init_phi} is to create a small cut with the correct material in certain corners. 
Depending on the problem, it is easy to see what should be the correct corner material distribution for the final design.

Another problem may appear at boundary points which are on the symmetry axis for symmetric problems.
Indeed, due to the smoothness of $\VV$ and the symmetry of the problem, we will get $\VV_x = 0$ or $\VV_y = 0$  at these points and the shape will not change there.
For instance in the symmetric cantilever case, this problem happens at the point $[0,l_y/2]$. 
There is no issues at $[l_x,l_y/2]$, since this is the point where the load is applied, so it must be fixed anyway.
This explains the term $\max(200(0.01 - x^2 - (y-l_y/2)^2),0)$ in  \eqref{init_phi}.
In line 46 we define the integration measure \f{dX}  used line 50 to integrate on all of $\hold$.
In line 47 the normal vector \f{n} to $\hold$ is introduced to define the boundary conditions for \f{av} in line 51.

\subsection{Domain update}
The main loop starts line 55. 
In line 57 we instantiate by \f{omega = Omega()}.
This either initializes \f{omega} or updates \f{omega} if \f{phi} has been updated inside the loop.
In line 59, \f{omega} is tagged with the number $1$, and the complementary of \f{omega} is tagged with  the number $0$ in line 58.
We define the integration measure for the subdomains $\Om$ and $\hold\setminus\Om$ using
\begin{lstlisting}[firstnumber = 60]
dx = Measure('dx')(subdomain_data=domains)
\end{lstlisting}
One assembles using \f{dx(1)} to integrate on $\Om$, and using \f{dx(0)} to integrate on $\hold\setminus\Om$; see for instance line 68. 
For details on how to integrate on specific subdomains and boundaries, we refer to the FEniCS Tutorial \cite{langtangen2017solving} available at \url{http://www.springer.com/gp/book/9783319524610}  and the FEniCS documentation \cite{fenics:book}.
\subsection{Solving the elasticity system}
Then we can compute \f{U}, the solution of the elasticity system in lines 61-62, using \f{_solve_PDE} in lines 116-127.
We use a LU solver to solve the system; see lines 125-126.
The surface load is applied pointwise using the function \f{PointSource}; see lines 122-124.
Note that \f{U} is a list since we consider the general case of several loads.
Thus the length of the list \f{U} is the length of \f{Load}; see line 30.
\subsection{Cost functional update}
In lines 64-70 we compute the compliance, the volume of \f{omega} and the cost functional \f{J} corresponding to \eqref{Jnum}.
Observe that the command for the calculation of the compliance is close to the mathematical notation, i.e. it resembles the following mathematical formula:
\begin{align*} 
J(\Om) & = \e\int_{\hold\setminus\Om} \Sb_1 \ddo D\VV + \int_\Om \Sb_1 \ddo D\VV ,\qquad \Sb_1  = 2\mu e(u)\ddo e(u) + \lambda \tr(e(u))^2.
\end{align*}
Recall that the compliance is a sum in the case of several loads, see Section \ref{sec:multiple_load}, hence the \f{for} loop in  line 65.
\subsection{Line search and stopping criterion}\label{sec:ls}
The line search starts at line 72. 
If the criterion 
\begin{lstlisting}[firstnumber = 72]
J[It] > J[It-1] 
\end{lstlisting}
is satisfied, then we reject the current step.
In this case we reduce the step size \f{beta} by multiplying it by \f{gamma} in line 74.
Also, we go back to the previous values of  \f{phi_mat} and \f{phi} which were stored in \f{phi_mat_old} and \f{phi_old}, see line 75.
Then we need to recalculate  \f{phi_mat} and \f{phi} in lines 76-77 using the new step size \f{beta}.

If the step is not rejected, then we go to the next iteration starting from line 85.
If the step was accepted in the first iteration of the line search, in order to speed up the algorithm 
we increase the reference step size \f{beta0} by setting
\begin{lstlisting}[firstnumber = 86]
beta0 = min(beta0 / gamma2, 1) 
\end{lstlisting}
to take larger steps, since \f{gamma2} is smaller than one.
Here \f{beta0} is kept below 1 to stabilize the numerical scheme for the Hamilton-Jacobi equation, see the time step \f{dt} in line 150. 
We chose \f{gamma2 = 0.8} in our examples; see line 21.

If the maximum number of line searches \f{ls_max} is reached, we decrease in line 85 the reference step size \f{beta0} by setting
\begin{lstlisting}[firstnumber = 85]
beta0 = max(beta0 * gamma2, 0.1*beta0_init)
\end{lstlisting}
We impose the lower limit \f{0.1*beta0_init} on \f{beta0} so that the step size does not become too small.
Note that \f{beta} is reseted to \f{beta0} in line 88.

\subsection{Descent direction}\label{subsec:descent}
FEniCS uses the Unified Form Language (UFL) for representing weak formulations of partial differential equations, which results in an intuitive notation, close to the mathematical one. 
This can be seen in lines 49-51, where we define the matrix \f{av} which corresponds to the bilinear form \eqref{VP_2}. 
Here \f{theta} and \f{xi} are functions in \f{Vvec}, \f{xi} is the test function in \eqref{VP_1}, while \f{theta} corresponds to $\VV$ in \eqref{VP_1}.
In our code, we use the notation \f{th} for $\VV$ when $\VV$ is the descent direction.
The coefficient \f{1.0e4} in the boundary conditions for \f{av}: 
\begin{lstlisting}[firstnumber=51]
1.0e4*(inner(dot(theta,n),dot(xi,n)) *(ds(0)+ds(1)+ds(2)))
\end{lstlisting}
forces $\VV\cdot n$ to be close to zero on $\partial\hold$, which corresponds to the constraint $\VV\in \Theta^k(\hold)$.
For cases where \f{dirBd2} is not defined, such as the cantilever case, the term \f{+ds(2)} has no effect.

We assemble the matrix for the PDE of \f{th} and define the LU solver in lines 50-53, before the start of the main loop. 
Indeed, the bilinear form $\mathcal B$ in \eqref{VP_2} is independent of $\Om$.
Thus we reuse the factorization of the LU solver to solve the PDE for \f{th} using the parameter \f{reuse_factorization} in line 53.
This allows to spare some calculations, but for grids larger than the ones considered in this paper, it would be appropriate to use more efficient approaches such as Krylov methods to solve the PDE. 
This can be done easily with FEniCS using one of the various available solvers.

In line 90, we compute the descent direction \f{th}.
The function \f{_shape_der} in lines 129-139 solves the PDE for $\VV$, i.e. it implements \eqref{VP_1}-\eqref{VP_2} using the volume expression of the shape derivative \eqref{shape_der_ersatz0} and \eqref{Sb1_ersatz_compliance}. 
The variational formulation used in our code for the case of one load is: find $\VV\in H^k_d(\hold)^m$   such that  
\begin{align} \label{eqth}
\int_\hold \alpha_1 D\VV\ddo D\xi +\alpha_2 \VV\cdot\xi
&= -d\mathcal{J}^{\tt{vol}}(\Om,\xi), \forall \xi\in H^1(\hold)^m
\end{align}
$$\mbox{ where }\qquad d\mathcal{J}^{\tt{vol}}(\Om,\xi) = \int_\hold (2Du^\transp A_\Om e(u)-  A_\Om e(u)\ddo e(u) I_d) : D\xi
 + \Lambda\int_{\Om} \divv{\xi},   $$
and with $\alpha_1 = 1$ and $\alpha_2 = 0.1$.
When several loads are applied, as in the case of \f{cantilever_twoforces}, the right-hand side in \eqref{eqth} should be replaced by a sum over the loads \f{u} in \f{u_vec}; see Section \ref{sec:multiple_load}.

The right-hand side of \eqref{eqth} is assembled in line 136, but we need to integrate separately on $\Om$ and $\hold\setminus\Om$ using \f{dx(1)} and \f{dx(0)}, respectively.
The system is solved line 138 using the solver defined line 52.

In line 90 we get the descent direction \f{th} in the space \f{Vvec}. 
To update \f{phi} we need \f{th} on the Cartesian grid. 
As we explained already, we just need to extract the appropriate values of \f{th} since the Cartesian grid is included in \f{mesh}. 
This is done in lines 91-97, and the corresponding function on the Cartesian grid is called \f{th_mat}.

\subsection{Update the level set function}\label{subsec:update_lsm}

Then, we proceed to update the level set function \f{phi_mat} using the subfunction \f{_hj}.
The subfunction \f{_hj} in lines 141-152 follows exactly the discretization procedure described in Section \ref{sec:discretization_hj}.
In lines 143-146, the quantities \f{Dxm}, \f{Dxp}, \f{Dyp}, \f{Dym} correspond to $p^-,p^+,q^+,q^-$, respectively.
In line 142, we take $10$ steps of the Hamilton-Jacobi update, which is a standard, although heuristic way, to accelerate the convergence.
In order to stabilize the numerical scheme, we choose the time step as
\begin{lstlisting}[firstnumber = 150]
dt  = beta*lx / (Nx*maxv) 
\end{lstlisting}
where \f{maxv} is equal to $\max_{\Om}(|\VV_1| + |\VV_2|)$, \f{lx/Nx} is the cell size, and the step size \f{beta} is smaller than $1$ at all time in view of line 86.
In line 99, we save the current versions of \f{phi} and \f{phi_mat} in the variables \f{phi_old}, \f{phi_mat_old}, for use in case the step gets rejected during the line search. 
In line 102, the function \f{phi} is extrapolated from \f{phi_mat} using \f{_comp_lsf}.

\subsection{Reinitialization of the  level set function}

Every $5$ iterations, we reinitialize the level set function in line 101.
This is achieved by the subfunction \f{_reinit} in lines 154-171.
The reinitialization follows the procedure described in Section \ref{reinit}.
In lines 161-164, \f{Dxm}, \f{Dxp}, \f{Dyp}, \f{Dym} correspond to $p^-,p^+,q^+,q^-$, respectively.
In lines 165 to 168, \f{Kp} and \f{Km} correspond to $K^+$ and $K^-$ from \eqref{Gplus}-\eqref{Gminus}, respectively.
Line 169 corresponds to \eqref{computeG} and line 170 to the update \eqref{reinit_phi_ij}.

The  function \f{signum}, computed in line 159, is the approximation of the sign function of $\phi$ corresponding to $S(\phi)$, defined in \eqref{signum}.
To compute \f{signum}, we use \f{lx/Nx} for $\e_s$, and $|\nabla \phi|$ is computed using symmetric finite differences $\phi$, which are given by \f{Dxs} and \f{Dys} in lines 155-158. 
\subsection{Stopping criterion and saving figures}
Finally in lines 104-105, we check if the stopping criterion
\begin{lstlisting}[firstnumber = 104]
if It>20 and max(abs(J[It-5:It]-J[It-1])) <2.0*J[It-1]/Nx**2: 
\end{lstlisting}
is satisfied.
Here, \f{It-1} is the current iteration.
This means that the algorithm stops when the maximum difference of the value of the cost functional at the current iteration with the values of the four previous iterations is below a certain threshold. 
In order to take smaller steps when the grid gets finer, we have determined heuristically the threshold \f{2.0*J[It-1]/Nx**2} which depends on the grid size \f{Nx}. 

Lines 107-113 are devoted to plotting the design. 
The filled contour of the zero level set of \f{phi_mat} is drawn using the \f{pyplot} function \f{contourf}; see the matplotlib documentation \url{http://matplotlib.org/} for details.

\section{Numerical results and case-dependent parameters}\label{sec:num_res}
In this section we discuss the case-dependent parameters in \f{init.py} such as boundary conditions,  load position, and Lagrangian $\Lambda$.

\subsection{Symmetric cantilever}
\begin{figure}
\begin{center}
\begin{subfigure}[b]{0.3\textwidth}
\includegraphics[width=\textwidth]{cantilever_init2_lagV=40_Nx=150_final_design.eps}
\caption{$(Nx,Ny) = (150,75)$.}
\label{fig:cantilever1}
\end{subfigure}
\begin{subfigure}[b]{0.3\textwidth}
\includegraphics[width=\textwidth]{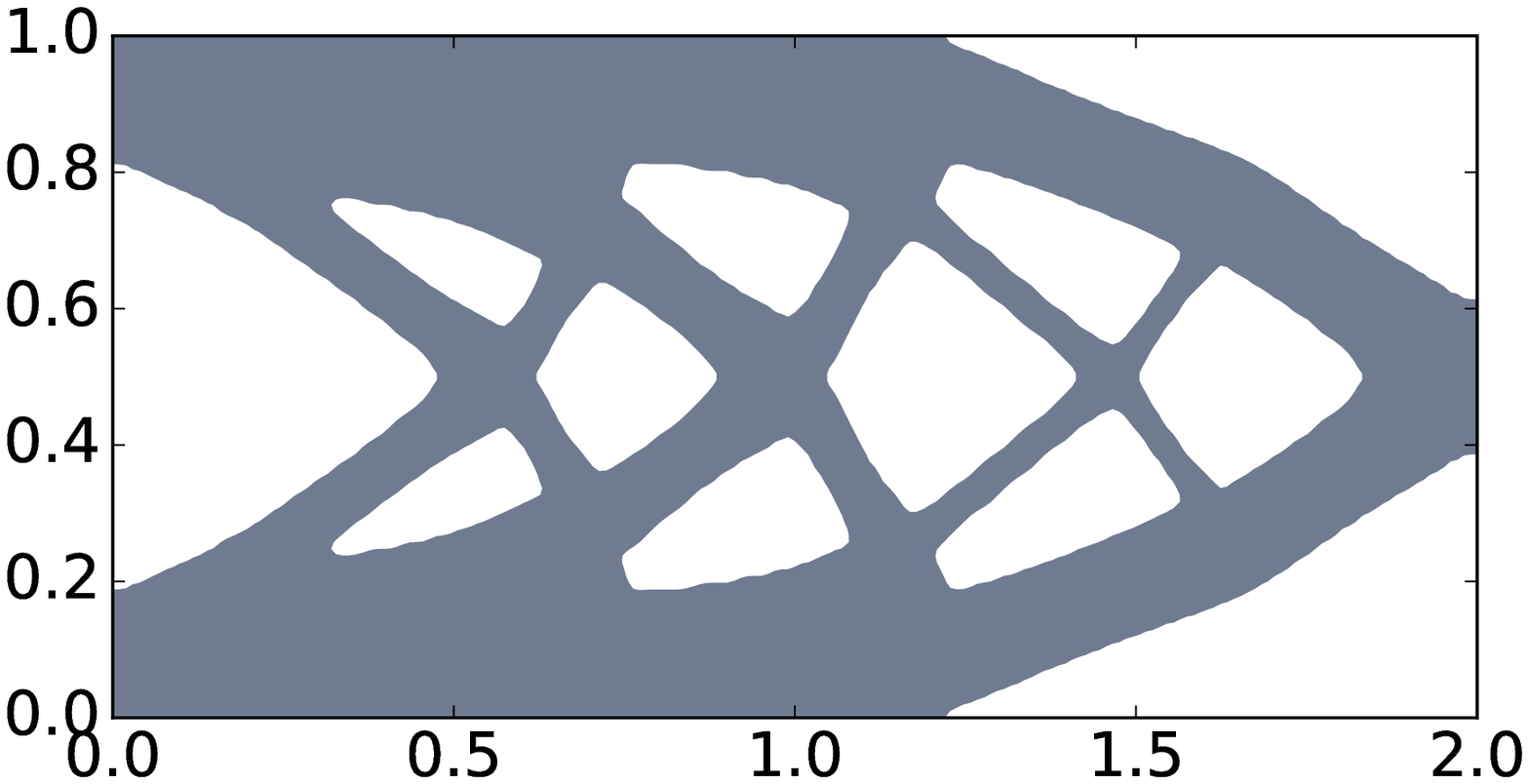}
\caption{$(Nx,Ny) = (202,101)$.}
\label{fig:cantilever2}
\end{subfigure}
\begin{subfigure}[b]{0.3\textwidth}
\includegraphics[width=\textwidth]{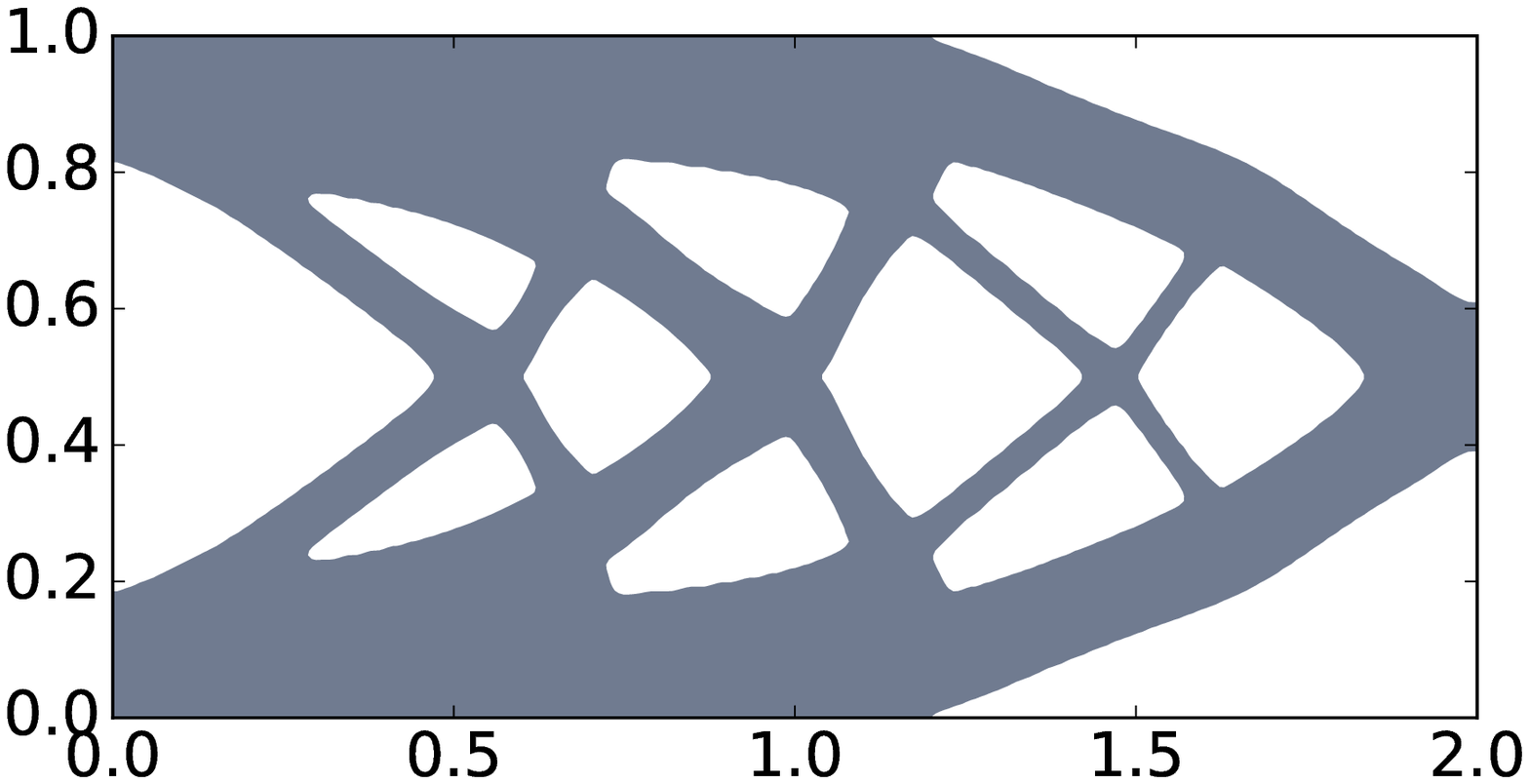}
\caption{$(Nx,Ny) = (302,151)$.}
\label{fig:cantilever3} 
\end{subfigure}
\caption{Optimal design for the symmetric cantilever,  with $\Lambda = 40$, and initialization \eqref{init_phi}.}
\label{fig:init_cantilevers1}
\end{center} 
\end{figure}

For the symmetric cantilever the load is placed at the point $(l_x,l_y/2)$, see the parameter \f{Load}.
The initialization for the symmetric cantilever is given by \eqref{init_phi}.
See Figure \ref{fig:init_cantilevers1} for the results of the symmetric cantilever for several grid sizes and $\Lambda = 40$. See also Figure \ref{fig:cantilever0} for a comparison of two different initializations.
We observe that the optimal set is independent of the mesh size, but depends on the initialization.

To obtain a short symmetric cantilever, one can set \f{lx = ly} and \f{Nx = Ny}. 
Still, one should choose an odd number for \f{Nx} in order to preserve the symmetry of the problem. 

\subsection{Asymmetric cantilever}
\begin{figure}
\begin{center}
\begin{subfigure}[b]{0.3\textwidth}
\includegraphics[width=\textwidth]{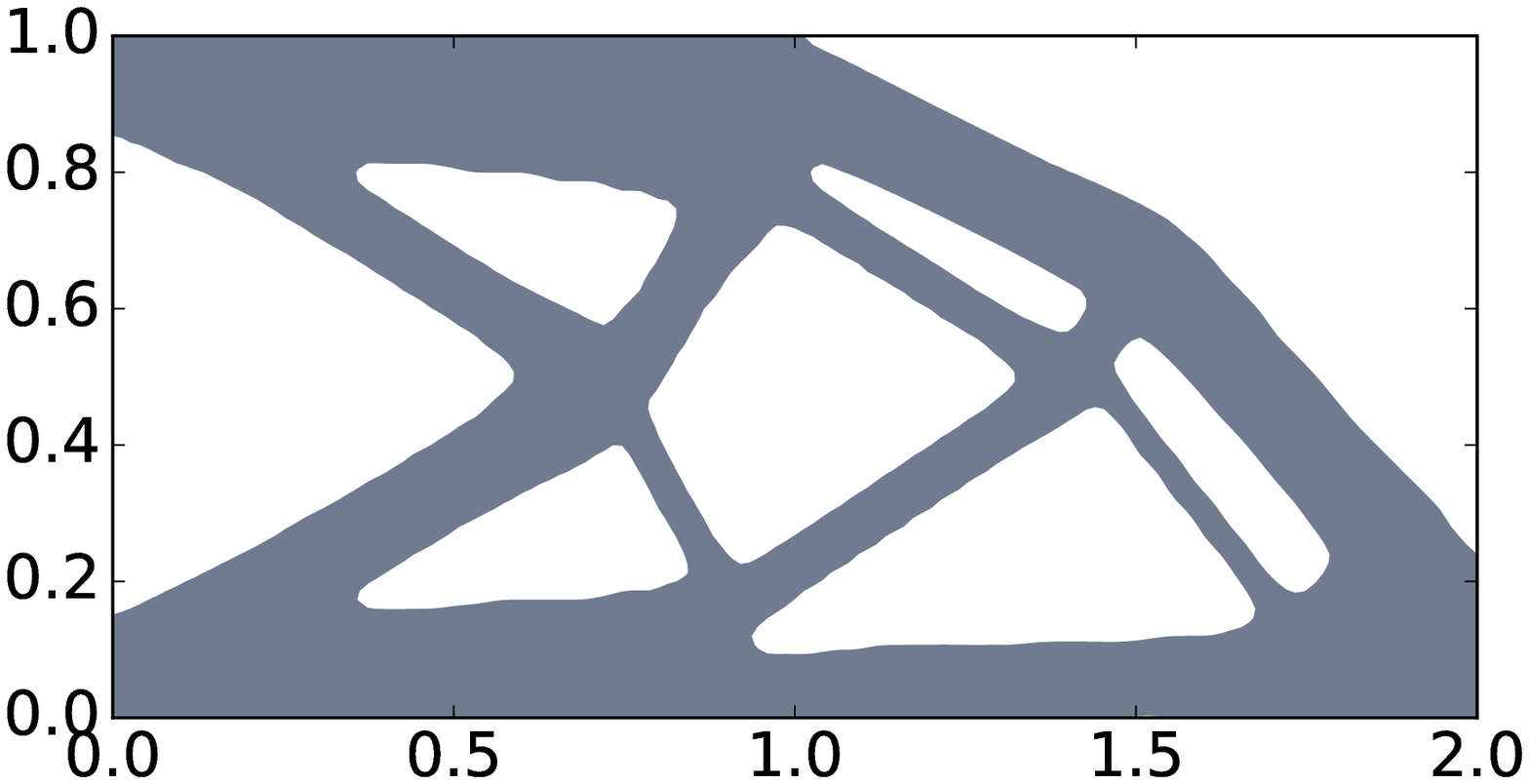}
\caption{$(Nx,Ny) = (150,75)$.}
\label{fig:cantilever_assym1}
\end{subfigure}
\begin{subfigure}[b]{0.3\textwidth}
\includegraphics[width=\textwidth]{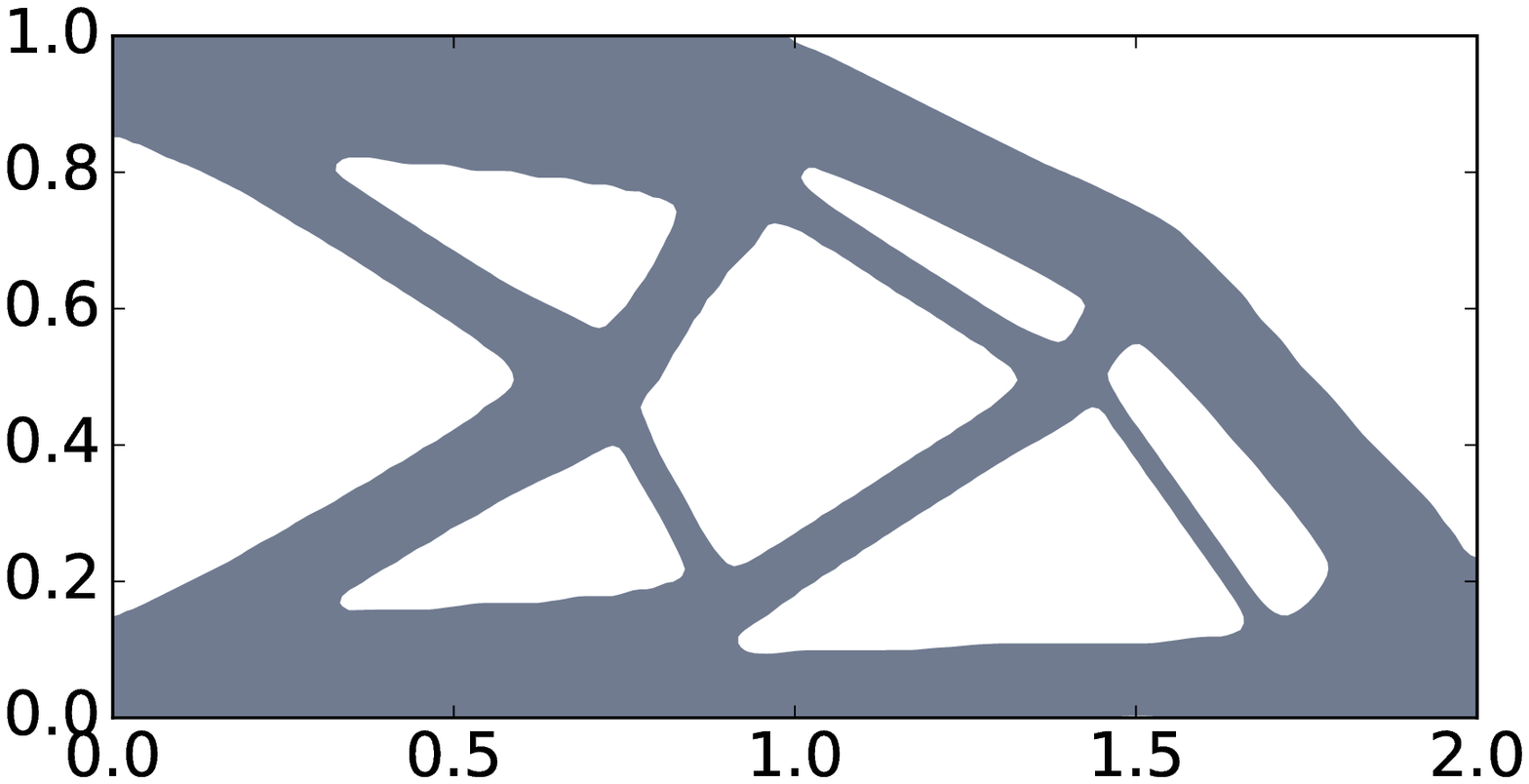}
\caption{$(Nx,Ny) = (202,101)$.}
\label{fig:cantilever_assym2}
\end{subfigure}
\begin{subfigure}[b]{0.3\textwidth}
\includegraphics[width=\textwidth]{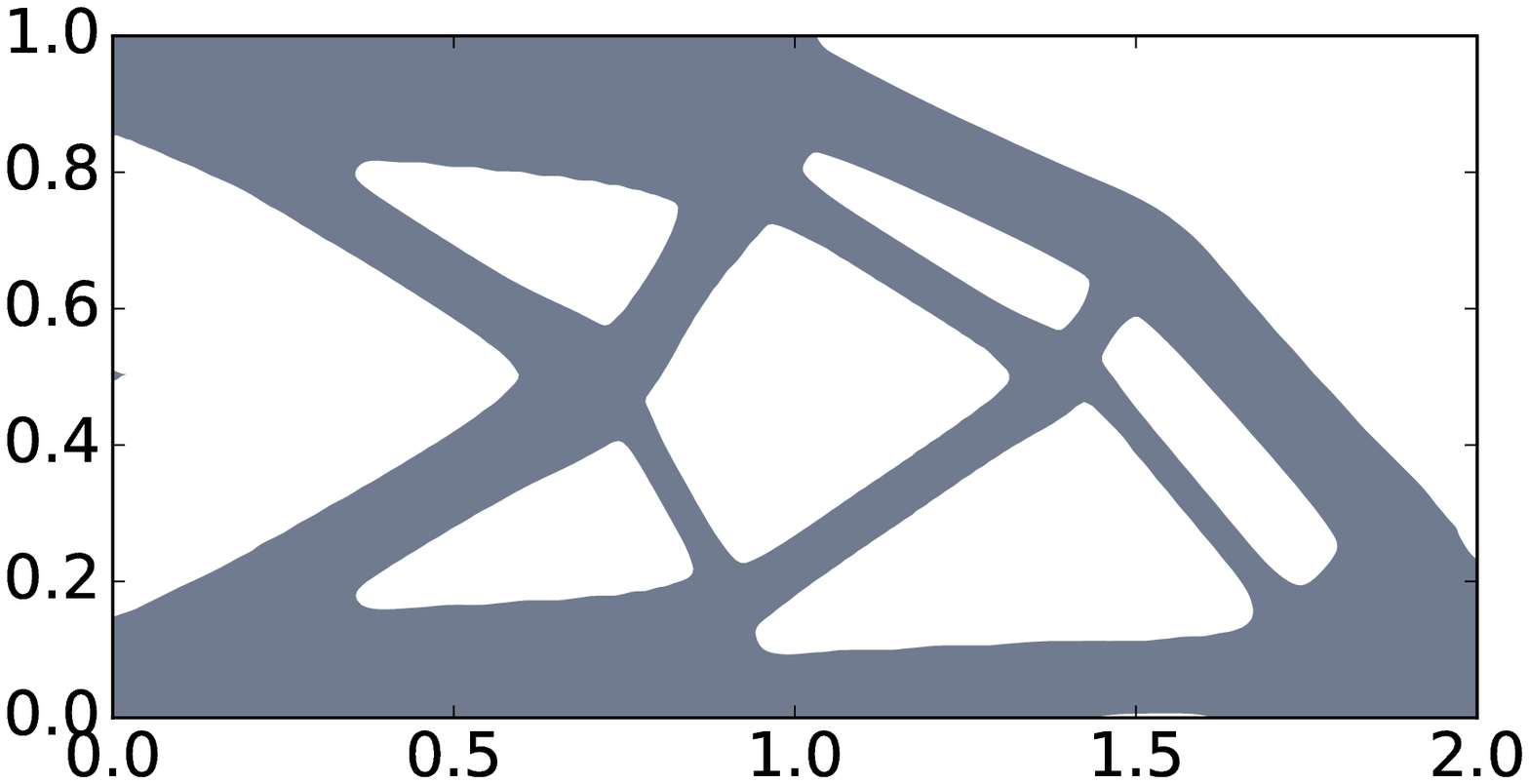}
\caption{$(Nx,Ny) = (302,151)$.}
\label{fig:cantilever_assym3}
\end{subfigure}
\begin{subfigure}[b]{0.3\textwidth}
\includegraphics[width=\textwidth]{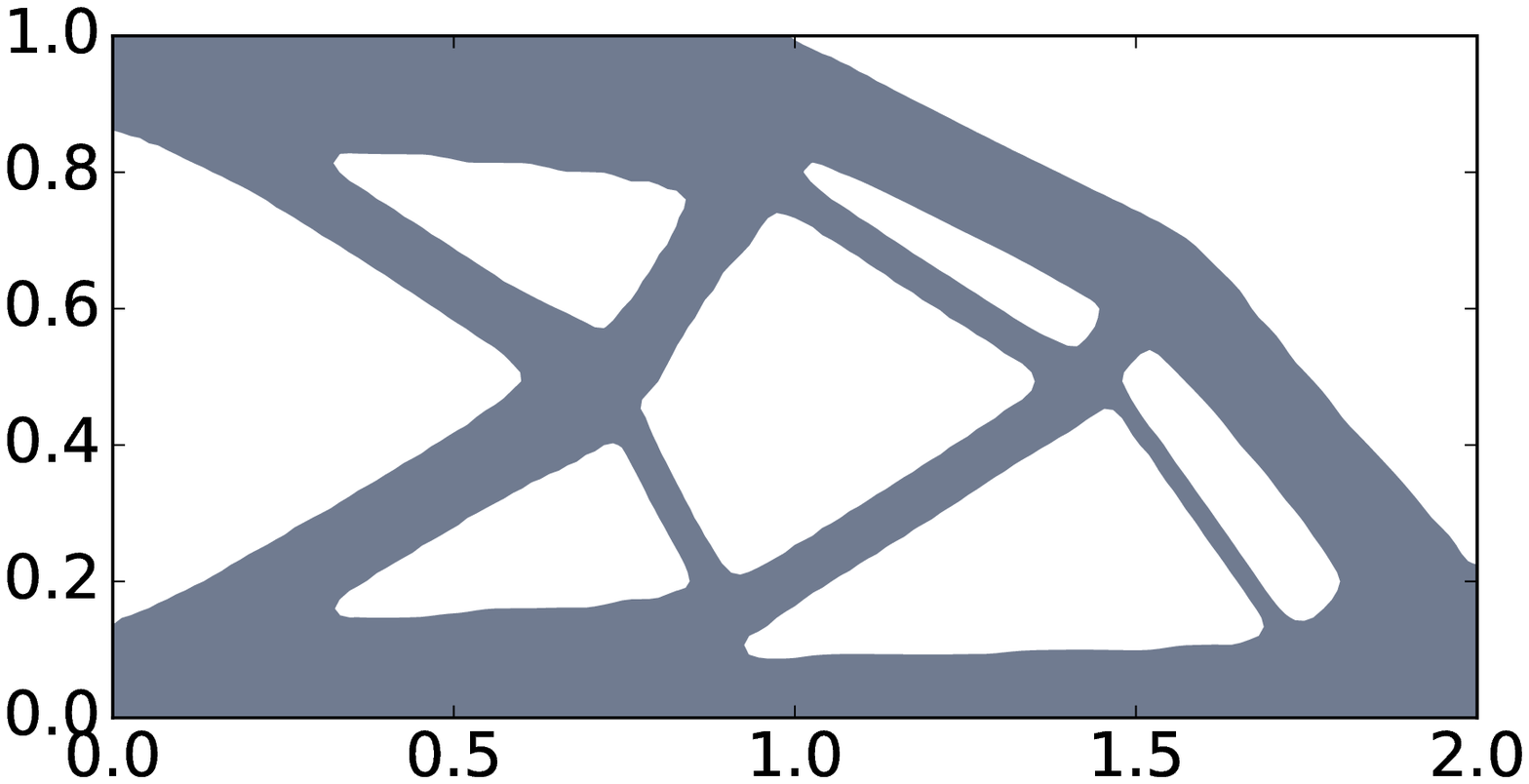}
\caption{$(Nx,Ny) = (150,75)$.}
\label{fig:cantilever_assym1_la70}
\end{subfigure}
\begin{subfigure}[b]{0.3\textwidth}
\includegraphics[width=\textwidth]{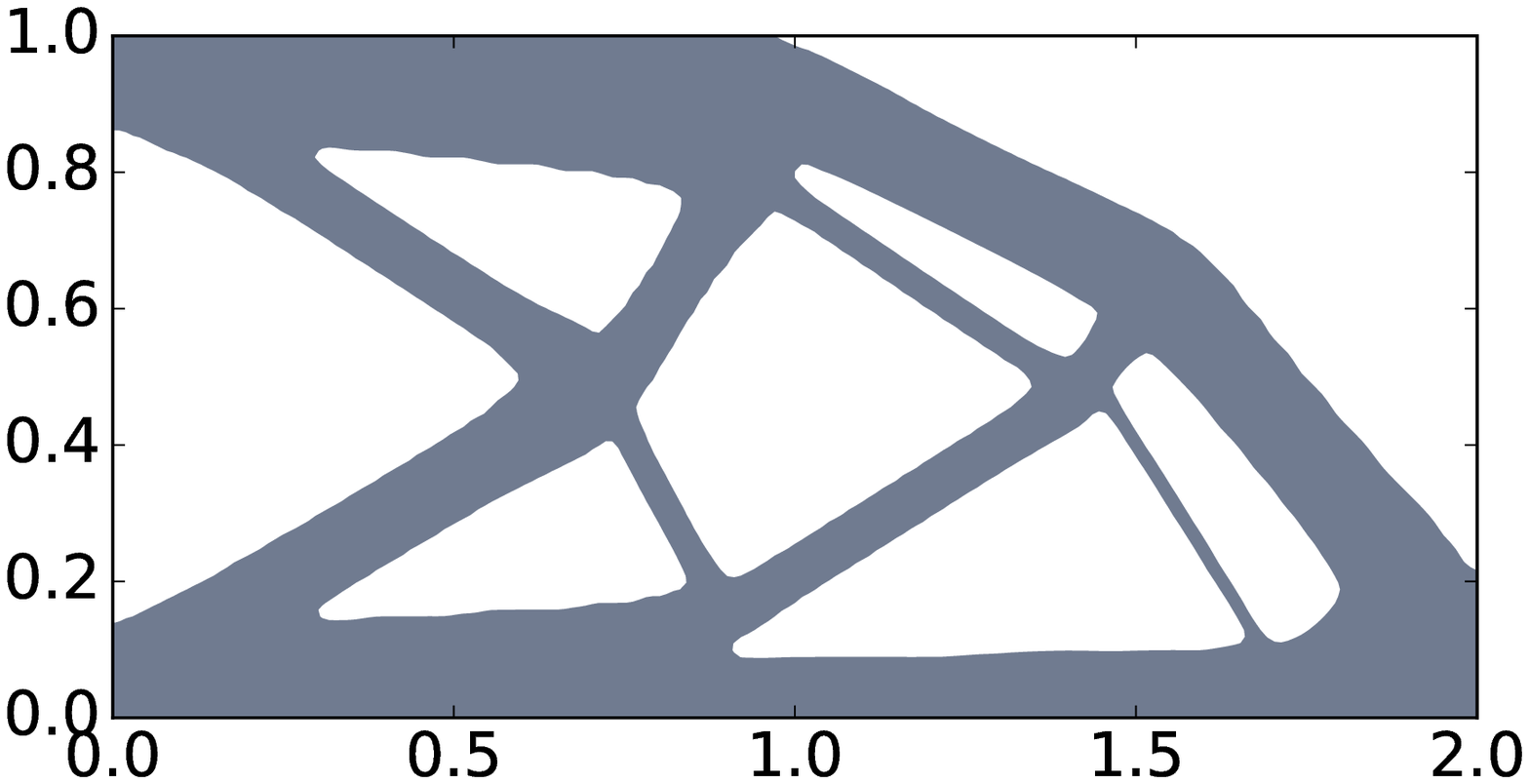}
\caption{$(Nx,Ny) = (202,101)$.}
\label{fig:cantilever_assym2_la70}
\end{subfigure}
\begin{subfigure}[b]{0.3\textwidth}
\includegraphics[width=\textwidth]{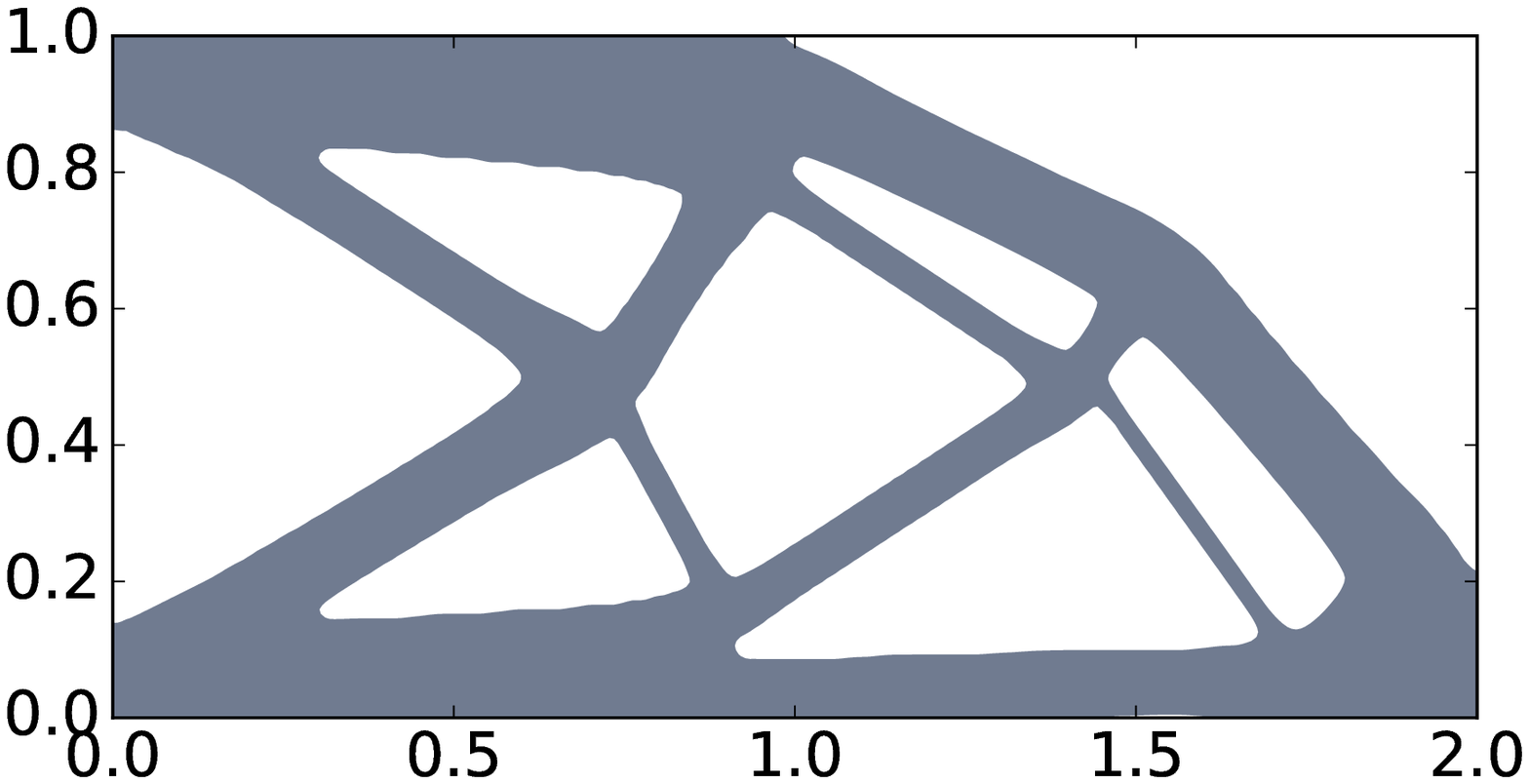}
\caption{$(Nx,Ny) = (302,151)$.}
\label{fig:cantilever_assym3_la70}
\end{subfigure}
\caption{Optimal design for the asymmetric cantilever, $\Lambda = 60$ (first row),  $\Lambda = 70$ (second row). }
\label{fig:cantilever_assym0} 
\end{center} 
\end{figure}

We take \f{Load = [Point(lx, 0.0)]} for the asymmetric cantilever, to have a load in the lower right corner.
The initialization is also changed, so as to start with the material phase where the load is applied, more precisely, we choose
\begin{align*}
\phi(x,y) =&  -\cos(6\pi x/l_x) \cos(4 \pi y) - 0.4
+ \max(100(x+y-l_x-l_y+0.1),0).  
\end{align*}
See Figure \ref{fig:cantilever_assym0} for the results of the asymmetric cantilever for several grid sizes, for $\Lambda = 60$ and $\Lambda = 70$.

\subsection{Half-wheel}\label{subsec:halfwheel}
For the half-wheel we have \f{lx,ly = [2.0,1.0]}. 
The position of the load is given by 
\f{Load = [Point(lx/2, 0.0)]}. 
In the corner $(0.0,0.0)$ we need pointwise Dirichlet conditions  and rolling conditions in $(l_x,0.0)$. 
For this we use the following boundaries in \f{init_py}:
\begin{lstlisting}[numbers=none,xleftmargin=0cm]
class DirBd(SubDomain):
    def inside(self, x, on_boundary):
        return abs(x[0])< tol and abs(x[1])< tol 
class DirBd2(SubDomain):
    def inside(self, x, on_boundary):
        return abs(x[0]-lx)<tol and abs(x[1])<tol
dirBd,dirBd2 = [DirBd(),DirBd2()]     
\end{lstlisting}
where \f{tol = 1E-14}.
Then the two boundary parts are tagged with different numbers 
\begin{lstlisting}[numbers=none,xleftmargin=0cm]
dirBd.mark(boundaries, 1)
dirBd2.mark(boundaries, 2),  
\end{lstlisting}
and we define the vector of boundary conditions as
\begin{lstlisting}[numbers=none,xleftmargin=0cm]
bcd  = [DirichletBC(Vvec, (0.0,0.0), dirBd, method='pointwise'),\
        DirichletBC(Vvec.sub(1), 0.0, dirBd2,method='pointwise')] 
\end{lstlisting}
The method \f{pointwise} is used since \f{dirBd2} is a single point.
Note here that the rolling boundary condition is achieved by setting the component \f{Vvec.sub(1)} to $0$, indeed \f{Vvec.sub(1)} represents the $y$-component of a vector function taken in the space \f{Vvec} .
In lines 50-51 of \f{compliance.py}, approximate Dirichlet conditions for $\VV$  are applied on \f{dirBd}  and \f{dirBd2} as these corners should be fixed.

The initialization should also change to fit the half-wheel case. 
We chose
\begin{lstlisting}[numbers=none,xleftmargin=0cm]
phi_mat = -np.cos((3.0*pi*(XX-1.0))) * np.cos(7*pi*YY) - 0.3 
+ np.minimum(5.0/ly *(YY-1.0) + 4.0,0) \
+ np.maximum(100.0*(XX+YY-lx-ly+0.1),.0) + np.maximum(100.0*(-XX+YY-ly+0.1),.0)  
\end{lstlisting}

\begin{figure}
\begin{center}
\begin{subfigure}[b]{0.3\textwidth}
\includegraphics[width=\textwidth]{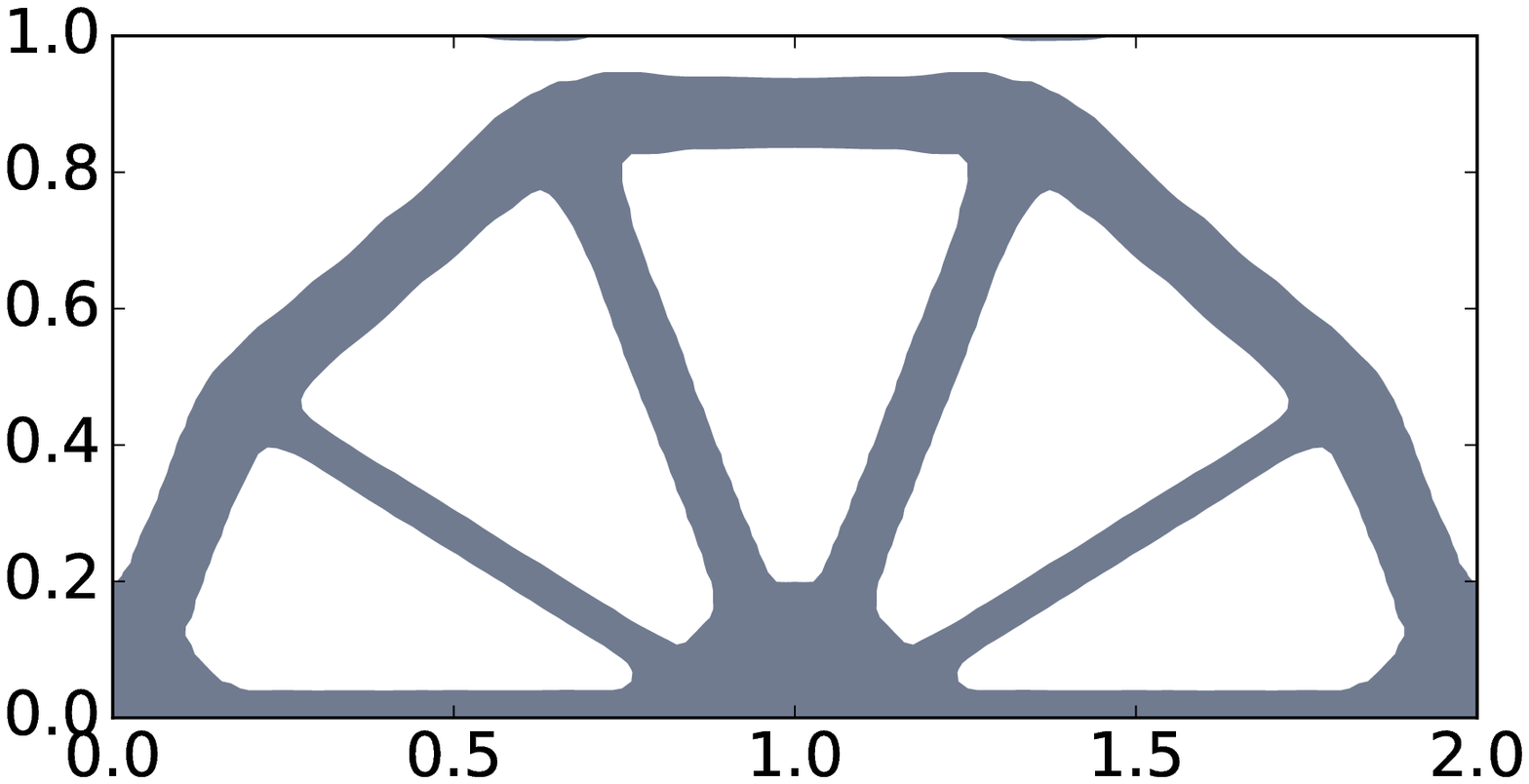}
\caption{$(Nx,Ny) = (150,75)$.}
\label{fig:wheel1}
\end{subfigure}
\begin{subfigure}[b]{0.3\textwidth}
\includegraphics[width=\textwidth]{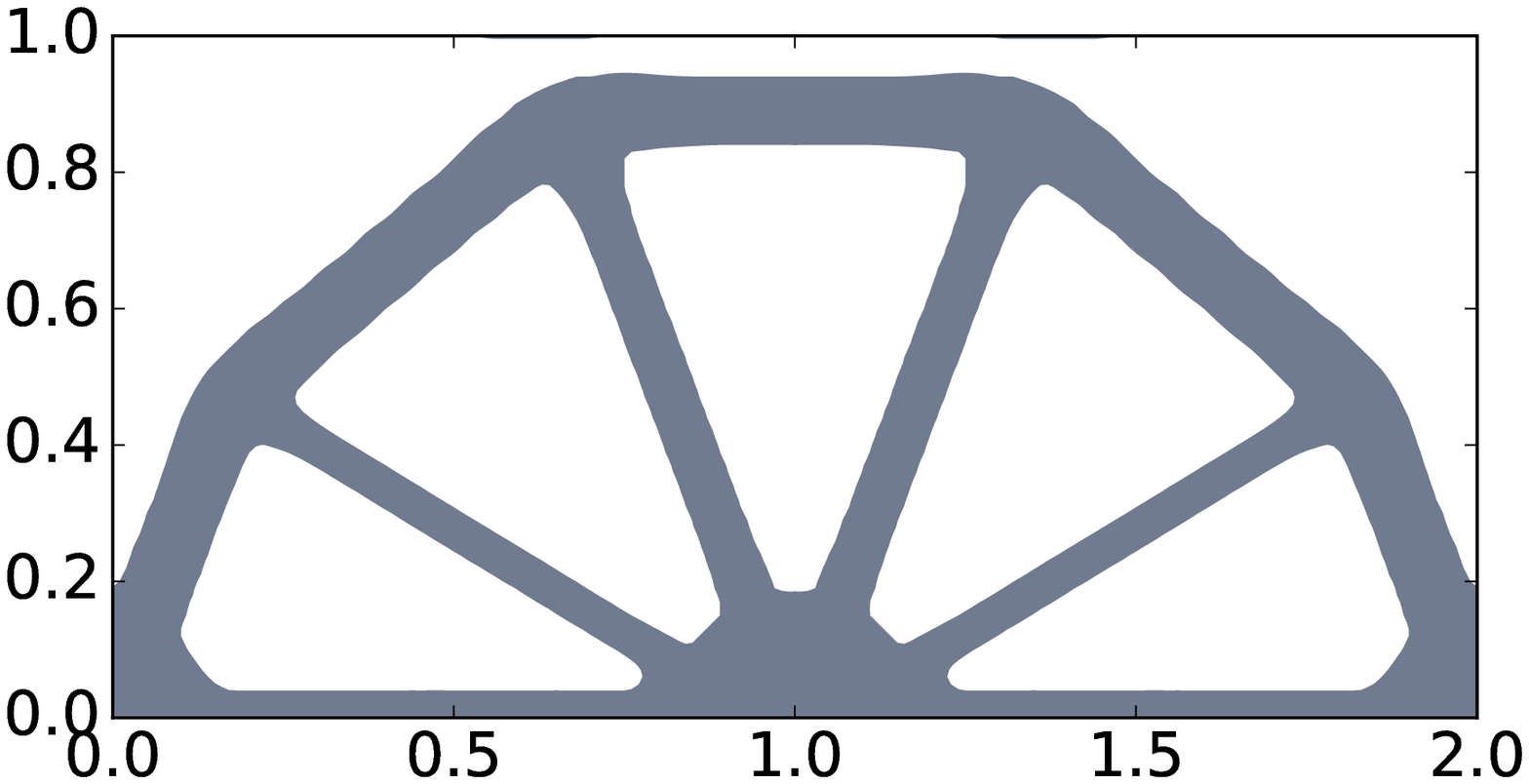}
\caption{$(Nx,Ny) = (200,100)$.}
\label{fig:wheel2}
\end{subfigure}
\begin{subfigure}[b]{0.3\textwidth}
\includegraphics[width=\textwidth]{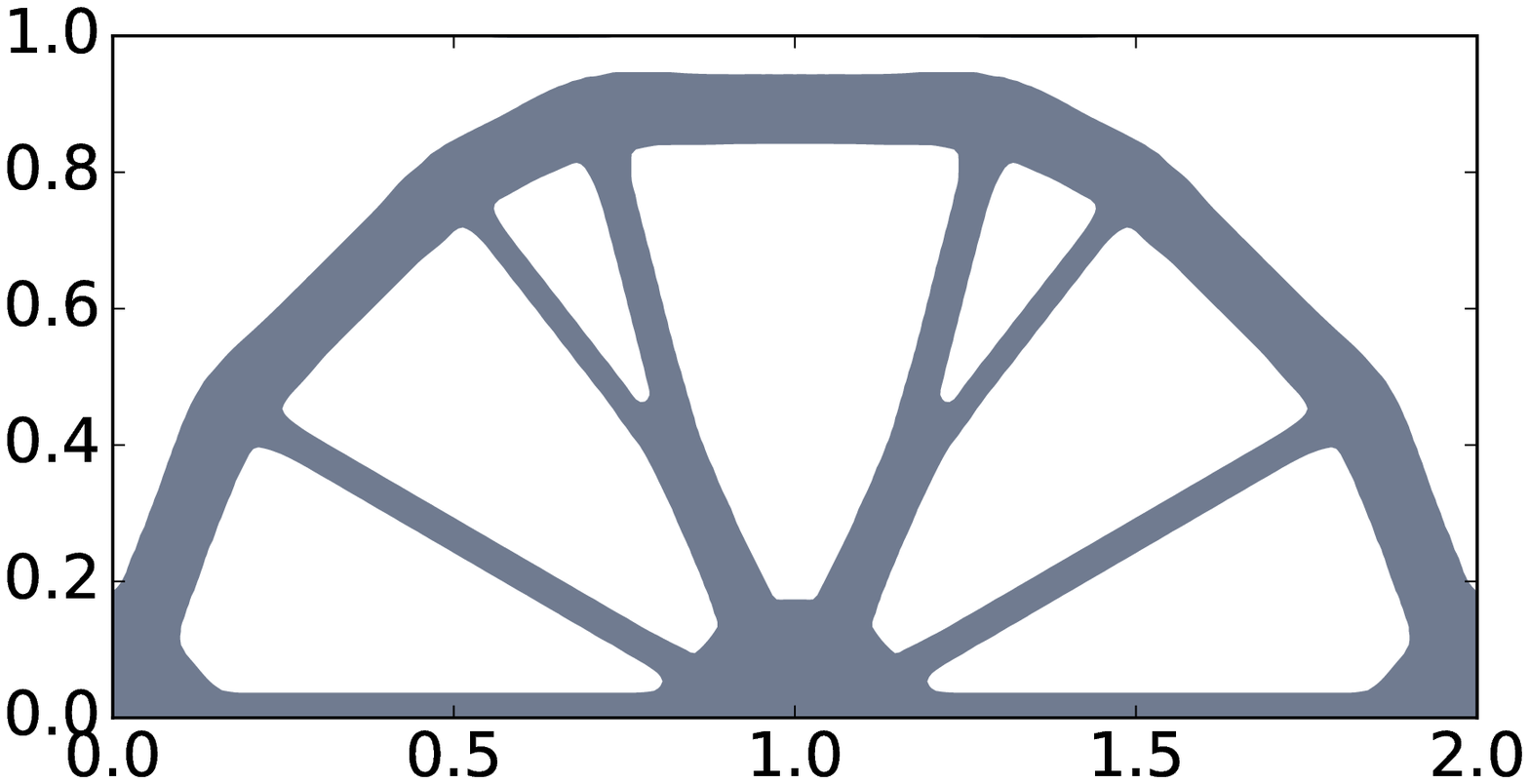}
\caption{$(Nx,Ny) = (300,150)$.}
\label{fig:wheel3}
\end{subfigure}\\
\begin{subfigure}[b]{0.3\textwidth}
\includegraphics[width=\textwidth]{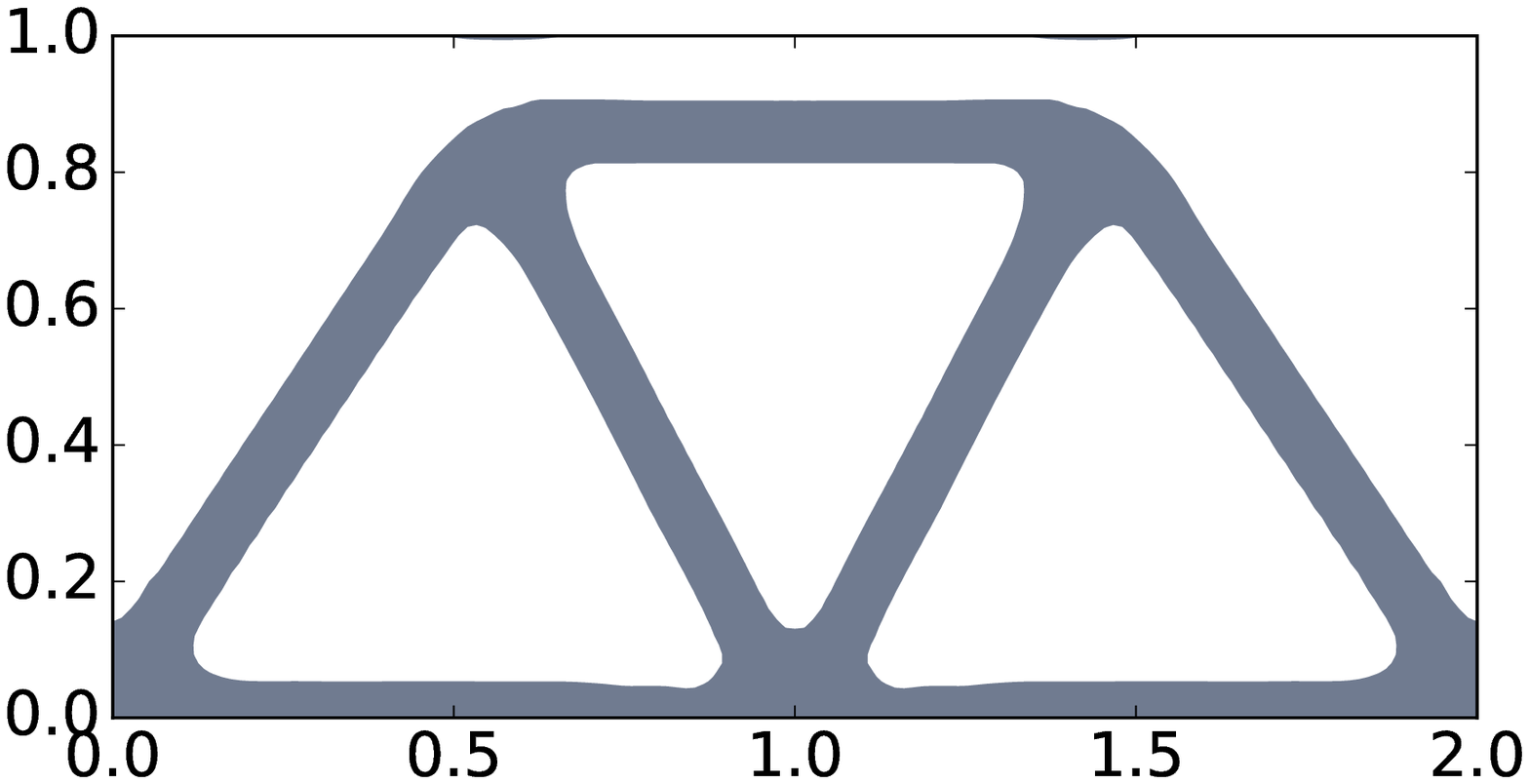}
\caption{$(Nx,Ny) = (150,75)$.}
\label{fig:wheel1_la50}
\end{subfigure}
\begin{subfigure}[b]{0.3\textwidth}
\includegraphics[width=\textwidth]{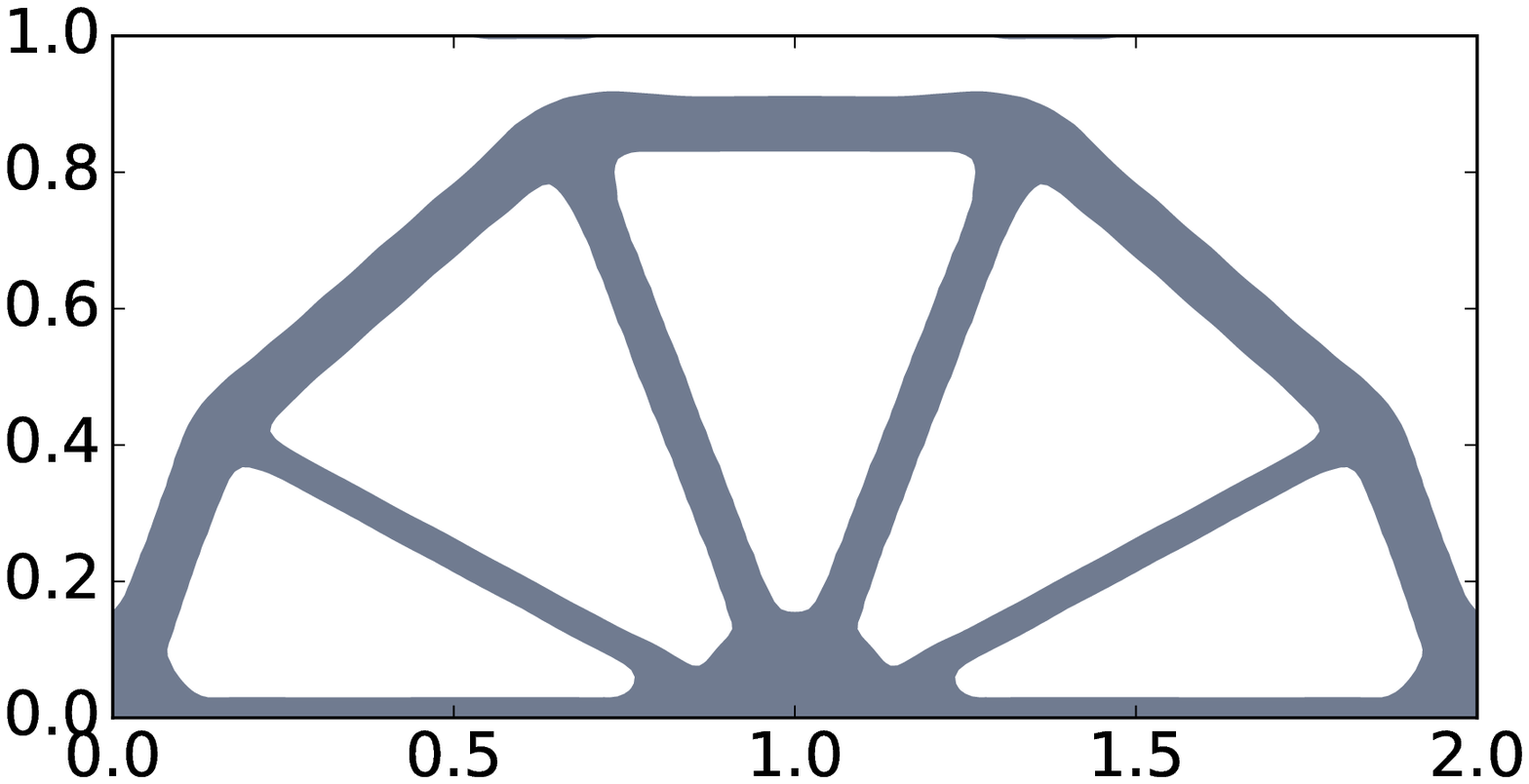}
\caption{$(Nx,Ny) = (200,100)$.}
\label{fig:wheel2_la50}
\end{subfigure}
\begin{subfigure}[b]{0.3\textwidth}
\includegraphics[width=\textwidth]{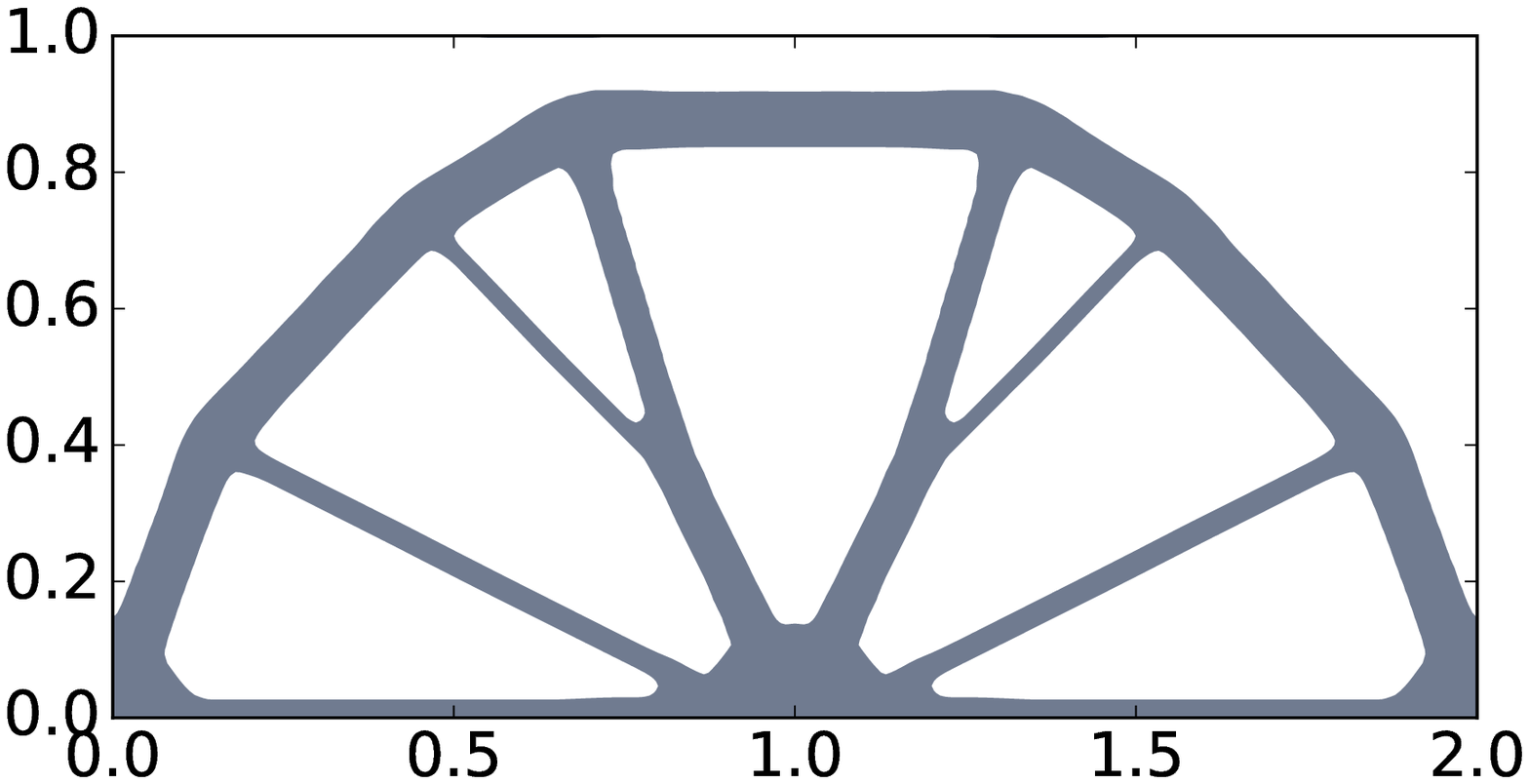}
\caption{$(Nx,Ny) = (300,150)$.}
\label{fig:wheel3_la50}
\end{subfigure}
\caption{Optimal design for the half-wheel,  $\Lambda = 30$ (first row),  $\Lambda = 50$ (second row).}
\label{fig:wheel0} 
\end{center} 
\end{figure}

In Figure \ref{fig:wheel0} we compare results obtained with $\Lambda = 30$ and $\Lambda = 50$.
\subsection{Bridge}
\begin{figure}
\begin{center}
\begin{subfigure}[b]{0.3\textwidth}
\includegraphics[width=\textwidth]{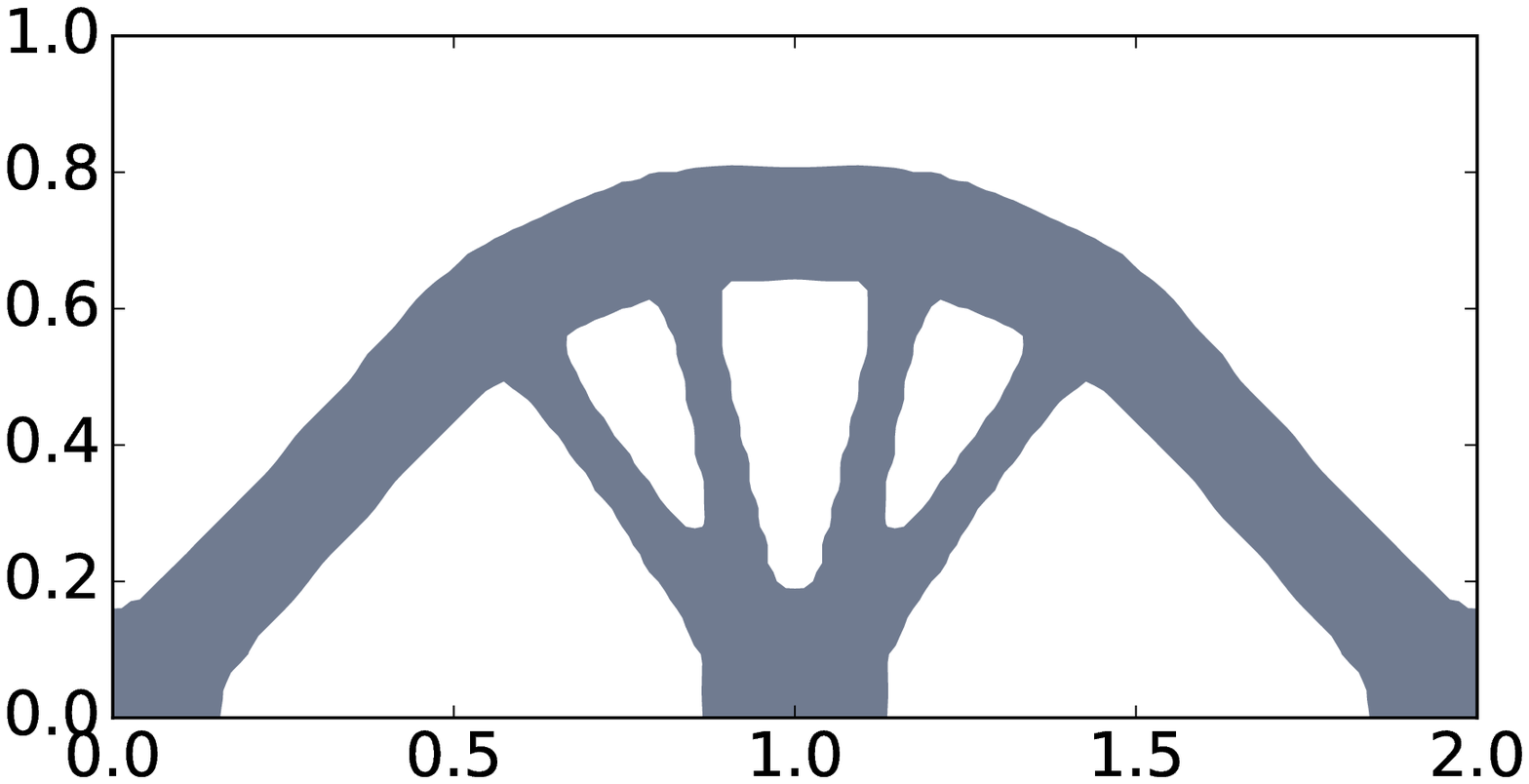}
\caption{$(Nx,Ny) = (150,75)$.}
\label{fig:brdige1}
\end{subfigure}
\begin{subfigure}[b]{0.3\textwidth}
\includegraphics[width=\textwidth]{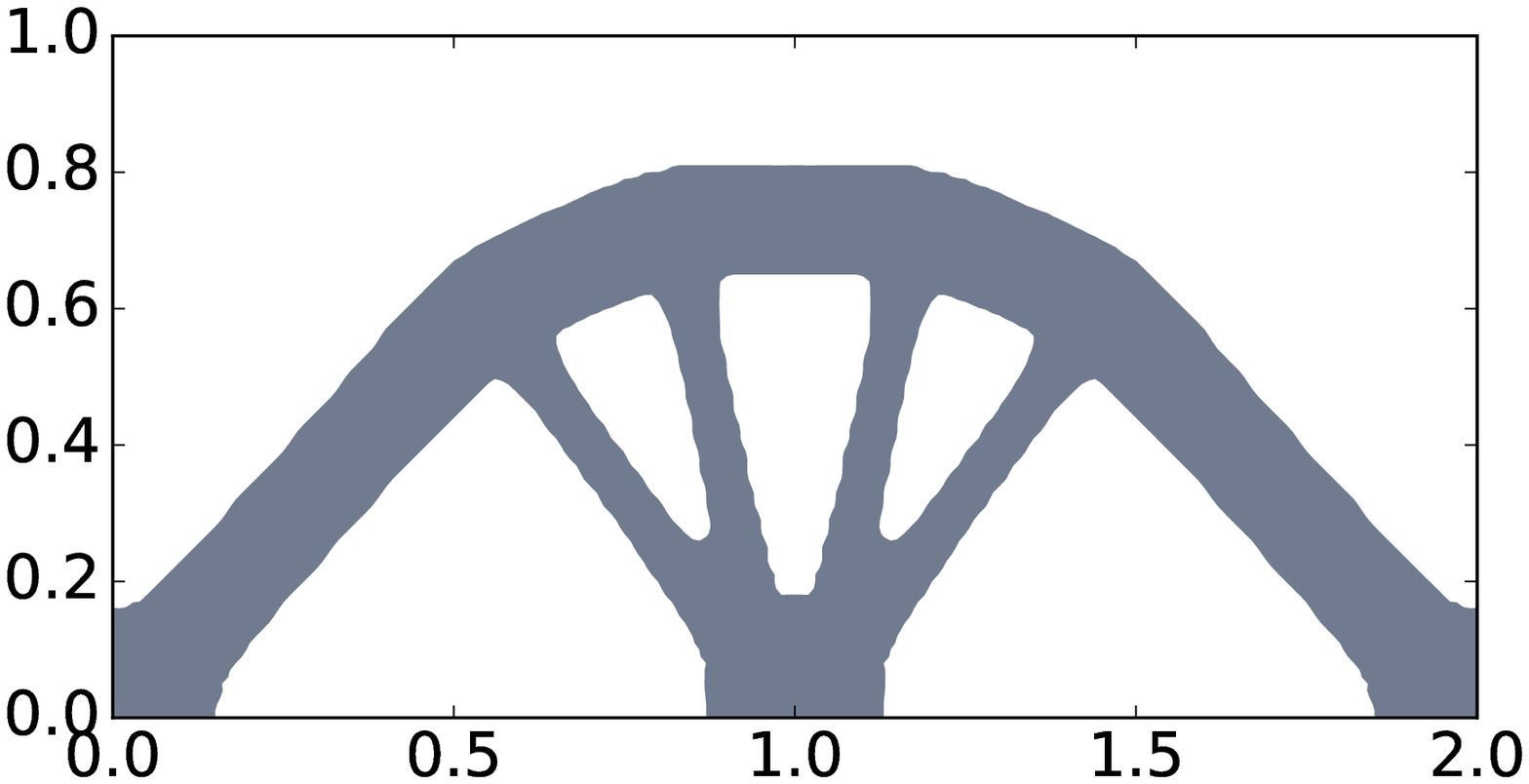}
\caption{$(Nx,Ny) = (200,100)$.}
\label{fig:bridge2}
\end{subfigure}
\begin{subfigure}[b]{0.3\textwidth}
\includegraphics[width=\textwidth]{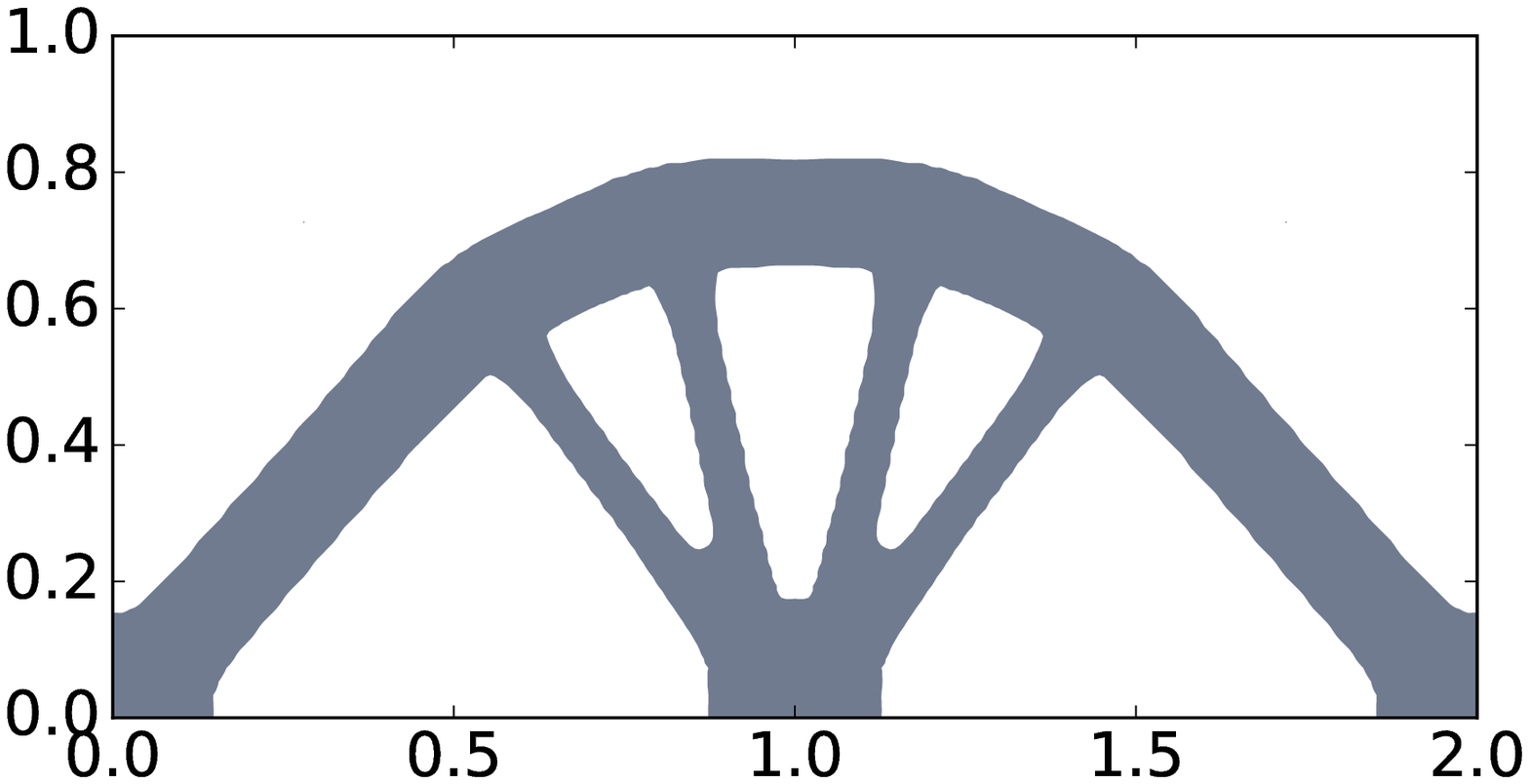}
\caption{$(Nx,Ny) = (300,150)$.}
\label{fig:bridge3} 
\end{subfigure}\\
\begin{subfigure}[b]{0.3\textwidth}
\includegraphics[width=\textwidth]{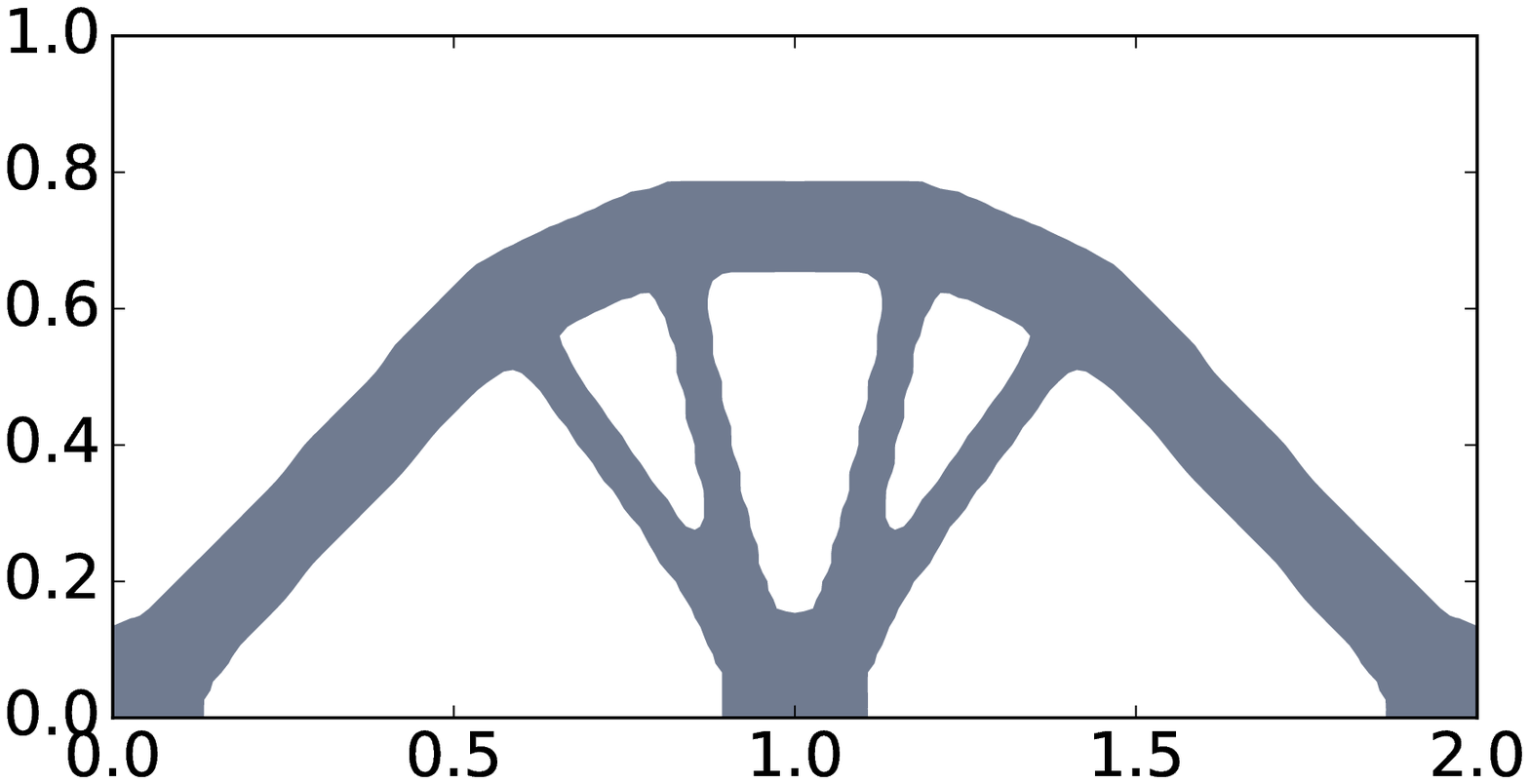}
\caption{$(Nx,Ny) = (150,75)$.}
\label{fig:brdige1b}
\end{subfigure}
\begin{subfigure}[b]{0.3\textwidth}
\includegraphics[width=\textwidth]{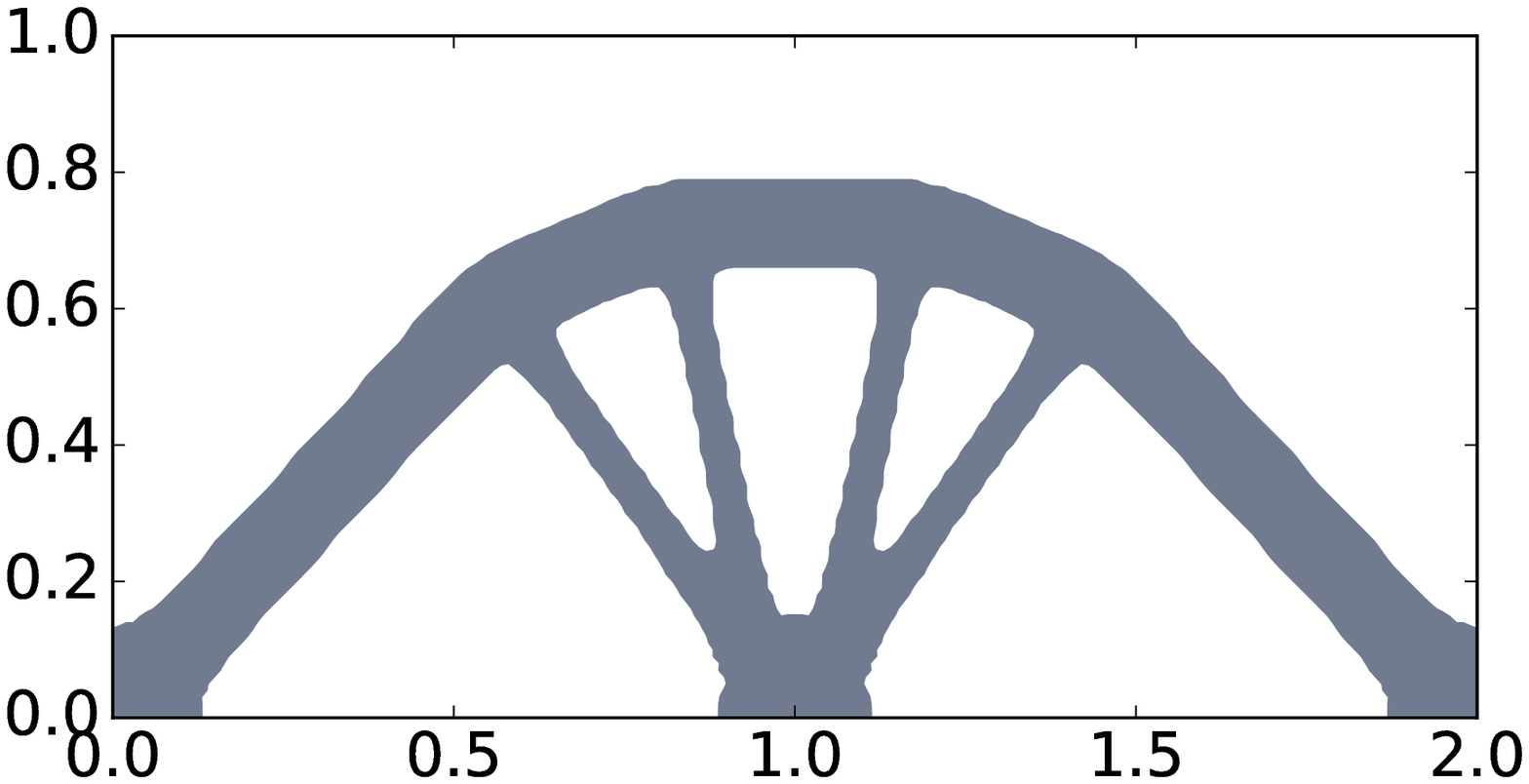}
\caption{$(Nx,Ny) = (200,100)$.}
\label{fig:bridge2b}
\end{subfigure}
\begin{subfigure}[b]{0.3\textwidth}
\includegraphics[width=\textwidth]{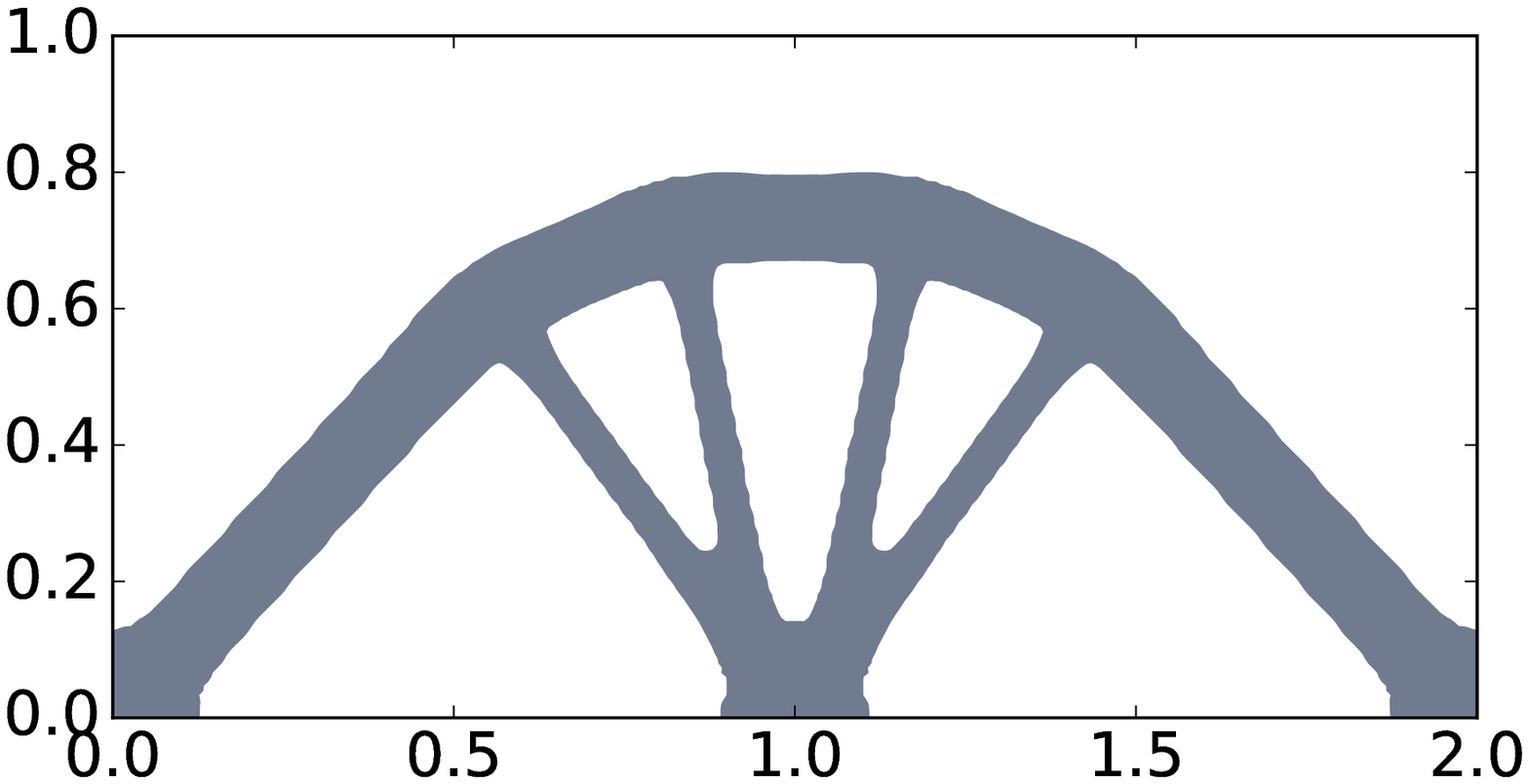}
\caption{$(Nx,Ny) = (300,150)$.}
\label{fig:bridge3b} 
\end{subfigure}
\caption{Optimal design for the bridge,   $\Lambda = 20$ (first row) and $\Lambda = 30$ (second row).}
\label{fig:init_bridges}
\end{center} 
\end{figure}

The case of the bridge is similar to the case of the half-wheel. 
The main difference is the pointwise Dirichlet condition in the lower right corner, which corresponds to
\begin{lstlisting}[numbers=none,xleftmargin=0cm]
DirichletBC(Vvec, (0.0,0.0), dirBd2, method='pointwise')
\end{lstlisting}
Also for the initialization we take 
\begin{lstlisting}[numbers=none, xleftmargin=0cm]
phi_mat = -np.cos((4.0*pi*(XX-1.0))) * np.cos(4*pi*YY) - 0.2 \
+ np.maximum(100.0*(YY-ly+0.05),.0)   
\end{lstlisting}
See Figure \ref{fig:init_bridges} for numerical results for the bridge, with $\Lambda =20$ and $\Lambda =30$. 

\subsection{MBB-beam}\label{subsec:mbb_beam}
We define the MBB beam as in the original paper \cite{Sigmund2014}.
We use the symmetry of the problem to compute the solution only on the right half of the domain. 
Thus we impose rolling boundary condition on the left side of the computational domain $\hold$, which corresponds to $u\cdot n =0$.
We take \f{lx=3.0} and \f{ly=1.0}, and \f{Nx}, \f{Ny} must be chosen accordingly, so as to keep a regular grid. 
For instance, we can choose \f{Nx=150}, \f{Ny=50}.
We take \f{Load = [Point(0.0, 1.0)]}.

\begin{figure}
\begin{center}
\begin{subfigure}[b]{0.3\textwidth}
\includegraphics[width=\textwidth]{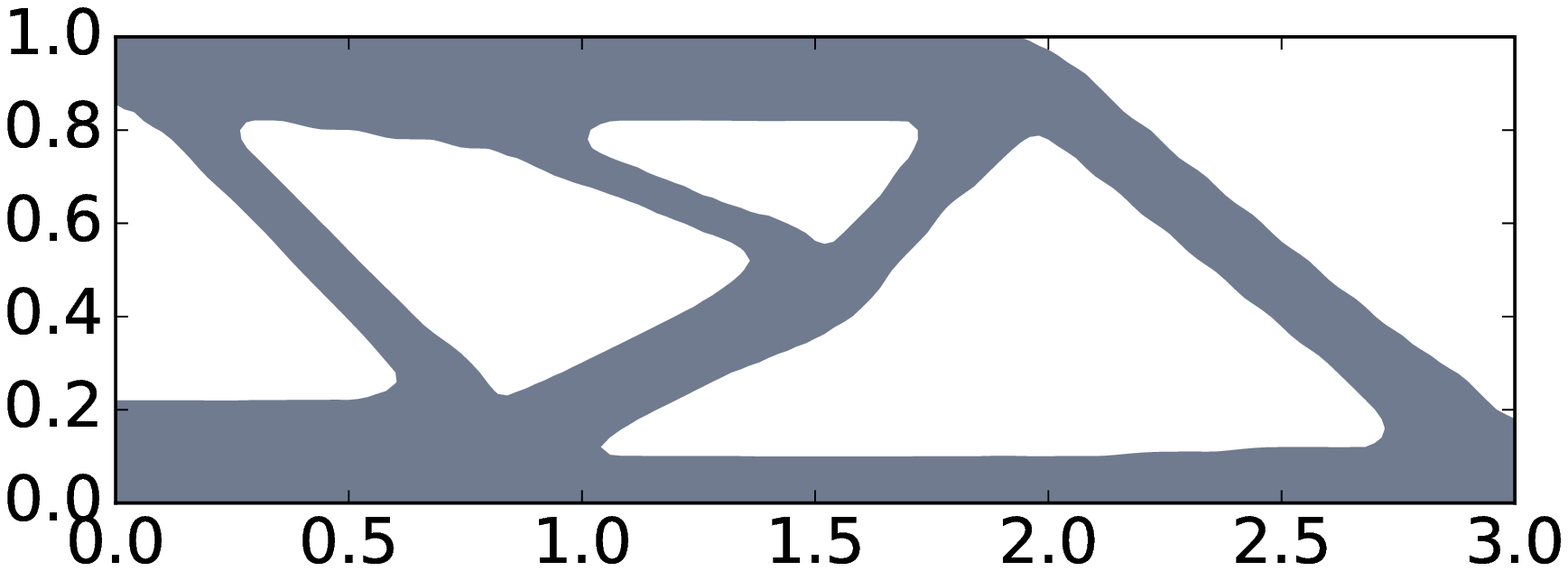}
\caption{$(Nx,Ny) = (150,50)$.}
\label{fig:MBB1}
\end{subfigure}
\begin{subfigure}[b]{0.3\textwidth}
\includegraphics[width=\textwidth]{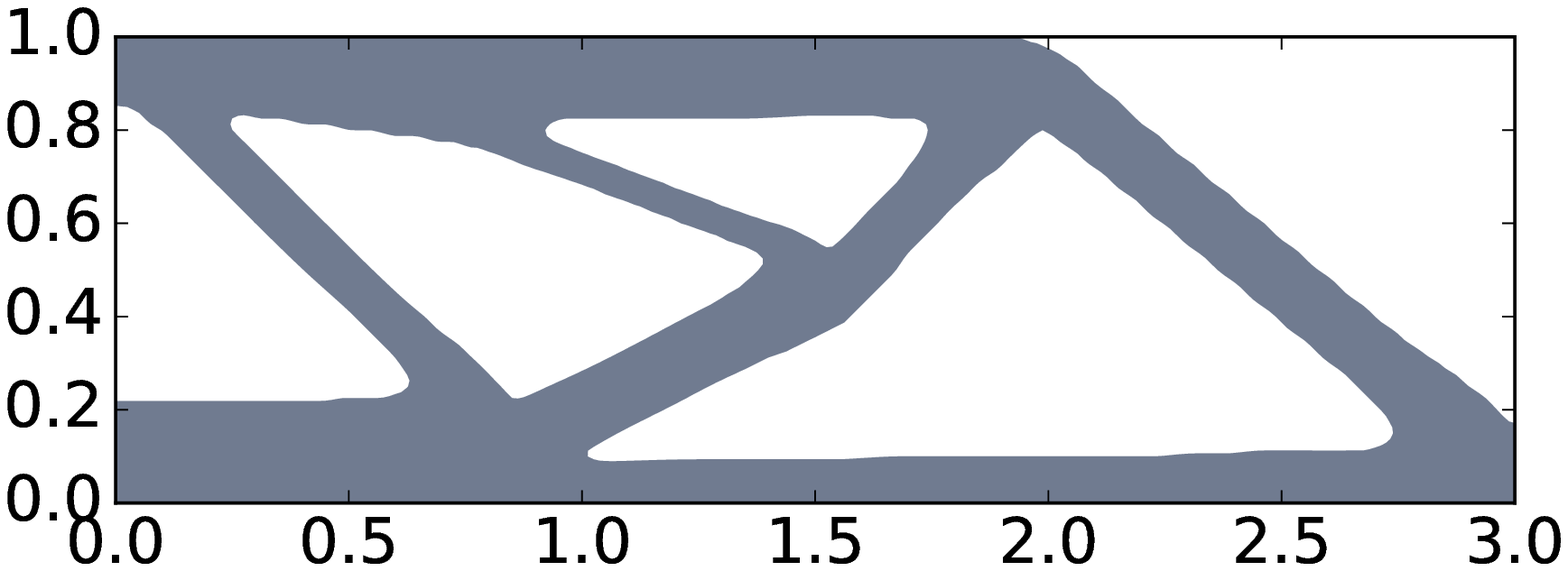}
\caption{$(Nx,Ny) = (240,80)$.}
\label{fig:MBB2}
\end{subfigure}
\begin{subfigure}[b]{0.3\textwidth}
\includegraphics[width=\textwidth]{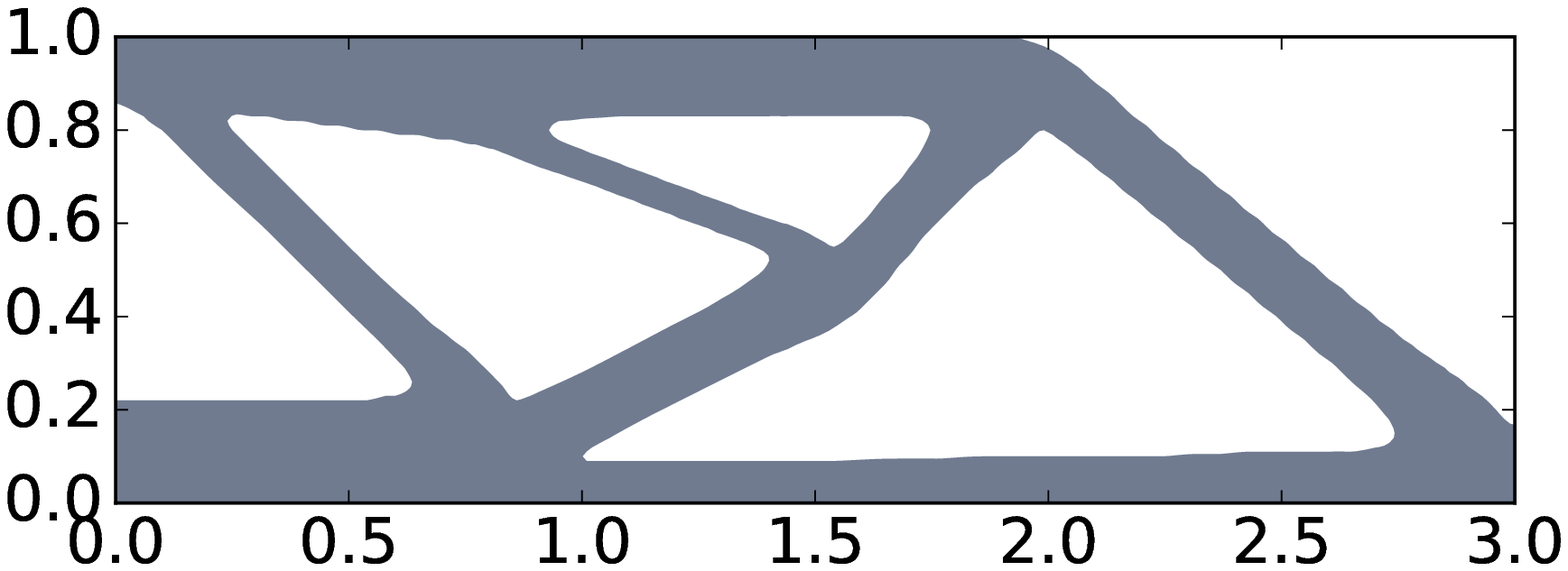}
\caption{$(Nx,Ny) = (300,100)$.}
\label{fig:MBB3}
\end{subfigure}
\caption{Optimal design for the MBB-Beam, $\Lambda = 130$.}
\label{fig:MBB0} 
\end{center} 
\end{figure}

We also have pointwise rolling boundary conditions on the lower right corner of $\hold$.
In \f{init.py} this corresponds to the following definitions of the boundaries:
\begin{lstlisting}[numbers=none,xleftmargin=0cm]
class DirBd(SubDomain):
    def inside(self, x, on_boundary):
        return near(x[0],.0)   
class DirBd2(SubDomain):
    def inside(self, x, on_boundary):
        return abs(x[0]-lx) < tol and abs(x[1])<tol                 
dirBd,dirBd2 = [DirBd(),DirBd2()]     
\end{lstlisting}
Then the boundaries are tagged with different numbers: 
\begin{lstlisting}[numbers=none,xleftmargin=0cm]
dirBd.mark(boundaries, 1)
dirBd2.mark(boundaries, 2)  
\end{lstlisting}
We define the boundary conditions on the two boundaries \f{dirBd} and \f{dirBd2}:
\begin{lstlisting}[numbers=none,xleftmargin=0cm]
bcd=[DirichletBC(Vvec.sub(0),0.0,boundaries,1),\
        DirichletBC(Vvec.sub(1), 0.0, dirBd2,method='pointwise')]   
\end{lstlisting}
Also, the term \f{+ds(2)} in lines 50-51  is active since \f{dirBd2} is not empty, 
as for the half-wheel case.
We also choose an appropriate initialization
\begin{lstlisting}[numbers=none,xleftmargin=0cm]
phi_mat = -np.cos(4.0/lx*pi*XX) *np.cos(4.0*pi*YY)-0.4\
 +np.maximum(100.0*(XX+YY-lx-ly+0.1),.0) +np.minimum(5.0/ly *(YY-1.0) + 4.0,0)  
\end{lstlisting}    
See Figure \ref{fig:MBB0}  for the MBB-beam case with $\Lambda = 130$.

\subsection{Multiple load cases}
We have for this case that 
\begin{lstlisting}[numbers=none, xleftmargin=0cm]
Load = [Point(lx, 0.0),Point(lx, 1.0)] 
\end{lstlisting}
is a list.
In line 30 of \f{compliance.py}, \f{U} is thus a list with two elements corresponding to the two loads.
This explains the \f{for} loop in  \f{_shape_der} (see line 132).

\begin{figure}
\begin{center}
\begin{subfigure}[b]{0.3\textwidth}
\includegraphics[width=\textwidth]{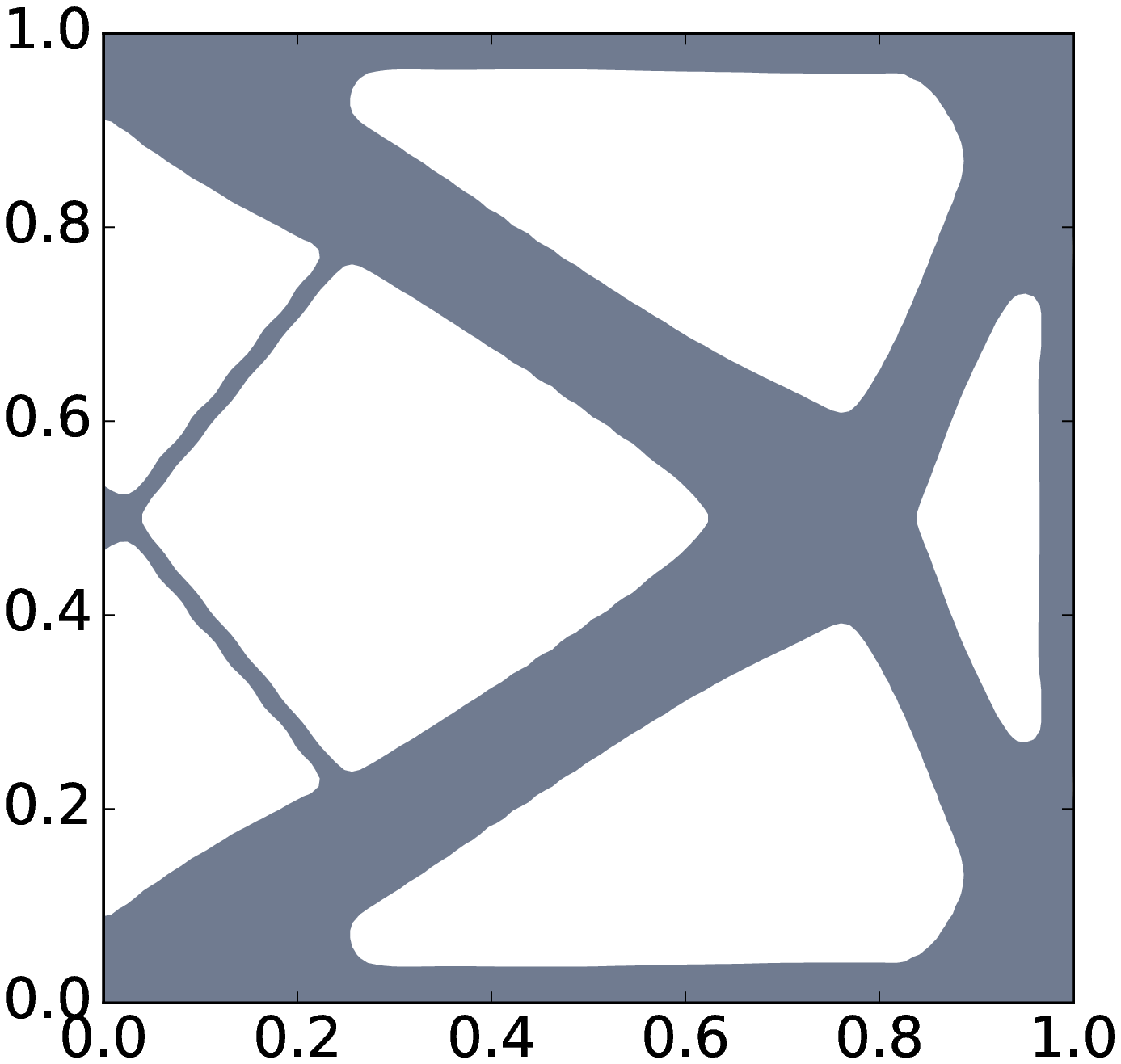}
\caption{$(Nx,Ny) = (121,121)$.}
\label{fig:cantilever_twoforces1}
\end{subfigure}
\begin{subfigure}[b]{0.3\textwidth}
\includegraphics[width=\textwidth]{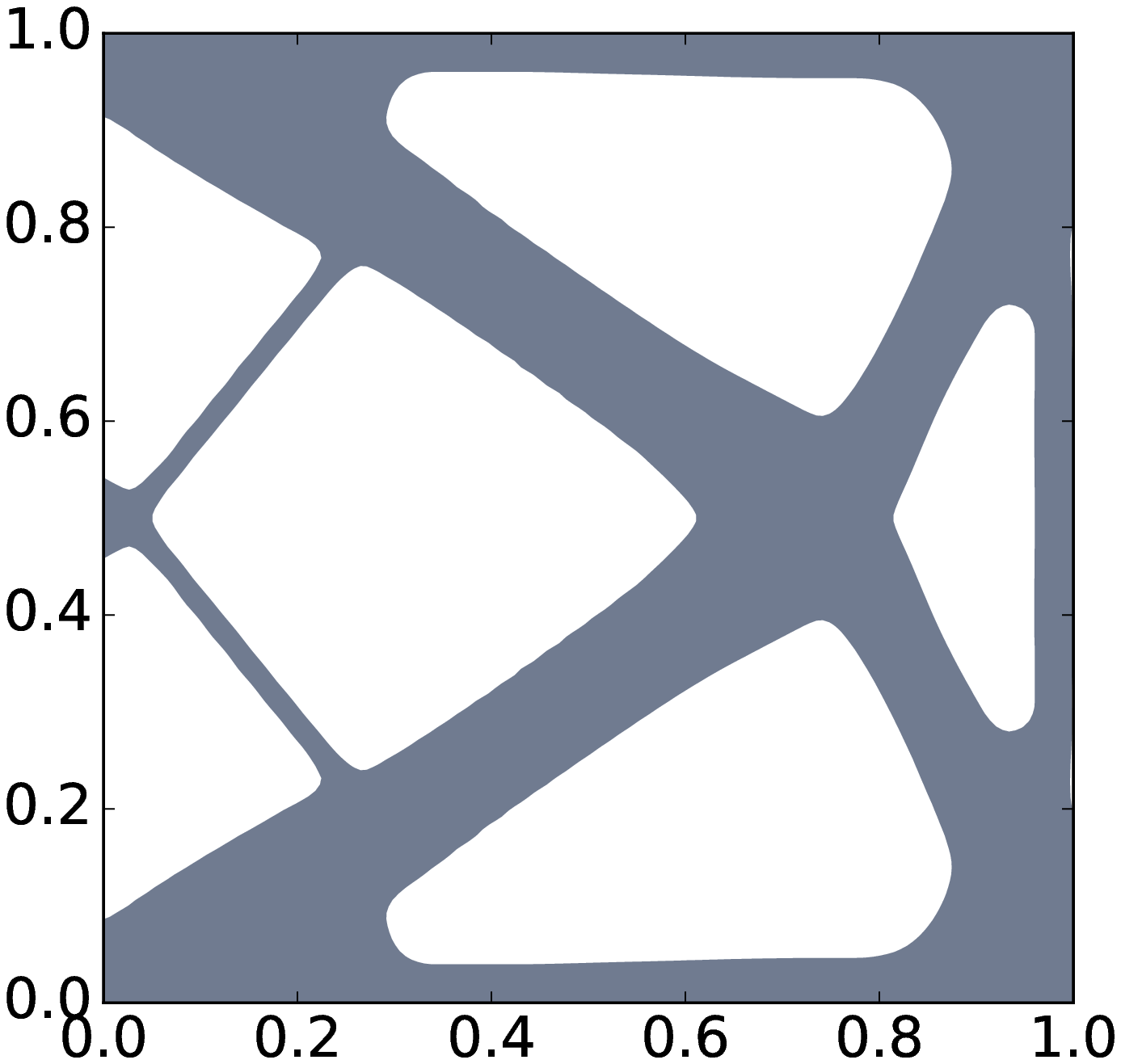}
\caption{$(Nx,Ny) = (151,151)$.}
\label{fig:cantilever_twoforces2}
\end{subfigure}
\begin{subfigure}[b]{0.3\textwidth}
\includegraphics[width=\textwidth]{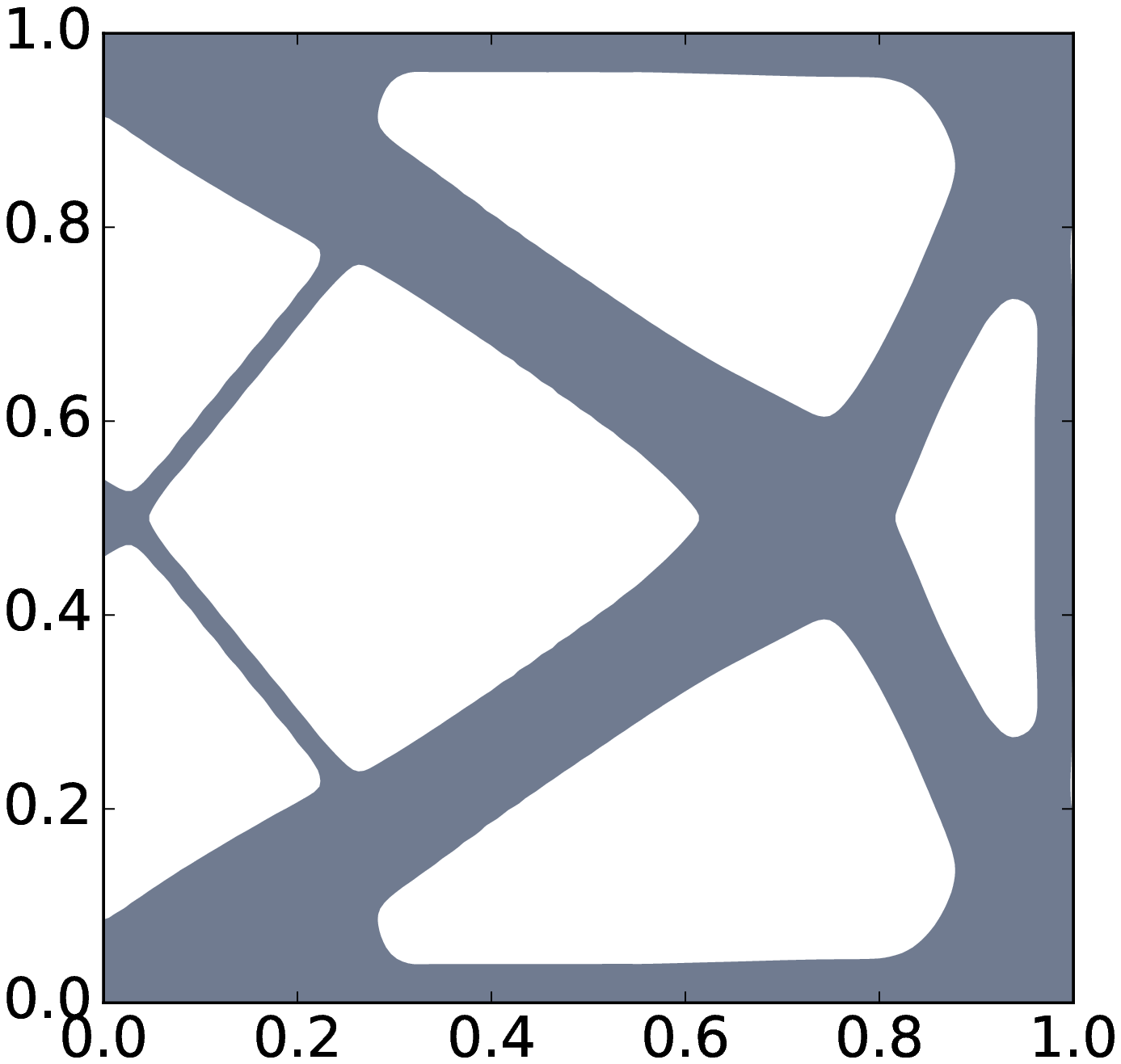}
\caption{$(Nx,Ny) = (175,175)$.}
\label{fig:cantilever_twoforces3}
\end{subfigure}
\caption{Optimal design for the cantilever with two loads, $\Lambda = 60$.}
\label{fig:cantilever_twoforces0} 
\end{center} 
\end{figure}

We illustrate multiple load cases with a cantilever problem with two loads applied at the bottom-right corner and the top-right corner, both with equal intensities to get a symmetric design. 
Here the Lagrangian is taken as $\Lambda = 60$. The results are shown in Figure \ref{fig:cantilever_twoforces0}.
We use the initialization
\begin{lstlisting}[numbers=none]
phi_mat = -np.cos(4.0*pi*(XX-0.5)) * np.cos(4.0*pi*(YY-0.5)) - 0.6 \
          -np.maximum(50.0*(YY-ly+0.1),.0)- np.maximum(50.0*(-YY+0.1),.0)    
\end{lstlisting}

\subsection{Inverter}

Mechanisms require additional modifications of the code, therefore we discuss here briefly the main differences and provide the code for the inverter separately as the file \f{mechanism.py}. 
The code can be downloaded at  \burl{http://antoinelaurain.com/compliance.htm}.
To run the inverter, type
\begin{lstlisting}[numbers=none,xleftmargin=0cm]
python mechanism.py inverter
\end{lstlisting}
Unlike the compliance, the case of compliance mechanisms is not self-adjoint, therefore we need to compute an adjoint given by \eqref{adjoint_inverter}.
For this we add a subfunction \f{_solve_adj} to compute the adjoint. 
The function \f{_solve_adj} works like \f{_solve_pde}, but implements the right-hand side corresponding to \eqref{adjoint_inverter}.
The modification of the objective functional and of \f{_shape_der} follows the description of Section \ref{sec:inverter} in a straightforward way.
In the \f{init.py} file, we define  the boundaries  \f{outputBd} and \f{inputBd} which correspond to $\Gamma_{out}$ and $\Gamma_{in}$ of Section \ref{sec:inverter}, respectively.

We take $\Gamma_{in} = \{0\}\times (0.47,0.53)$ and $\Gamma_{out} = \{1\}\times (0.43,0.57)$. 
In order to keep the regions around $\Gamma_{out}$ and $\Gamma_{in}$ fixed, we  define and tag the following small region in \f{mechanism.py}
\begin{lstlisting}[numbers=none,xleftmargin=0cm]
class Fixed(SubDomain):
    def inside(self, x, on_boundary):
        return (between(x[0], (.0,.05)) and between(x[1], (.48,.52)))\
         or  (between(x[0], (.9,1.0)) and between(x[1], (.43,.57)) )  
fixed = Fixed()  
fixed.mark(domains, 2) 
\end{lstlisting}
This is used in the following definition
\begin{lstlisting}[numbers=none,xleftmargin=0cm]
av = assemble((inner(grad(theta),grad(xi)) +0.1*inner(theta,xi))*dx(0)\
     +1.0e5*inner(theta,xi) * dx(2)\
     +1.0e5*(inner(dot(theta,n),dot(xi,n)) * (ds(0)+ds(1)+ds(2)+ds(3))) )
\end{lstlisting}
The large coefficient \f{1.0e5} in the subdomain \f{fixed} forces \f{th} to  be close to zero during the entire process.

We add a volume term to the objective functional with the coefficient $\Lambda = 0.01$. 
We choose the parameters $k_s = 0.01$,  $\epsilon =0.01$, $E = 20$, $\eta_{in}=2,\eta_{out}=1$,   \f{lx=1.0}, \f{ly=1.0}, \f{beta0_init = 1.0}, \f{ItMax = int(2.0*Nx)} and 
\f{delta = PointSource(V.sub(0), Load, 0.05)}
in function \f{_solve_pde}.
The other parameters are the same as in \f{compliance.py}. 
For the initialization of \f{phi_mat} we refer to the file \f{init.py}.
See Figure \ref{fig:mechanism_fig} and Section \ref{sec:inverter} for a description of the design domain, boundary conditions and optimal design.
\section{Initialization and computation time}\label{sec:init_time}
\subsection{Influence of initialization}
\begin{figure}
\begin{center}
\begin{subfigure}[b]{0.3\textwidth}
\includegraphics[width=\textwidth]{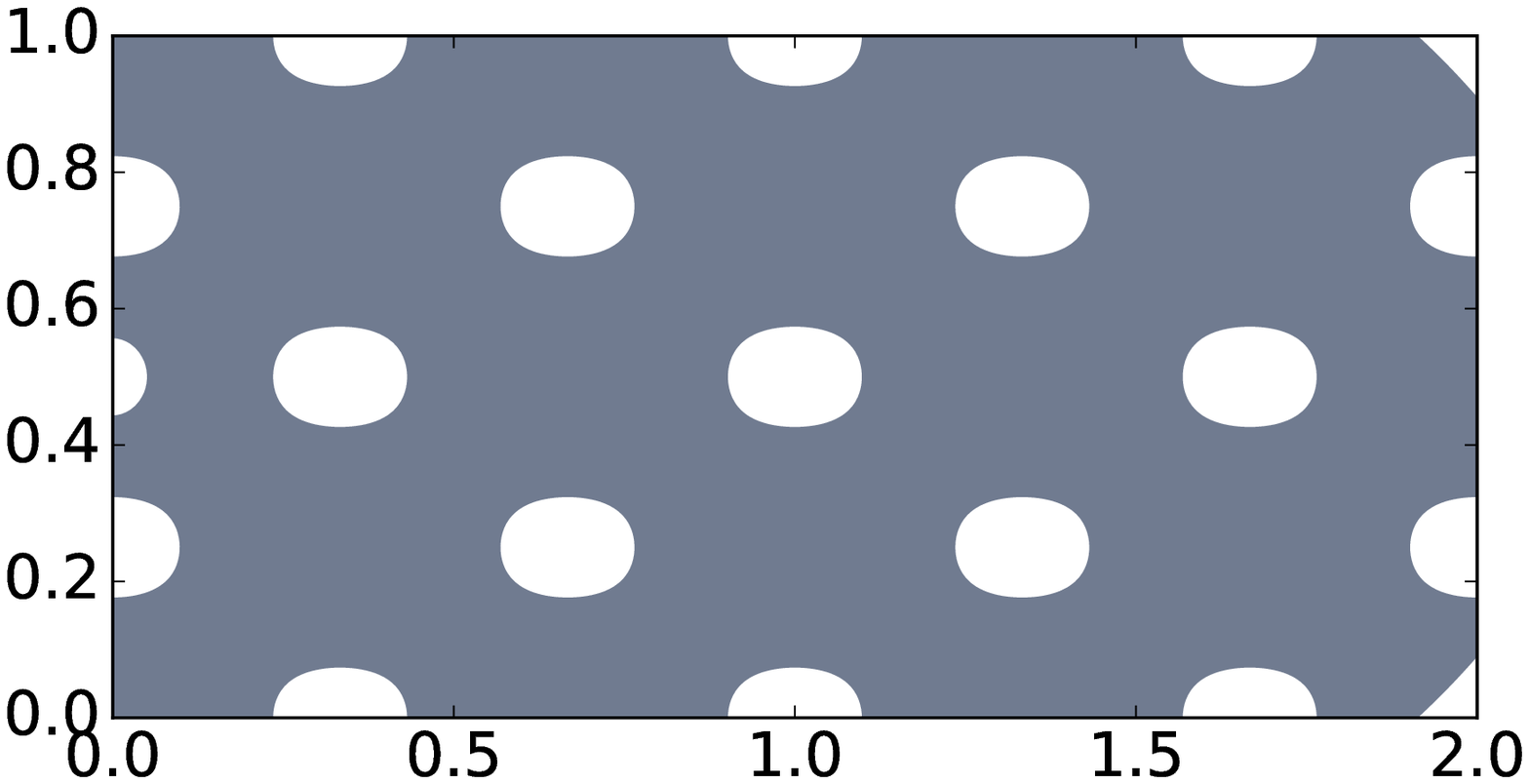}
\label{fig:cantilever1_init1}
\end{subfigure}
\begin{subfigure}[b]{0.3\textwidth}
\includegraphics[width=\textwidth]{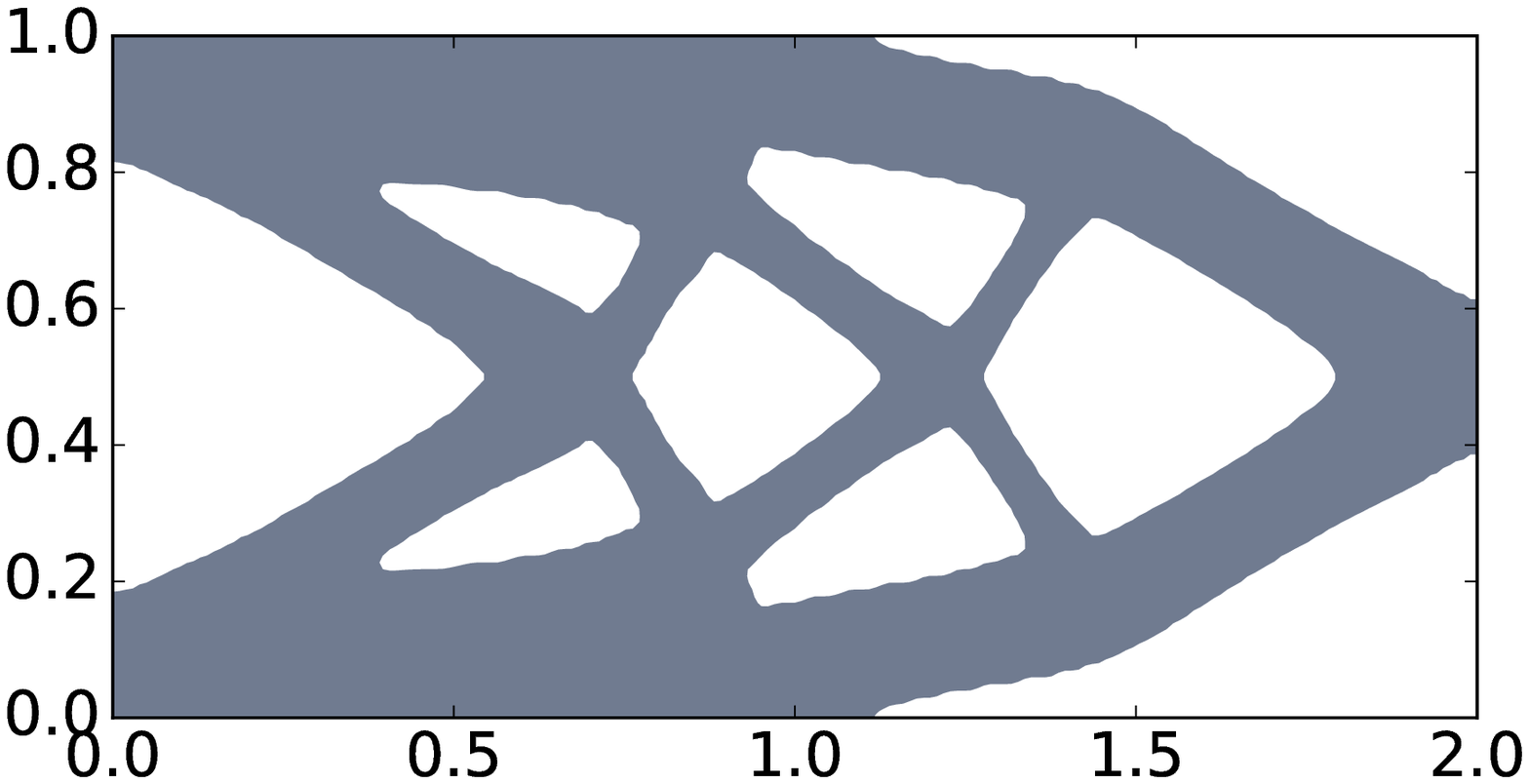}
\label{fig:cantilever2_init1}
\end{subfigure}
\begin{subfigure}[b]{0.3\textwidth}
\includegraphics[width=\textwidth]{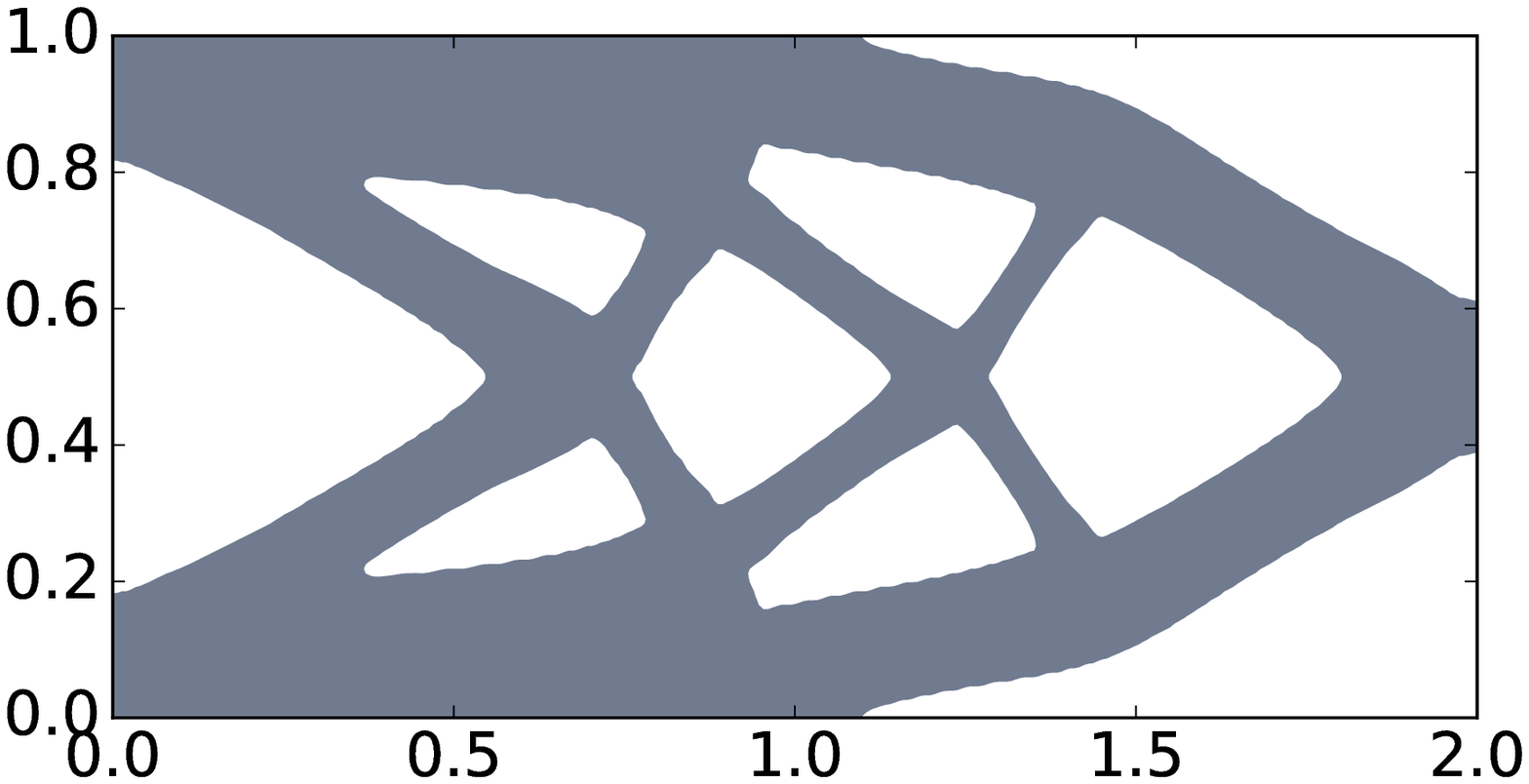}
\label{fig:cantilever3_init1}
\end{subfigure}
\begin{subfigure}[b]{0.3\textwidth}
\includegraphics[width=\textwidth]{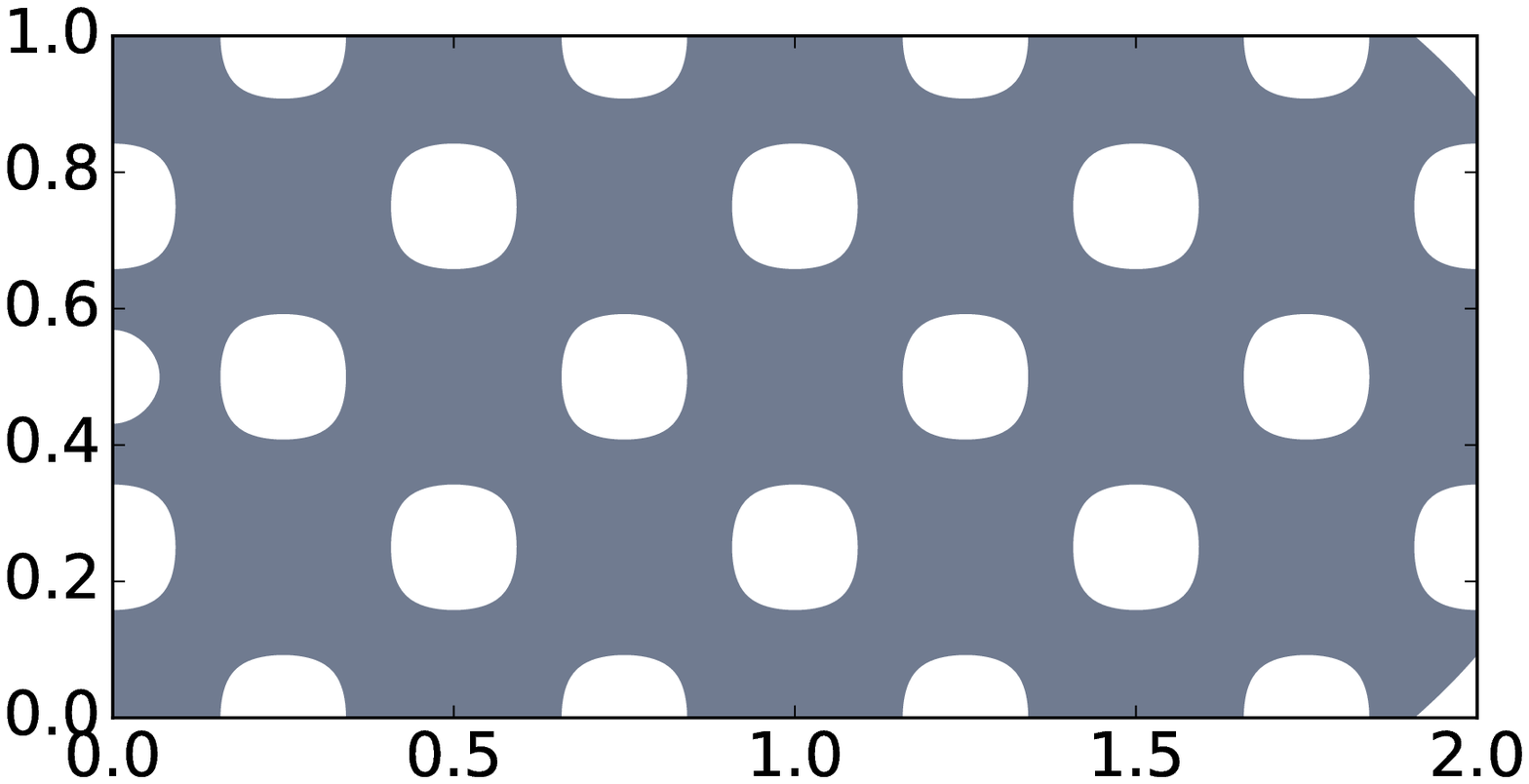}
\label{fig:cantilever1_init2}
\end{subfigure}
\begin{subfigure}[b]{0.3\textwidth}
\includegraphics[width=\textwidth]{cantilever_init2_lagV=40_Nx=202_final_design.eps}
\label{fig:cantilever2_init2}
\end{subfigure}
\begin{subfigure}[b]{0.3\textwidth}
\includegraphics[width=\textwidth]{cantilever_init2_lagV=40_Nx=302_final_design.eps}
\label{fig:cantilever3_init2}
\end{subfigure}
\caption{Optimal design for the symmetric cantilever,  with $\Lambda = 40$ and two different initial guesses. First column: initial guess, second column: optimal design for $(Nx,Ny) = (202,101)$, third column: optimal design for $(Nx,Ny) = (302,151)$. The first line uses initialization \eqref{init_phi2}, the second line uses initialization \eqref{init_phi}.}
\label{fig:cantilever0} 
\end{center} 
\end{figure}

It is known that the final result may depend on the initial guess for the minimization of the compliance.
We observe this phenomenon in our algorithm, as illustrated in Figure \ref{fig:cantilever0}, where two different initializations provide two different optimal designs. 
We compare initialization \eqref{init_phi} with
\begin{align} 
\label{init_phi2}
\begin{split}
\phi(x,y) =& -\cos(6\pi x/l_x)\cos(4\pi y) - 0.6
+ \max(200(0.01 - x^2 - (y-l_y/2)^2),0)\\
& + \max(100(x +y -l_x-l_y+0.1),0)
+ \max(100(x -y -l_x+0.1),0).
\end{split}
\end{align}
The choice $\phi(x,y) = -\cos(6\pi x/l_x)\cos(4\pi y) - 0.6$ corresponds to a standard choice of seven holes inside the domain for the $2\times 1$ cantilever; see \cite{MR2414920}. 
It can be seen in Figure \ref{fig:cantilever0} that the initialization with the higher number of holes provides an optimal design $\hold\setminus\Om$ with more connected components.

\subsection{Computation time}

The numerical tests were run on a PC with four processors Intel Core2 Q9400, 2.66 GHz, 3.8 GB memory, with LinuxMint 17 and FEniCS 2017.1. 
In Tables \ref{tab4} and \ref{tab5} we show the computation time for the symmetric and asymmetric cantilevers.
The average time for one iteration is computed by averaging over all iterations. 
However, it is not counting the time spent by \f{init.py}, and the time spent when the line search is performed, i.e. the time is recorded only for steps which are accepted.

Since we use a mesh with crossed elements, the number of elements is 
\f{dofsV_max = (Nx+1)*(Ny+1) + Nx*Ny}, see line 37.  
We solve two partial differential equations during each iteration (again, without counting the line search, and assuming we have only one load), one to compute \f{U}, and one to compute \f{th}.
Since these are vectors, the number of degrees of freedom for solving each of these PDEs is
\begin{lstlisting}[numbers=none,xleftmargin=0cm]
2*dofsV_max = 2*((Nx+1)*(Ny+1) + Nx*Ny) 
\end{lstlisting}
For instance for the case \f{(Nx,Ny)=(302,151)}, as in Table \ref{tab4}, we get $91,658$ elements and  $183,316$ degrees of freedoms.
When comparing with the computational time for an algorithm such as the one described in \cite{Andreassen2010}, one should use the number of elements as the basis for comparison. 
For example, a $300\times 100$ mesh in \cite{Andreassen2010} gives 30,000 elements, corresponding approximately to a grid of $170\times 85$,  which gives $29,156$ elements for our code. 

The computation time is comparable with the results in \cite{Andreassen2010}, although slightly slower, for the same number of elements. 
Comparing with the educational code from \cite{Challis2009}, which is also based on the level set method, our code is significantly faster. 
Indeed, it was observed in \cite{Gain2013} that the code of \cite{Challis2009} takes a long time to converge if the mesh discretization is greater than $5,000$ elements. 
In \cite{Gain2013} the authors have improved its efficiency by using a sparse matrix assembly, but they did not report on computation time.

\begin{table}
\caption{Computation time for the asymmetric cantilever benchmark for $\Lambda = 60$. }
\begin{center}
\begin{tabular}{@{}lccc@{}}
\toprule
Mesh size                      & $102\times 51$   & $202\times 101$ &  $302\times 151$\\
\toprule
Number of  elements                         & 10,558  & 41,108 & 91,658  \\ \midrule
\lstinline[]$ItMax$                & 153 &  303 & 453 \\ \midrule
Total iterations                 & 72  & 303 & 377  \\ \midrule
Average time per iteration (s)   & 2.02  & 7.90 & 18.00  \\ \midrule
Total time (h:m:s)               & 0:04:27 &  1:07:13 & 2:12:54 \\ \bottomrule
\end{tabular}
\end{center}
\label{tab4}
\end{table}

\begin{table}
\caption{Computation time for the symmetric cantilever benchmark for $\Lambda = 40$ and initialization \eqref{init_phi2}. }
\begin{center}
\begin{tabular}{@{}lccc@{}}
\toprule
Mesh size                                & $102\times 51$   & $202\times 101$ &  $302\times 151$ \\
\toprule
Number of elements               & 10,558  & 41,108 & 91,658  \\ \midrule
\lstinline[]$ItMax$                     & 153  & 303 & 453 \\ \midrule
Total iterations      & 60  & 82 & 247  \\ \midrule
Average time per iteration (s)   & 2.00  & 8.25 & 18.20 \\ \midrule
Total time (h:m:s)     & 0:04:56 &  0:18:00 & 2:17:37  \\ \bottomrule
\end{tabular}
\end{center}
\label{tab5}
\end{table}

\section{Conclusion}
We have presented a FEniCS code for structural optimization based on the level set method. 
The principal feature of the code is to rely on the notion of distributed shape derivative, which is easy to implement with FEniCS, and on the corresponding reformulation of the level set equation. 
We have shown how to compute the distributed shape derivative for a fairly general functional which can be used for compliance minimization and compliant mechanisms in particular. 
Various benchmarks of compliance minimization were tested, as well as an example of inverter mechanism.

We encourage students and newcomers to the field to experiment with new examples and parameters. 
The code can be used as a basis for more advanced problems. 
One could take advantage of the versatility of FEniCS to solve various types of PDEs, and adapt the code for multiphysics problems. 
An extension to three dimensions is also relatively easy using FEniCS, since the variational formulation is independent on the dimension. 
The main effort for extending the present code to three dimensions resides in adapting the numerical scheme for the level set part.\\

\noindent {\bf Acknowledgements.} The author acknowledges the support of the Brazilian National Council for Scientific and Technological Development  (Conselho Nacional de Desenvolvimento Cient\'ifico e Tecnol\'ogico - CNPq), through the program  ``Bolsa de Produtividade em Pesquisa - PQ 2015'', process: 302493/2015-8.
The author also acknowledges the support of FAPESP (Funda\c{c}\~ao de Amparo \`a Pesquisa do Estado de S\~ao Paulo), process: 2016/24776-6.

\section{Appendix: FEniCS code {\tt compliance.py}}

\begin{lstlisting}[frame=none,
  language=Python,
  aboveskip=3mm,
  belowskip=3mm,
  showstringspaces=false,
  columns=flexible,
  basicstyle={\small\ttfamily},
  numbers=left,
  numberstyle=\tiny\color{gray},
  keywordstyle=\color{blue},
  commentstyle=\color{dkgreen},
  stringstyle=\color{mauve},
  breaklines=true,
  breakatwhitespace=true,
  tabsize=3]
# ----------------------------------------------------------------------
# FEniCS 2017.1 code for level set-based structural optimization.
# Written by Antoine Laurain, 2017
# ----------------------------------------------------------------------
from dolfin import *
from init import *
from matplotlib import cm,pyplot as pp
import numpy as np, sys, os
pp.switch_backend('Agg')
set_log_level(ERROR)
# ----------------------------------------------------------------------
def _main():
    Lag,Nx,Ny,lx,ly,Load,Name,ds,bcd,mesh,phi_mat,Vvec=init(sys.argv[1])      
    eps_er, E, nu = [0.001, 1.0, 0.3]  # Elasticity parameters
    mu,lmbda = Constant(E/(2*(1 + nu))),Constant(E*nu/((1+nu)*(1-2*nu)))
    # Create folder for saving files
    rd = os.path.join(os.path.dirname(__file__),\
     Name +'/LagVol=' +str(np.int_(Lag))+'_Nx='+str(Nx))
    if not os.path.isdir(rd): os.makedirs(rd) 
    # Line search parameters
    beta0_init,ls,ls_max,gamma,gamma2 = [0.5,0,3,0.8,0.8]   
    beta0 = beta0_init
    beta  = beta0
    # Stopping criterion parameters    
    ItMax,It,stop = [int(1.5*Nx), 0, False] 
    # Cost functional and function space
    J = np.zeros( ItMax )  
    V = FunctionSpace(mesh, 'CG', 1)
    VolUnit = project(Expression('1.0',degree=2),V) # to compute volume 
    U = [0]*len(Load) # initialize U
    # Get vertices coordinates 
    gdim     = mesh.geometry().dim()
    dofsV    = V.tabulate_dof_coordinates().reshape((-1, gdim))    
    dofsVvec = Vvec.tabulate_dof_coordinates().reshape((-1, gdim))     
    px,py    = [(dofsV[:,0]/lx)*2*Nx, (dofsV[:,1]/ly)*2*Ny]
    pxvec,pyvec = [(dofsVvec[:,0]/lx)*2*Nx, (dofsVvec[:,1]/ly)*2*Ny]  
    dofsV_max, dofsVvec_max =((Nx+1)*(Ny+1) + Nx*Ny)*np.array([1,2])    
    # Initialize phi  
    phi = Function( V )  
    phi = _comp_lsf(px,py,phi,phi_mat,dofsV_max) 
    # Define Omega = {phi<0}     
    class Omega(SubDomain):
        def inside(self, x, on_boundary):
            return .0 <= x[0] <= lx and .0 <= x[1] <= ly and phi(x) < 0             
    domains = CellFunction("size_t", mesh)     
    dX = Measure('dx') 
    n  = FacetNormal(mesh) 
    # Define solver to compute descent direction th
    theta,xi = [TrialFunction(Vvec), TestFunction( Vvec)]     
    av = assemble((inner(grad(theta),grad(xi)) +0.1*inner(theta,xi))*dX\
         + 1.0e4*(inner(dot(theta,n),dot(xi,n)) * (ds(0)+ds(1)+ds(2))) ) 
    solverav = LUSolver(av)
    solverav.parameters['reuse_factorization'] = True                  
    #---------- MAIN LOOP ----------------------------------------------
    while It < ItMax and stop == False:
        # Update and tag Omega = {phi<0}, then solve elasticity system.  
        omega = Omega()
        domains.set_all(0)
        omega.mark(domains, 1)
        dx = Measure('dx')(subdomain_data = domains)   
        for k in range(0,len(Load)):   
            U[k] = _solve_pde(Vvec,dx,ds,eps_er,bcd,mu,lmbda,Load[k])      
        # Update cost functional 
        compliance = 0
        for u in U:
            eU = sym(grad(u))
            S1 = 2.0*mu*inner(eU,eU) + lmbda*tr(eU)**2
            compliance += assemble( eps_er*S1* dx(0) + S1*dx(1) )  
        vol = assemble( VolUnit*dx(1) )          
        J[It]  = compliance + Lag * vol                            
        # ------- LINE SEARCH ------------------------------------------
        if It > 0 and J[It] > J[It-1] and ls < ls_max:
            ls   += 1
            beta *= gamma            
            phi_mat,phi = [phi_mat_old,phi_old]
            phi_mat = _hj(th_mat, phi_mat, lx,ly,Nx, Ny, beta)
            phi     = _comp_lsf(px,py,phi,phi_mat,dofsV_max) 
            print('Line search iteration : %s' % ls)           
        else:          
            print('************ ITERATION NUMBER %s' % It)                   
            print('Function value        : %.2f' % J[It])
            print('Compliance            : %.2f' % compliance)
            print('Volume fraction       : %.2f' % (vol/(lx*ly))) 
            # Decrease or increase line search step
            if ls == ls_max: beta0 = max(beta0 * gamma2, 0.1*beta0_init)  
            if ls == 0:      beta0 = min(beta0 / gamma2, 1)
            # Reset beta and line search index    
            ls,beta,It = [0,beta0, It+1]
            # Compute the descent direction th           
            th = _shape_der(Vvec,U,eps_er,mu,lmbda,dx,solverav,Lag)
            th_array = th.vector().array()
            th_mat = [np.zeros((Ny+1,Nx+1)),np.zeros((Ny+1,Nx+1))]          
            for dof in xrange(0, dofsVvec_max,2):
                if np.rint(pxvec[dof]) %2 == .0:
                    cx,cy= np.int_(np.rint([pxvec[dof]/2,pyvec[dof]/2]))
                    th_mat[0][cy,cx] = th_array[dof]
                    th_mat[1][cy,cx] = th_array[dof+1]                           
            # Update level set function phi using descent direction th
            phi_old, phi_mat_old = [phi, phi_mat]            
            phi_mat = _hj(th_mat, phi_mat, lx,ly,Nx,Ny, beta)
            if np.mod(It,5) == 0: phi_mat = _reinit(lx,ly,Nx,Ny,phi_mat)     
            phi = _comp_lsf(px,py,phi,phi_mat,dofsV_max)                      
            #------------ STOPPING CRITERION ---------------------------
            if It>20 and max(abs(J[It-5:It]-J[It-1]))<2.0*J[It-1]/Nx**2: 
                stop = True  
            #------------ Plot Geometry --------------------------------  
            if np.mod(It,10)==0 or It==1 or It==ItMax or stop==True:   
                pp.close()     
                pp.contourf(phi_mat,[-10.0,.0],extent = [.0,lx,.0,ly],\
                 cmap=cm.get_cmap('bone'))
                pp.axes().set_aspect('equal','box')
                pp.show()
                pp.savefig(rd+'/it_'+str(It)+'.pdf',bbox_inches='tight')             
    return
# ----------------------------------------------------------------------
def _solve_pde(V, dx, ds, eps_er, bcd, mu, lmbda, Load):
    u,v = [TrialFunction(V), TestFunction(V)]
    S1 = 2.0*mu*inner(sym(grad(u)),sym(grad(v))) + lmbda*div(u)*div(v)
    A = assemble( S1*eps_er*dx(0) + S1*dx(1) )   
    b = assemble( inner(Expression(('0.0', '0.0'),degree=2) ,v) * ds(2))    
    U = Function(V)
    delta = PointSource(V.sub(1), Load, -1.0)
    delta.apply(b) 
    for bc in bcd: bc.apply(A,b)    
    solver = LUSolver(A)
    solver.solve(U.vector(), b)    
    return U
#-----------------------------------------------------------------------    
def _shape_der(Vvec, u_vec , eps_er, mu, lmbda, dx, solver, Lag):   
    xi = TestFunction(Vvec)  
    rv = 0.0 
    for u in u_vec:       
        eu,Du,Dxi = [sym(grad(u)),grad(u),grad(xi)]
        S1 = 2*mu*(2*inner((Du.T)*eu,Dxi) -inner(eu,eu)*div(xi))\
         + lmbda*(2*inner( Du.T, Dxi )*div(u) - div(u)*div(u)*div(xi) )
        rv += -assemble(eps_er*S1*dx(0) + S1*dx(1) + Lag*div(xi)*dx(1))
    th = Function(Vvec)                  
    solver.solve(th.vector(), rv)
    return th
#-----------------------------------------------------------------------        
def _hj(v,psi,lx,ly,Nx,Ny,beta): 
    for k in range(10):
        Dym = Ny*np.repeat(np.diff(psi,axis=0),[2]+[1]*(Ny-1),axis=0)/ly 
        Dyp = Ny*np.repeat(np.diff(psi,axis=0),[1]*(Ny-1)+[2],axis=0)/ly
        Dxm = Nx*np.repeat(np.diff(psi),[2]+[1]*(Nx-1),axis=1)/lx 
        Dxp = Nx*np.repeat(np.diff(psi),[1]*(Nx-1)+[2],axis=1)/lx          
        g = 0.5*( v[0]*(Dxp + Dxm) + v[1]*(Dyp + Dym)) \
          - 0.5*(np.abs(v[0])*(Dxp - Dxm) + np.abs(v[1])*(Dyp - Dym)) 
        maxv = np.max(abs(v[0]) + abs(v[1]))
        dt  = beta*lx / (Nx*maxv)
        psi = psi - dt*g
    return  psi         
#-----------------------------------------------------------------------         
def _reinit(lx,ly,Nx,Ny,psi):      
    Dxs = Nx*(np.repeat(np.diff(psi),[2]+[1]*(Nx-1),axis=1) \
          +np.repeat(np.diff(psi),[1]*(Nx-1)+[2],axis=1))/(2*lx) 
    Dys = Ny*(np.repeat(np.diff(psi,axis=0),[2]+[1]*(Ny-1),axis=0)\
          +np.repeat(np.diff(psi,axis=0),[1]*(Ny-1)+[2],axis=0))/(2*ly)                  
    signum = psi / np.power(psi**2 + ((lx/Nx)**2)*(Dxs**2+Dys**2),0.5)
    for k in range(0,2):      
        Dym = Ny*np.repeat(np.diff(psi,axis=0),[2]+[1]*(Ny-1),axis=0)/ly 
        Dyp = Ny*np.repeat(np.diff(psi,axis=0),[1]*(Ny-1)+[2],axis=0)/ly
        Dxm = Nx*np.repeat(np.diff(psi),[2]+[1]*(Nx-1),axis=1)/lx 
        Dxp = Nx*np.repeat(np.diff(psi),[1]*(Nx-1)+[2],axis=1)/lx               
        Kp = np.sqrt((np.maximum(Dxm,0))**2 + (np.minimum(Dxp,0))**2 \
           + (np.maximum(Dym,0))**2 + (np.minimum(Dyp,0))**2)
        Km = np.sqrt((np.minimum(Dxm,0))**2 + (np.maximum(Dxp,0))**2 \
           + (np.minimum(Dym,0))**2 + (np.maximum(Dyp,0))**2)           
        g  = np.maximum(signum,0)*Kp + np.minimum(signum,0)*Km
        psi  = psi - (0.5*lx/Nx)*(g - signum)       
    return psi   
#-----------------------------------------------------------------------
def _comp_lsf(px,py,phi,phi_mat,dofsV_max):               
    for dof in range(0,dofsV_max):              
        if np.rint(px[dof]) %2 == .0: 
            cx,cy = np.int_(np.rint([px[dof]/2,py[dof]/2]))                                            
            phi.vector()[dof] = phi_mat[cy,cx]
        else:
            cx,cy = np.int_(np.floor([px[dof]/2,py[dof]/2]))                      
            phi.vector()[dof] = 0.25*(phi_mat[cy,cx] + phi_mat[cy+1,cx]\
              + phi_mat[cy,cx+1] + phi_mat[cy+1,cx+1])    
    return phi            
# ----------------------------------------------------------------------
if __name__ == '__main__':
    _main()  
\end{lstlisting}

\bibliographystyle{plain}
\bibliography{struct_SMO.bib}
\end{document}